\documentclass[12pt]{amsart}
\usepackage[latin1]{inputenc}

\usepackage{xspace}
\usepackage{amsfonts,amssymb}
\usepackage{amsmath}
\usepackage{amsthm}
\usepackage{tikz}
\usepackage[latin1]{inputenc}
\usepackage[cyr]{aeguill}
\usepackage[english]{babel}
\usepackage{graphics}
\usepackage{vmargin}
\usepackage{xspace}
\usepackage{amsmath,amsfonts}

\def\footnoterule{\kern -15pt \hrule width 2truein \kern10.6pt\relax}%

\newtheorem{thm}{Theorem}[section]

\newtheorem{claim}{Claim}

\newtheorem{prop}[thm]{Proposition}
\newtheorem{df}[thm]{Definition}
\newtheorem{lm}[thm]{Lemma}
\newtheorem{cor}[thm]{Corollary}
\newtheorem{rmq}[thm]{Remark}
\newtheorem{prop-def}[thm]{Proposition and Definition}
\newtheorem{cj}{Conjecture}[section]
\newtheorem{cjapart}{Conjecture}

\def\rev{(M_i\rightarrow M)_{i\in \nn}}
\def\disque{\mathbb{D}^2}
\def\sphere{\mathbb{S}^2}

\def\cercle{\mathbb{S}^1}
\def\H3{\mathbb{H}^3}
\def\zz{\mathbb{Z}}
\def\z2z{\mathbb{Z}/_{2\mathbb{Z}}}

\def\nn{\mathbb{N}}

\def\rr{\mathbb{R}}

\def\bb{\mathbb{B}}

\def\vol{{\rm Vol}}

\def\inj{{\rm Inj}}
\def\dist{{\rm dist}}
\def\d{{\rm d}}
\def\diam{{\rm diam}}
\def\lgr{{\rm length}\,}

\def\div{\rightarrow+\infty}

\newcommand{\incl}[1][r]
  {\ar@<-0.2pc>@{^(-}[#1] \ar@<+0.2pc>@{-}[#1]}

\title[Detecting surface bundles in finite covers.]{Detecting surface bundles in
finite covers of hyperbolic closed 3-manifolds.}

\author{Claire \textsc{Renard}.}

\begin{document}

\begin{abstract}

The main theorem of this article provides sufficient conditions
for a degree $d$ finite cover $M'$ of a hyperbolic 3-manifold $M$
to be a surface-bundle. Let $F$ be an embedded, closed and
orientable surface of genus $g$, close to a minimal surface in the
cover $M'$, splitting $M'$ into a disjoint union of $q$
handlebodies and compression bodies. We show that there exists a
fiber in the complement of $F$ provided that $d$, $q$ and $g$
satisfy some inequality involving an explicit constant $k$
depending only on the volume and the injectivity radius of $M$. In
particular, this theorem applies to a Heegaard splitting of a
finite covering $M'$, giving an explicit lower bound for the genus
of a strongly irreducible Heegaard splitting of $M'$. Applying the
main theorem to the setting of a circular decomposition associated
to a non trivial homology class of $M$ gives sufficient conditions
for this homology class to correspond to a fibration over the
circle. Similar methods lead also to a sufficient condition for an
incompressible embedded surface in $M$ to be a fiber.

\end{abstract}

\maketitle

\today

\section*{Introduction}

Thurston conjectured that every complete hyperbolic, connected and
orientable 3-manifold of finite volume virtually fibers over the
circle, i.e. such a manifold has a finite covering that is a
surface bundle over the circle.

This conjecture received a great deal of attention during the past
few years, culminating with the announcement of its proof by Ian
Agol very recently (thanks to works of Daniel Wise, Jeremy Kahn
and Vladimir Markovic, Fr\'ed\'eric Haglund, Nicolas Bergeron, and
many other people). The proof is based on Daniel Wise's program.

The aim of this article is to provide explicit criteria for a
given finite cover of a closed hyperbolic 3-manifold to be a
surface bundle. More explicitly, given a cover $M'\rightarrow M$
of $M$ with finite degree $d$, a natural question is to wonder
whether $M'$ contains an embedded surface that is a fiber, and to
give an upper bound for its genus. The idea is to start with
surfaces that already exist in $M'$, like Heegaard surfaces.

The method is inspired by Lackenby's program to find surface
bundles in towers of finite coverings of a given closed hyperbolic
3-manifold.

~

Let us be more precise. If $C$ is a handlebody or a compression
body, set $\chi_-(C):=\chi_-(\partial_+C)$. If $S$ is a union of
connected components of $\partial_-C$, the definition implies that
$\chi_-(S)\leq \chi_-(C)$.

\begin{df}

An embedded surface $S$ in a Riemannian 3-manifold $M$ is called
\textbf{pseudo-minimal} if it is orientable, closed, and $S$ is a
minimal surface or the boundary of a regular neighborhood of a
minimal non-orientable surface, possibly with a little tube
attached vertically in the $I$-bundle structure.

\end{df}

The main result of this article is the following theorem.

~

\noindent \textbf{Theorem A.} \emph{Let $M$ be a hyperbolic, connected,
oriented and closed 3-manifold. Denote by $\inj(M)$ the injectivity radius of
$M$ and set $\epsilon = \frac{1}{2} {\rm Inj}(M)$.}

\emph{There exists an explicit constant $k=k(\epsilon,\vol(M))>0$, depending
only on $\epsilon$ and the volume $\vol(M)$ satisfying the following
properties.}

\emph{Let $M'\rightarrow M$ be a cover of finite degree $d$ which contains
a closed, orientable, embedded and pseudo-minimal surface $F$, splitting $M'$
into a disjoint union of $q$ handlebodies and compression bodies
$C_1,\ldots,C_{q}$. Suppose that:}

\begin{enumerate}
 \item \emph{the union $F^{-}$ of the components of $F$ corresponding to
the negative boundary components of $C_j$ is a union of incompressible
surfaces, and}
  \item \emph{the inequality $\boldsymbol{k\,c\ln c<\ln\ln  \frac{d}{q}}$ holds,
where
$c =\max_{j=1,\ldots,q} \{\chi_-(C_j)\}$.}
\end{enumerate}

\emph{Then one of the components of $F^{-}$ is the fiber of a
surface-bundle structure for  $M'$ (corresponding to a bundle over
the circle or a twisted $I$-bundle).}

~

The proof of this theorem leads to the following corollary.

\begin{cor}\label{cor-thmA-volhandlebodies}

Under the assumptions of theorem A, the volume of a handlebody
$C_j$ (i.e. such that $\partial_-C_j=\emptyset$), among the $q$
compression bodies, must satisfy $\vol(C_j) < \vol(M) d/q$.

\end{cor}

The topological assumption $(1)$ of theorem A may not be
necessary. We conjecture that:

\begin{cjapart}[$*$]

Theorem A is still true even if assumption $(1)$ is not a priori satisfied.

\end{cjapart}

~

If $N$ is a connected, compact and orientable 3-manifold, the
Heegaard Euler characteristic $\chi_-^h(N)$ of $N$  is the minimum
over all Heegaard surfaces $F$ of the negative part
$\chi_-^{}(F)=\min\{-\chi(F),0\}$ of the Euler characteristic of
$F$. Likewise, the strong Heegaard Euler characteristic
$\chi_-^{sh}(N)$ is the minimum of $\chi_-^{}(F)$ over all the
strongly irreducible Heegaard surfaces $F$ of $N$. By convention,
if the manifold $N$ does not contain any strongly irreducible
Heegaard surface, $\chi_-^{sh}(N)=+\infty$. For further
definitions and details about the theory of Heegaard splittings,
see section \ref{rappels}.

As a Heegaard surface divides a 3-manifold into two compression
bodies, after some work, this general theorem applies in the
setting of Heegaard splittings. A consequence is the following
result, which gives a stronger and explicit version of a theorem
of J. Maher \cite{Mah}, stating that an infinite tower of finite
coverings of $M$ with a uniform bound on the Heegaard genus does
contain surface bundles. This theorem of Maher and its proof were
the starting point of this work.

\begin{thm}\label{thm-gHetfibration}

Let $M$ be a hyperbolic, connected, oriented and closed
3-manifold. Denote by $\inj(M)$ the injectivity radius of $M$ and
set $\epsilon = \frac{1}{2} {\rm Inj}(M)$.

\begin{enumerate}

\item There exists an explicit constant
$\bar{k}=\bar{k}(\epsilon,\vol(M))$ such that for every covering
$M' \rightarrow M$ with finite degree $d$ such that
$\boldsymbol{\bar{k}\, \chi_-^h(M') \ln \chi_-^h(M') \leq \ln \ln
d}$, $M'$ is a surface bundle with fiber of genus at most $g(M')$.

\item Moreover, for every covering $M' \rightarrow M$ with finite
degree $d$, one always has $\boldsymbol{\bar{k}\, \chi_-^{sh}(M')
\ln \chi_-^{sh}(M') > \ln \ln d}$.

\end{enumerate}

\end{thm}

Lackenby \cite{Lac} in his program introduced the notion of
Heegaard gradient.

\begin{df}\cite[p. 319 et 320]{Lac}

Let $M$ be a compact, connected and orientable 3-manifold. One
defines the \textbf{infimal Heegaard gradient} of the collection
of finite coverings $\rev$ with degree $d_i$ as:
$$\nabla^h(\rev)=\inf_{i\in I}\left\{\frac{\chi_-^h(M_i)}{d_i}\right\}.$$

Likewise, the \textbf{infimal strong Heegaard gradient} of the
collection $\rev$ is:
$$\nabla^{sh}(\rev)=\inf_{i\in I}\left\{\frac{\chi_-^{sh}(M_i)}{d_i}\right\},$$
where $\chi_-^{sh}(M_i)$ is the strong Heegaard Euler
characteristic of the finite covering $\rev$.

If the family of finite covers is not specified, by convention it
is the family of all finite covers of $M$. The corresponding
gradients are called the Heegaard gradient of $M$, denoted by
$\nabla^h(M)$, and the strong Heegaard gradient of $M$, denoted by
$\nabla^{sh}(M)$.

\end{df}

Results of Lackenby show that those two quantities provide
information about the existence of incompressible surfaces in
finite covers of a manifold $M$ with sufficiently large degrees.
They led Lackenby to formulate the following conjectures.

\begin{cj}[Heegaard gradient Conjecture]\cite[p.
320]{Lac}\label{cj-GH}

The Heegaard gradient of a compact, connected and orientable
hyperbolic 3-manifold is zero if and only if the manifold $M$
virtually fibers over the circle $\cercle$.

\end{cj}

This conjecture would follow immediately from the announcement of
Thurston's virtual fibration conjecture.

A second conjecture deals with the strong Heegaard genus, and
remains open.

\begin{cj}[Strong Heegaard gradient Conjecture]\cite[p.
320]{Lac}\label{cj-GFH}

The strong Heegaard gradient of a closed, connected and orientable
hyperbolic 3-manifold is always strictly positive.

\end{cj}

Theorem \ref{thm-gHetfibration} leads to a "sub-logarithmic"
version of conjecture \ref{cj-GH} (for given collections of finite coverings)
and conjecture \ref{cj-GFH}. As there exist infinite towers of
non-fibered finite coverings of a hyperbolic 3-manifold $M$ (see
\cite{BoWa}), it makes sense to asks for a condition to ensure
that a given collection $\rev$ of finite covers of $M$ contains
surface bundles.

\begin{df}

Let $\eta\in(0,1)$.

The \textbf{$\eta$-sub-logarithmic Heegaard gradient} associated
to a sequence of finite covers $\rev$ with finite degrees $d_i$ is
defined by~:
$$\nabla_{log,\eta}^h(\rev)=\inf
\left\{\frac{\chi_-^h(M_i)}{(\ln \ln d_i)^\eta}\right\}.$$

One can also define the \textbf{strong $\eta$-sub-logarithmic
Heegaard gradient} of $M$ by
$$\nabla_{log,\eta}^{sh}(M)=\inf\left\{ \frac{\chi_-^{sh}(M_i)}
{(\ln\ln d_i)^\eta}\right\},$$ where the infimum is over the
(countable) collection of all finite covers $\rev$ of $M$.

\end{df}

\begin{cor}\label{cor-gHslog}

\begin{enumerate}

\item If the $\eta$-sub-logarithmic Heegaard gradient
$\nabla_{log,\eta}^h(\rev)$ is zero, then for infinitely many
$i\in\nn$ the finite covering $M_i$ is a surface bundle.

\item The strong $\eta$-sub-logarithmic Heegaard gradient of
$M$ is always positive: $\nabla_{log,\eta}^{sh}(M) >0$.

\end{enumerate}

\end{cor}

~

Theorem A also applies in the setting of circular decompositions
associated to a non-trivial cohomology class, to give a sufficient
condition for this class to be fibered.

\begin{df}\label{def-deccirc}

Let $M$ be a hyperbolic, connected, oriented and closed
3-manifold. If $\alpha\in H^1(M)=H^1(M,\zz)$ is a non-trivial
cohomology class, let us denote by $\|\alpha \|$ the Thurston
norm of $\alpha$. By definition,
$$\|\alpha\|=\min\{\chi_-(R),\, [R]=\mathcal{P}(\alpha)\},$$
where $R$ is an embedded surface and $\mathcal{P}(\alpha)$ the
Poincar\'e-dual class of $\alpha$. We will call such a surface $R$
realizing the Thurston norm of $\alpha$ a
\textbf{$\|\alpha\|$-minimizing surface.}

If $R$ is a non-separating and $\|\alpha\|$-minimizing surface for
a given non-trivial cohomology class $\alpha\in H^1(M)$, take
$\mathcal{N}(R)\cong R\times (-1,1)$ a regular neighborhood of $R$
in $M$, and denote by $M_R=M\setminus \mathcal{N}(R)$. Set
$$h(M,\alpha,R)=\min\{\chi(R)-\chi(S)\},$$ where $S$ is a Heegaard
surface for $(M_R,R\times\{1\},R\times\{-1\})$. Said differently,
$\frac{1}{2}h(M,\alpha,R)$ is the minimal number of $1$-handles we
need to attach to a regular neighborhood of $R\times\{1\}$ in
$M_R$ to get the first compression body of a Heegaard splitting of
$(M_R,R\times\{1\},R\times\{-1\})$. Set $$h(M,\alpha)=h(\alpha)=
\min\{h(M,\alpha,R),\,[R]=\mathcal{P}(\alpha),\,
\chi_-(R)=\|\alpha\|\}.$$ For each non-trivial cohomology class
$\alpha\in H^1(M)$, let $\chi_-^c(\alpha)=\|\alpha\|+h(\alpha)$ be
the \textbf{circular characteristic} of $\alpha$. It is the
negative part of the Euler characteristic of a minimal genus
Heegaard surface for $M_R$, where $R$ is a $\|\alpha\|$-minimizing
surface such that the number $h(M,\alpha,R)$ is minimal among all
$\|\alpha\|$-minimizing surfaces.

\end{df}

In this setting, theorem A leads to the following corollary, which
is analogous to theorem \ref{thm-gHetfibration} for circular
decompositions associated to a non-trivial cohomology class.

\begin{cor}\label{cor-dec-circ}

Let $M$ be a hyperbolic, connected, oriented and closed 3-manifold.
Set $\epsilon = \inj(M)/2$. There exists an explicit constant
$\ell'=\ell'(\epsilon,\vol(M))$, depending only on $\epsilon$ and
the volume of $M$, satisfying the following property. Let
$M'\rightarrow M$ be a covering of $M$ of finite degree $d$, and
$\alpha'\in H^1(M')$ a non-trivial cohomology class.

If $\ell'\,\chi_-^c(\alpha')\ln \chi_-^c(\alpha') \leq \ln \ln d$,
then every $\|\alpha'\|$-minimizing surface $R'$ in $M'$ is a
fiber.

\end{cor}

Thus we have a criterion to ensure that a non-trivial homology
class can be represented by a fiber. If $R$ is an incompressible
embedded surface in $M$, its homology class is trivial if and only
if $R$ is separating. We have also established a sufficient
condition for an incompressible surface $R$ to be a virtual fiber.

\begin{df}\label{def-caracteristique-surface}

Let $M$ be a hyperbolic, connected, oriented and closed
3-manifold. Suppose that $R$ is an incompressible, orientable and
connected embedded surface in $M$. If $R$ is non-separating, the
homology class $[R] \in H_2(M)$ is non-trivial. Let the
\textbf{Heegaard characteristic of the surface $R$} be the minimum
of $\left|\chi(S)\right|$, where $S$ is a Heegaard surface for
$(M_R:=M\setminus \mathcal{N}(R), R\times \{1\}, R\times\{-1\})$.

If the surface $R$ is separating, the manifold $M_R:=M\setminus
\mathcal{N}(R)$ is the disjoint union of two connected components
$M_l$ and $M_r$. Let the \textbf{Heegaard characteristic of the
surface $R$} be the maximum of $\{\chi_-^h(M_l),\chi_-^h(M_r)\}$.

In both cases, let us denote by $\chi_-^h(R)$ the Heegaard
characteristic of the incompressible surface $R$.

\end{df}

In the following corollary, the surface $R$ can either be separating
or non-separating.

\begin{cor}\label{cor-surf-incomp&fibre-virtuelle}

Let $M$ be a hyperbolic, connected, oriented and closed
3-manifold, and set $\epsilon = \inj(M)/2$. There exists an
explicit constant $\ell''=\ell''(\epsilon,\vol(M))$, depending
only on $\epsilon$ and $\vol(M)$ and satisfying the following
property. Let $R$ be an incompressible, connected, orientable and
closed embedded surface in $M$. Let $M'\rightarrow M$ be a
covering of $M$ of finite degree $d$. Let also $R'$ be a connected
component of the preimage of $R$ in $M'$.

If $\ell''\,\chi_-^h(R')\ln \chi_-^h(R') \leq \ln \ln d$, then the
incompressible surface $R$ is a fiber. Moreover, if $R'$
is non-separating, $R'$ is the fiber of a bundle over the circle, and the same
holds for $R$ if it is non-separating. Otherwise, it is the fiber of a twisted
$I$-bundle.

\end{cor}

\begin{rmq}

The explicit expression of constants $k$, $\bar{k}$, $\ell'$ and
$\ell''$ involved in theorem A and corollaries of this work allows
us to study their behavior. If the volume $\vol(M)$ is fixed and
that $\epsilon$ tends to zero, or if $\epsilon$ is fixed and
$\vol(M)$ tends to infinity, all those constants tend to infinity.
Thus, the sufficient conditions become more and more difficult to
satisfy when the injectivity radius decreases, or if the volume
grows.

\end{rmq}

\noindent \textbf{Outline of the paper:} After some definitions
and generalities about Heegaard splittings in the first section,
we prove theorem A in the second section. The third section is
dedicated to the application of theorem A to Heegaard splittings
and the proof of theorem \ref{thm-gHetfibration} and corollary
\ref{cor-gHslog}. The last section deals with applications to
circular decompositions and the proof of corollaries
\ref{cor-dec-circ} and \ref{cor-surf-incomp&fibre-virtuelle}.

~

\noindent \textbf{Acknowledgement:} I would like to thank warmly my advisor,
Michel Boileau, whose encouragements, kindness and patience were
essential ingredients in this work. I am grateful to Juan Souto,
Nicolas Bergeron, Fr\'ed\'eric Paulin, Joan Porti, Jean-Marc Schlenker,
Jean-Pierre Otal, Vincent
Guirardel, Cyril Lecuire, Steven Boyer, David Gabai, Dick Canary, Thomas
Haettel, Anne Berry and Genevi\`eve Simonet for very helpful
conversations during the
elaboration of this paper.

\section{Background on Heegaard splittings.}\label{rappels}

In this section, we briefly summarize the theory of Heegaard
splittings. We also refer to \cite{Scha} for a survey on the
subject.

A \textbf{handlebody} is the regular neighborhood of a connected
graph. Its boundary is a connected, orientable and closed surface.
The genus $g$ of this surface is called the \textbf{genus} of the
handlebody. The original graph is called a \textbf{spine} for the
handlebody. If an orientable 3-manifold $M$ is closed, a
\textbf{Heegaard splitting} of $M$ is a decomposition of $M$ as
the union of two handlebodies with the same genus, glued together
by a diffeomorphism of their boundaries. A \textbf{compression
body} is a connected and orientable 3-manifold $H$ with boundary,
obtained from a regular neighborhood $S\times[0,1]$ of a closed
surface $S$, not necessarily connected. One glues some 1-handles
to the surface $S\times\{1\}$ to get the compression body $H$. The
surface $S\times\{0\}$, denoted by $\partial_-H$, is called the
\textbf{negative boundary} of the compression body $H$. The
boundary of $H$ minus the negative boundary $\partial_-H$ is a
connected surface $\partial_+H$, called the \textbf{positive
boundary} of $H$. The genus of the closed surface $\partial_+H$ is
called the \textbf{genus} of the compression body $H$ and denoted
by $g(H)$. By convention, a handlebody as defined above is a
compression body $H$ for which $\partial_- H=\emptyset$. A
\textbf{spine} for a compression body $H$ is the union $\Gamma$ of
the negative boundary $\partial_-H$ together with a graph whose
vertices lie on $\partial_-H$, such that $H$ deformation retracts
on $\Gamma$.

\begin{df}[Heegaard Splittings]

Let $(M,N_0,N_1)$ be a cobordism of $M$, with possibly
$N_0=\emptyset$ or $N_1=\emptyset$. In particular, if $M$ is
closed, $N_0=N_1=\emptyset$. A \textbf{Heegaard splitting of $M$}
associated to the cobordism $(M,N_0,N_1)$ is a decomposition of
$M$ into two compression bodies $H_0$ and $H_1$ such that:
\begin{enumerate}
\item $\partial_-H_0=N_0$, $\partial_-H_1=N_1$,
\item $\partial_+H_0\cong\partial_+H_1\cong S$ where $S$ is a closed
surface, and
\item $M=H_0\cup_S H_1$ is obtained from $H_0$ and
$H_1$ by gluing their positive boundaries by a homeomorphism of $S$.
\end{enumerate}
The surface $S$ is called a \textbf{Heegaard surface} for $M$ and
its genus is called the \textbf{genus} of the Heegaard splitting
$M=H_0\cup_S H_1$.

\end{df}

Every compact and orientable 3-manifold $M$ admits a Heegaard
splitting. The \textbf{Heegaard genus} of the manifold $M$,
denoted by $g(M)$, is the minimal genus of all Heegaard splittings
of $M$. The \textbf{Heegaard Euler characteristic} of $M$ is
$\chi_-^h(M)=2g(M)-2$, the negative part of the Euler
characteristic of a minimal genus Heegaard surface for $M$.

A \textbf{meridian disc} for a Heegaard splitting of $M$ is a
properly embedded disc in one of the compression bodies, which
bounds an essential curve in the Heegaard surface. A Heegaard
splitting (or a Heegaard surface) is said to be \textbf{strongly
irreducible} if there does not exist any pair of disjoint meridian
discs, one in each compression body. In other words, in a strongly
irreducible Heegaard splitting, the boundaries of any two meridian
discs each in one side of the Heegaard surface necessarily
intersect. For any orientable 3-manifold $M$, one defines the
\textbf{strong Heegaard Euler characteristic} $\chi_-^{sh}(M)$ of
$M$ as the minimum over all strongly irreducible Heegaard surfaces
$F$ of the negative part $\chi_-(F)$ of the Euler characteristic
of $F$. If the manifold $M$ does not
have any strongly irreducible Heegaard splitting, then
$\chi_-^{sh}(M)=+\infty$.

Note that in the case of hyperbolic 3-manifolds, the Heegaard
Euler characteristics and the strong Heegaard Euler
characteristics are always at least $2$.

A Heegaard splitting can be seen as a handle decomposition for a
closed 3-manifold $M$. Starting from a collection of 0-handles,
one attaches some 1-handles to them, then a collection of
2-handles, to finish by 3-handles. The first handlebody
corresponds to the 0- and 1-handles, to which the 2- and 3-handles
that compose the second handlebody are attached. More generally, a
\textbf{generalized Heegaard splitting} for a 3-manifold $M$
corresponds to a handle decomposition: starting from 0-handles and
possibly collars of some boundary components of $M$, one attaches
some 1-handles, then a collection of 2-handles, then another
collection of 1-handles, and so on, alternating 1- and 2-handles,
to finish after the last collection of 2-handles with a collection
of 3-handles. If one stops during the process, the object obtained
after attaching the $j$-th batch of 1- or 2-handles is a
3-manifold embedded in $M$. Let $F_j$ be its boundary, after
discarding any 2-sphere component that bounds a 0- or a 3-handle.
After a small isotopy to make all the surfaces $F_j$ disjoint, one
gets a collection of $2n-1$ disjoint surfaces $F_j$ in $M$. The
surfaces $F_{2j}$, called the \textbf{even surfaces}, separate the
manifold $M$ into $n$ 3-manifolds, for which the surfaces
$F_{2j-1}$, called the \textbf{odd surfaces}, form Heegaard
surfaces.

\begin{center}
\begin{tikzpicture}[scale=.4]

\draw (-2,0) -- (2,0) -- (3,2) -- (2,4) -- (3,6) -- (2,8) --
(3,10) -- (2,12)-- (3,14) -- (2,16) -- (-2,16) -- (-3,14) --
(-2,12) -- (-3,10) -- (-2,8) -- (-3,6) -- (-2,4) -- (-3,2) --
(-2,0);

\draw (-3,2) -- (3,2);

\draw (-2,4) -- (2,4);

\draw (-3,6) -- (3,6);

\draw (-2,8) -- (2,8);

\draw (-3,10) -- (3,10);

\draw (-2,12) -- (2,12);

\draw (-3,14) -- (3,14);

\node at (4,0) {$N_0.$};

\node at (4,2) {$F_1$};

\node at (4,4) {$F_2$};

\node at (4,6) {$F_3$};

\node at (4,8) {$F_4$};

\node at (4,10) {$\vdots$};

\node at (4.5,12) {$F_{2n-2}$};

\node at (4.5,14) {$F_{2n-1}$};

\node at (4,16) {$N_1$};

\draw (-4,0) -- (-4,16);

\draw (-4.1,0) -- (-3.9,0);

\draw (-4.1,4) -- (-3.9,4);

\draw (-4.1,8) -- (-3.9,8);

\draw (-4.1,12) -- (-3.9,12);

\draw (-4.1,16) -- (-3.9,16);

\node at (-5,2) {$M_1$};

\node at (-5,6) {$M_2$};

\node at (-5,10) {$\vdots$};

\node at (-5,14) {$M_n$};

\end{tikzpicture}
\end{center}

To each 1- and 2-handles of a generalized Heegaard splitting, one
can associate a meridian disc. If the splitting of the region
between two even surfaces is not strongly irreducible, two
disjoint meridian discs can be used to change the order in which
the handles are attached. A 2-handle corresponding to one of the
meridian discs can be attached before a 1-handle corresponding to
the other meridian disc. We will call this operation a
\textbf{surgery of generalized Heegaard splitting}.

Let $F$ be a closed and orientable surface. If $F$ is connected,
one defines the \textbf{complexity} of $F$ as $c(F)=0$ if $F$ is
the 2-sphere $\sphere$, and $c(F)=2g(F)-1=1-\chi(F)$ otherwise. If
$F$ is not connected, the complexity of $F$ is the sum over all
components of $F$ of the complexity of the component.

If $\mathcal{H}=\{F_1,F_2,\ldots,F_{2n-1}\}$ is a generalized
Heegaard splitting of $M$, the \textbf{width} of this
decomposition is the set
$w(\mathcal{H})=\{c(F_1),\ldots,c(F_{2n-1})\}$ of the complexities
of the odd surfaces, with repetitions and arranged in
monotonically non-increasing order. Widths are compared using the
lexicographic order.

Starting from a generalized Heegaard splitting
$\mathcal{H}=\{F_1,F_2,\ldots,F_{2n-1}\}$ in which at least one of
the surfaces $F_{2i-1}$ is not strongly irreducible, one can do a
surgery of generalized Heegaard splittings to change the order in
which the 1- and 2-handles are attached, to get a new generalized
Heegaard splitting $\mathcal{H}'$ with
$w(\mathcal{H}')<w(\mathcal{H})$.

\begin{center}
\begin{tikzpicture}[scale=.7]

\draw (-4,6) -- (4,6);

\draw [thick] (-2.25,6) -- (-2.25,8);

\draw [thick] (2.25,6) -- (2.25,4);

\node at (5,6) {$F_{2i-1}$};

\node at (-1.75,7) {$D_1$};

\node at (2.75,5) {$D_0$};

\draw [->] (0,3.8) -- (0,2.2);

\draw [very thick] (-4,0) -- (-2.5,0) -- (-2.5,2);

\draw [very thick] (-2,2) -- (-2,0) -- (2,0) -- (2,-2);

\draw [very thick] (2.5,-2) -- (2.5,0) -- (4,0);

\draw (-4,.5) -- (-3,.5) -- (-3,2);

\draw (-1.5,2) -- (-1.5,.5) -- (4,.5);

\draw (4,-.5) -- (3,-.5) -- (3,-2);

\draw (1.5,-2) -- (1.5,-.5) -- (-4,-.5);

\draw [dashed] (-2.25,0) -- (-2.25,2);

\draw [dashed] (2.25,0) -- (2.25,-2);

\node at (4.8,0) {$G_2$};

\node at (4.8,.5) {$G_3$};

\node at (4.8,-.5) {$G_1$};

\node at (-6,6) {$\mathcal{H}$};

\node at (-6,0) {$\mathcal{H}'$};

\end{tikzpicture}
\end{center}

If $\mathcal{H}$ is a generalized Heegaard splitting for $M$, let
$\mathcal{S}_{\mathcal{H}}$ be the set of all generalized Heegaard
splittings obtained from $\mathcal{H}$ by surgery. An element
$\mathcal{H}'\in \mathcal{S}_{\mathcal{H}}$ of minimal width is
called an \textbf{$\mathcal{H}$-thin generalized Heegaard
splitting}.

\begin{prop}\label{prop-amincir}

Let $M$ be a connected, oriented and compact 3-manifold, and
$\mathcal{H}$ a generalized Heegaard splitting for $M$.

Every $\mathcal{H}$-thin generalized Heegaard splitting
$\mathcal{H'}=(F_1,\ldots, F_{2n-1})$ satisfies the following
properties.

\begin{enumerate}
\item The odd surfaces $F_{2i-1}$ correspond to strongly irreducible
Heegaard surfaces.
\item The even surfaces $F_{2i}$ are incompressible surfaces in $M$.
\item Furthermore, if the manifold $M$ is irreducible, then no component
of any even surface is a 2-sphere.
\end{enumerate}

\end{prop}\qed

The proof of this proposition is a consequence of the definition
of a surgery of generalized Heegaard splittings. See for example
\cite{CaGo} and \cite{ST2}.

A generalized splitting of minimal width among all generalized
Heegaard splittings of $M$ is called a \textbf{thin position} (see
\cite{ST2}).

\section{Finding a fibration.}

\subsection{Main theorem.}\label{par_thmprincipal}

The aim of this section is to prove theorem A.

If $S$ is a surface, let us denote by
$\chi_-(S)=\max\{0,-\chi(S)\}$ the negative part of the Euler
characteristic of $S$.

If $C$ is a compression body, set
$\chi_-(C):=\chi_-(\partial_+C)$. If $S$ is a union of connected
components of $\partial_-C$, the definition implies that
$\chi_-(S)\leq \chi_-(C)$.

\begin{df}

An embedded surface $S$ in a Riemannian 3-manifold $M$ is called
\textbf{pseudo-minimal} if it is orientable, closed, and $S$ is a
minimal surface or the boundary of a regular neighborhood of a
minimal non-orientable surface, possibly with a little tube
attached vertically in the $I$-bundle structure.

\end{df}

~

\noindent \underline{Proof of theorem A.}
\nopagebreak

Suppose that there exists a finite cover $M'\rightarrow M$ of
degree $d$ satisfying the hypotheses of theorem A. The proof relies on
three key propositions, the proof of which we postpone to the
three next subsections. Let us denote by $g=\frac{c}{2}+1$. It is
an upper bound for the genus of the compression bodies of $M'$.

\begin{lm}\label{lm-trouver-corps-a-anses}

There exists a compression body $C$ among the $q$ compression
bodies $C_1,\ldots,C_{q}$ of $M'$ such that
$$\vol(C)\geq \vol(M)\, \frac{d} {q}.$$

\end{lm}

\noindent \underline{Proof of lemma \ref{lm-trouver-corps-a-anses}.}
\nopagebreak

The proof is straightforward, as there are $q$ compression bodies
$C_1,\ldots,C_{q}$, and $\vol(M')=d \,\vol(M)$.\qed

~

Let $C$ be a compression body as in lemma \ref{lm-trouver-corps-a-anses}.

\begin{lm}\label{lm-chgts-negligeables}

Let $k_0=\max\left\{ \frac{\ln(4(2\epsilon a'+1))} {2\ln 2},
1+\frac{\ln(1+\ln(12V_3 /\vol(M))}{2\ln 2}\right\}$, where $V_3$ is the maximal
volume of an ideal hyperbolic tetrahedron in $\H3$, and
$a'=6(21/4+3/4\pi+3/4\epsilon+2/\sinh^2(\epsilon/4))$.

If $k_0\,\chi_-(C)\ln \chi_-(C)\leq \ln \ln d/q $ and $\vol(M)\geq
\pi/2$, then there is a way of replacing the boundary surfaces of $C$ by
simplicial surfaces, to obtain a new compression body $C''$ with:
$$\vol(C'')\geq \frac{1}{4}\vol(C)\geq \frac{\vol(M) d}{4q}.$$

\end{lm}

This lemma is proven in subsection \ref{par-modifications-cea}. To simplify
notations, this new compression body $C''$ will still be denoted by $C$.

\begin{df}\label{df-separer}

Let $x$ be a point in $C$ and $S$ an immersed surface in $C$.
We say that \textbf{$S$ separates $x$ from $\partial_+C$} if
every oriented path from $x$ to $\partial_+C$ has its algebraic
intersection number with $\partial_+C$ equal to $+1$.

If two surfaces $S$ and $T$ immersed in $C$ are such that $S$
separates every point of $T$ from $\partial_+C$, we say that
\textbf{$T$ separates $S$ from $\partial_+C$}. In this case, the
surfaces $S$ and $T$ are said to be \textbf{nested}.

\end{df}

We will denote the ceil function of the real number $x$ by $\lceil
x\rceil$, i.e. the smallest integer not less than $x$. Similarly,
$\lfloor x \rfloor$ is the floor function of $x$, and represents
the largest integer no greater than $x$. By convention, we set
$\lceil x \rceil$ and $\lfloor x \rfloor$ equal to zero if $x$ is
non-positive.

The following proposition is a step towards the construction of a
certain amount of parallel surfaces in the compression body $C$.
It is an adaptation of Lemma 4.5 p. 2251 of \cite{Mah}. We
postpone its proof to section \ref{dem-lm4.5}.

~

\noindent \textbf{Proposition B (of Nested Surfaces).}
\nopagebreak

\emph{Let $\delta$ be the diameter of the compression body $C$ of
$M'$, $\epsilon= \inj(M)/2$,
$K=4\left(3+\frac{1}{\sinh{\frac{\epsilon}{8}}^2}\right)g(C)-10$
and $K'=2a'\chi_-(C)$. Moreover, suppose that $k_0\,\chi_-(C)\\ \ln
\chi_-(C) \leq \ln\ln \frac{d}{q}$.}

\emph{Under those assumptions, there exist at least $n=\lceil
\frac{\delta}{36\epsilon K}-\frac{2}{9}- \frac{K'}{3K}\rceil$ orientable,
disjoint and nested surfaces, immersed in $C$. All of those
surfaces are homotopic to components of surfaces obtained by
compressing $\partial_+C$. Moreover, the $\epsilon$-diameter
of those surfaces in  $M'$ is bounded from above by $K$ and they
are separated from each other by a distance greater than or equal
to $10\epsilon K$.}

~

With this proposition, we obtain at least $n=\lceil
\frac{\delta}{36\epsilon K}-\frac{2}{9}- \frac{K'}{3K}\rceil$
nested immersed surfaces in the handlebody $C$. Those surfaces are
all disjoint and homotopic to components of surfaces obtained from
$\partial_+C$ by surgery. This implies that the the genus of those
surfaces is between 0 and $g(C)$, the genus of $C$ (which is, by
assumption, less than or equal to $g$). We can thus find at least
$\lfloor\frac{n}{g(C)+1}\rfloor$ such nested immersed surfaces of
the same genus. The next step is then to replace those nested
immersed surfaces by parallel embedded surfaces.

~

\noindent \textbf{Proposition C (of Parallel Surfaces).}
\nopagebreak

\emph{Let $\delta$ be the diameter of the compression body $C$
in $M'$, $\epsilon= \inj(M)/2$,
$K=4\left(3+\frac{1}{\sinh{\frac{\epsilon}{8}}^2}\right)g(C)-10$
and $K'=2a'\chi_-(C)$. Suppose that $k_0\,\chi_-(C)\ln \chi_-(C)
\leq \ln\ln \frac{d}{q}$.}

\emph{Under those assumptions, there exists at least $m=(\lfloor
\frac{1}{g(C)+1}\lceil \frac{\delta}{36\epsilon K}-\frac{2}{9}-
\frac{K'}{3K}\rceil\rfloor -4)$ orientable, parallel and connected surfaces
embedded in $C$, separated from each other by
a distance greater than or equal to $\epsilon K$, and each
of those surfaces can be covered by at most $K$ embedded balls in
$M'$ of radius $2\epsilon$. In particular, their diameter in the
manifold $M'$ is uniformly bounded from above by $4\epsilon K$.}

~

For the proof of this proposition, see section
\ref{dem-prop-surf//}.

Let
 \begin{eqnarray*}
 a&=&2\left(3+\frac{1}
{\sinh^2({\frac{\epsilon}{8}})}\right),\\
 b&=&2\left(1+\frac{2}
{\sinh^2({\frac{\epsilon}{8}})}\right)\,{\rm and}\\
 a'&=&6\left(\frac{21}{4}+\frac{3}{4\pi} + \frac{3}{4\epsilon}
 +\frac{2}{\sinh^2(\frac{\epsilon}{4})}\right).
 \end{eqnarray*}

\begin{lm}\label{lm-minoration-m}

Under assumptions of theorem A, according to proposition C, there
exist $m$ parallel surfaces embedded in the compression body $C$
of $M'$, with
$$m\geq  \frac{2}{\chi_-(C)+4}\left(
\frac{\ln\left(\frac{d}{q}\right)+
\ln\left(\frac{\vol(M)}{2\pi}\right)}{72\epsilon
(a\chi_-(C)+b)}-\frac{2}{9}(1+ \frac{3a'}{a})\right) -5.$$

\end{lm}

\noindent \underline{Proof of lemma \ref{lm-minoration-m}.}
\nopagebreak

The number $m$ of embedded parallel surfaces in $C$ obtained by
proposition C is equal to:
$$m=\lfloor \frac{1}{g(C)+1}\lceil \frac{\delta}{36\epsilon
K}-\frac{2}{9}- \frac{K'}{3K}\rceil \rfloor-4,$$ where
\begin{eqnarray*}
K&=& 4\left(3+\frac{1}
{\sinh{\frac{\epsilon}{8}}^2}\right)g(C)-10\\
&=& a\chi_-(C)+b,
\end{eqnarray*}
 and $$K'=2a'\chi_-(C).$$

~

The diameter $\delta$ of the compression body $C$ and the ratio
$d/q$ are related. On the one hand,
$$\vol(C)\leq\vol\left(\bb_{\mathbb{H}^3}\left(\delta
\right)\right)
=\pi(\sinh(2\delta)-2\delta)\leq\frac{\pi}{2}e^{2\delta}.$$

\begin{rmq}
The second logarithm of the expression $\ln\ln \frac{d}{q}$ comes
from this estimation linking the diameter with the volume of a
hyperbolic 3-manifold.

\end{rmq}

On the other hand, lemmas \ref{lm-trouver-corps-a-anses} and
\ref{lm-chgts-negligeables}  give the lower bound
$$\vol(C)\geq
\vol(M)\,\frac{d}{4q},$$ which leads to the inequality
\begin{equation}\label{diametre}
\delta\geq\frac{1}{2}\ln\left(\frac{d}{q}\right)+ \frac{1}{2}
\ln\left(\frac{\vol(M)}{2\pi}\right).
\end{equation}
In particular, if $d/q$ tends to infinity, $\delta$ tends also to
infinity.

The expression of $m$ involves the ratio $\frac{K'}{3K}$. Now,
\begin{eqnarray*}
\frac{K'}{3K}&=&\frac{2a'\chi_-(C)} {3a\chi_-(C)+3b}\\
&=& \frac{2a'}{3a+3b/\chi_-(C)}\\
&\leq &\frac{2a'}{3a}.
\end{eqnarray*}
Replacing the ratio $K'/3K$ by $2a'/3a$ and taking inequality
(\ref{diametre}) into account, one obtains

\begin{eqnarray*}
m&\geq& \lfloor \frac{2}{\chi_-(C)+4}\lceil
\frac{\ln\left(\frac{d}{q}\right)+
\ln\left(\frac{\vol(M)}{2\pi}\right)}{72\epsilon
(a\chi_-(C)+b)}-\frac{2}{9}- \frac{2a'}{3a}\rceil \rfloor-4\\
&\geq & \frac{2}{\chi_-(C)+4}\left(
\frac{\ln\left(\frac{d}{q}\right)+
\ln\left(\frac{\vol(M)}{2\pi}\right)}{72\epsilon
(a\chi_-(C)+b)}-\frac{2}{9}(1+ \frac{3a'}{a})\right) -5,
\end{eqnarray*}
which ends the proof of lemma \ref{lm-minoration-m}.\qed

~

Those $m$ parallel surfaces obtained by proposition C are candidates
for a fiber. But we still have to select some of them to
get a virtual fibration of the base manifold $M$.

Let $\mathcal{D}$ be a Dirichlet fundamental polyhedron for $M$ in
its universal cover $\widehat{M}\simeq \mathbb{H}^3$. Translates
of $\mathcal{D}$ by the covering transformations give a tiling of
the universal cover $\widehat{M}$. This tiling descends to a
tiling of the finite cover $M'$ by $d$ copies of $\mathcal{D}$.
Each of the $m$ parallel, connected and embedded surfaces in
$M'$ obtained by proposition C intersects a finite
and connected set of copies of $\mathcal{D}$. We call such a set a
\textbf{pattern of fundamental domains}. We can suppose that each
of the embedded surfaces is transverse to the 2-skeleton
of the tiling. More precisely, we can suppose that each surface
does not meet the vertices of the fundamental polyhedra, that it
intersects the edges in isolated points and it is transverse to
the 2-dimensional faces of the polyhedra. Thus a pattern of
fundamental domains is a connected set that is the union of copies
of $\mathcal{D}$ glued along some of their 2-dimensional faces.

\begin{lm}\label{lm-nb-de-motifs}

Let $\mathcal{D}$ be a Dirichlet fundamental polyhedron for $M$ in
$\mathbb{H}^3$. Let $\alpha$ be the number of faces of
$\mathcal{D}$ of dimension two.

For each $\ell\in\nn$, the number of possibilities to glue
together at most $\ell$ copies of $\mathcal{D}$ to form a pattern
of $\ell$ fundamental domains is less than or equal to
$\left(\alpha \sqrt{2}\; \ell \right)^{\alpha \ell}$.

\end{lm}

\noindent \underline{Proof of lemma \ref{lm-nb-de-motifs}.}
\nopagebreak

For every $\ell\in\nn$, let us denote by $B(\ell)$ the number of
possibilities to glue together $\ell$ copies of $\mathcal{D}$ to
form a pattern of $\ell$ fundamental domains. We have to find an
upper bound for the number of possibilities to identify pairwise
some 2-dimensional faces of at most $\ell$ Dirichlet polyhedra.

First, let us notice that there are at most $\alpha\ell$ such
2-dimensional faces. Thus, there are at most $(\alpha\ell)!\leq
(\alpha\ell)^{\alpha\ell}$ ways to match pairwise those faces.

If $(F_1,F_2)$ is a pair of such faces, we can choose to glue them
together by an orientation-reversing isometry
$h\,:\,F_1\longrightarrow F_2$ (if such an isometry between those
two faces exists). This isometry corresponds to a "pairing
transformation" (see for example \cite[Proposition 3.5.1 p.
117]{Mar}). It is a reflection in $\H3$ and its hyperplane contains
one of the faces of $\mathcal{D}$. Thus, if such an isometry
exists, it is unique. We can also decide not to glue those two
faces together: by convention, we will say that we glue them by
the empty gluing. Therefore, there are at most $2$ ways to glue
$F_1$ and $F_2$ together, including the empty gluing.

Thus there are at most $(\alpha\ell)!2^{
\frac{\alpha\ell}{2}}\leq \left(\alpha \sqrt{2}\; \ell
\right)^{\alpha \ell}$ ways to glue together at most $\ell$ copies
of fundamental domains to form a pattern of fundamental domains,
which ends the proof of lemma \ref{lm-nb-de-motifs}.\qed

~

The following lemma is a way to bound the number $\alpha$ of
2-faces of a fundamental polyhedron $\mathcal{D}$ and its diameter
in $\H3$ in terms of the volume of the manifold $M$ and a lower
bound for its injectivity radius.

\begin{lm}\label{lm-estimation_Dalpha&L}

Let $\mathcal{D}$ be a Dirichlet fundamental polyhedron for the
manifold $M$, embedded in the universal cover $\widetilde{M}\simeq
\mathbb{H}^3$. Let $D$ be an upper bound for the diameter of
$\mathcal{D}$ in $\H3$ and $\alpha$ be the number of its 2-faces.
We have the following estimates:

\begin{equation}\label{estimation_D}
\diam(\mathcal{D})\leq \frac{8\epsilon \vol(M)}{\pi(\sinh
(2\epsilon)-2\epsilon)}=D,
\end{equation} and

\begin{equation}\label{estimation_alpha}
\alpha\leq \frac{\pi(\sinh(4D)-4D)}{\vol(M)}-1.
\end{equation}

If $S$ is an embedded surface in the finite cover $M'$ of $M$,
which can be covered by at most $K$ embedded balls in $M'$
of radius $\epsilon'\leq \inj(M)$, then $S$ intersects at most $L$
images of $\mathcal{D}$ in $M'$, with

\begin{equation}\label{estimation_L}
L=\lfloor \frac{\pi(\sinh(2D+2\epsilon')-2D-2\epsilon')}
{\vol(M)}K\rfloor.
\end{equation}

\end{lm}

\noindent \underline{Proof of lemma \ref{lm-estimation_Dalpha&L}.}
\nopagebreak

To prove inequality (\ref{estimation_D}), first notice that
$\diam(\mathcal{D})\leq 2\,\diam(M)$. To prove it, recall that
there exists $w\in\mathbb{H}^3$ such that
$\mathcal{D}=\{x\in\mathbb{H}^3\,,\,\d(\gamma(w),x)\geq\d(w,x)\;\forall
\gamma\in \pi_1(M)\}$. If $x$ and $y\in\mathcal{D}$ satisfy
$\d(x,y)=\diam(\mathcal{D})$, then
$$\diam(\mathcal{D})=\d(x,y)\leq\d(x,w)+\d(y,w)\leq 2\,\diam(M).$$

Take $x$ and $y\in M$ such that $\d(x,y)=\diam(M)$, and let
$\gamma$ be a minimizing geodesic from $x$ to $y$. By definition,
$\lgr(\gamma)=\diam(M)$. Let $\mathcal{B}$ be a collection of
points in $\gamma$ which is maximal among collections of points of
$\gamma$ such that two balls of radius $\epsilon$ and whose
centers are two distinct points of $\mathcal{B}$ have disjoint
interiors. Then, by maximality of $\mathcal{B}$, the union of
balls with centers in $\mathcal{B}$ and radius $2\epsilon$ cover
the geodesic $\gamma$.

Thus, $\left|\mathcal{B}\right| \geq \frac{\lgr(\gamma)}
{4\epsilon}$. As balls of centers in $\mathcal{B}$ and radius
$\epsilon$ have disjoint interiors, considering volumes, we
deduce:

\begin{eqnarray*}
\vol(M)&\geq & \sum_{u\in\mathcal{B}}\vol(B(u,\epsilon))\\
&\geq & \frac{\lgr(\gamma)}{4\epsilon} \vol(B_{\mathbb{H}^3}(\epsilon))\\
&\geq &\frac{\diam(M)}{4\epsilon} \pi(\sinh(2\epsilon)
-2\epsilon),
\end{eqnarray*}
proving inequality (\ref{estimation_D}).

Let us show inequality (\ref{estimation_alpha}). To each 2-face of
$\mathcal{D}$, one can associate a unique translate
$g_F(\mathcal{D})$ of $\mathcal{D}$ adjacent to $\mathcal{D}$
along $F$. As the diameter of $g_F(\mathcal{D})$ is also
$\diam(\mathcal{D})\leq D$, every point $x\in g_F(\mathcal{D})$
lies at distance at most $\dist(x,F)+\diam(\mathcal{D})\leq 2D$
from $w\in \mathcal{D}$. Thus, the ball of center $w$ and radius
$2D$ contains the fundamental polyhedron $\mathcal{D}$ together
with the union of all its translates $g_F(\mathcal{D})$, where $F$
is a 2-face of $\mathcal{D}$. As those polyhedra have disjoint
interiors, for volumes, we obtain:
$$(\alpha+1)\vol(\mathcal{D})\leq \vol(B_{\mathbb{H}^3}(w, 2D)),$$ and thus
$$\alpha \leq \frac{\pi(\sinh(4D)-4D)}{\vol(M)}-1.$$

The proof of inequality (\ref{estimation_L}) is similar. Denote by
$\mathcal{B}$ the set of the centers of a collection of $K$
embedded balls in $M'$ of radius $\epsilon'$ covering the surface
$S$. Let $\mathcal{N}=\cup_{x\in\mathcal{B}}B(x,D+\epsilon')$.
Those balls are not necessarily isometric to hyperbolic embedded
balls in $\H3$ as $D+\epsilon'>\inj(M)$. However, let us show that
$\mathcal{N}$ contains every fundamental polyhedron of $M'$
intersecting $S$.

To prove it, let $x$ be a point in a fundamental polyhedron of
$M'$ intersecting $S$. Take $y\in S$ such that
$\d(x,y)=\dist(x,S)\leq D$. As $y$ is a point of $S$, there exists $x\in
\mathcal{B}$ such that the ball $B(x,\epsilon')$ contains $y$.
Therefore $\d(z,x)\leq \d(z,y)+\d(y,x) \leq D+\epsilon'$, showing
that $z\in B(x,\epsilon'+D)\subset \mathcal{N}$.

Comparing volumes, we get:
\begin{eqnarray*}
L\,\vol(\mathcal{D})&\leq &\vol(\mathcal{N})\\
L\,\vol(M)&\leq & \left|\mathcal{B}\right|\vol(B_{\H3}(\epsilon'+D))\\
L&\leq
&\frac{\pi(\sinh(2\epsilon'+2D)-2\epsilon'-2D)}{\vol(M)}K,
\end{eqnarray*}
proving inequality (\ref{estimation_L}), as
$L$ is a natural integer.\qed

~

The following key proposition is a quantitative version of Lemma
4.12 p. 2258 of \cite{Mah}. We postpone its proof to section
\ref{dem-prop-motifs-dom-fond}.

~

\noindent \textbf{Proposition D (Pattern Proposition).}
\nopagebreak

\emph{Assume that in the cover $M'$ we have $m$ connected, orientable, embedded,
disjoint and parallel surfaces,
at distance at least $r>0$ from each other. Moreover, suppose that
each of those surfaces can be covered by at most $K$
embedded balls in $M'$ of radius $\epsilon'\leq \inj(M)$.}

\emph{Let $\mathcal{D}$ be a Dirichlet fundamental domain for the
manifold $M$ in its universal cover $\widehat{M}\simeq
\mathbb{H}^3$. Let us denote by $D$ an upper bound for the
diameter of $\mathcal{D}$ and $\alpha$ an upper bound for the
number of its 2-dimensional faces.}

\emph{For all $\ell\in\nn$, let $B(\ell)$ be an upper bound for the
number of possibilities of patterns obtained by gluing together at
most $\ell$ fundamental domains that intersect a connected,
orientable and embedded surface. Let $L=\lfloor \frac{\pi(\sinh(2D
+2\epsilon')-2D-2\epsilon')}{\vol(M)}K\rfloor$.}

\emph{If $r/ (2D+1)\geq 1$ and $\frac{m}{\alpha^2L^2B(L)}\geq 4$,
or if $r/(2D+1)\leq 1$ and $\left(\frac{r}{2D+1}m-1\right)
\frac{1}{\alpha^2L^2B(L)}\geq4$, then the manifold $M$ virtually
fibers over the circle $\cercle$, and the $m$ parallel surfaces in
$M'$ are fibers of a bundle over the circle or of a twisted $I$-bundle.}

\begin{rmq}

The first logarithm in the expression $\ln\ln\frac{d}{q}$ and the
function of the complexity $c\ln(c)$ in the assumption
$k\,c\ln(c)<\ln \ln \frac{d}{q}$ arise from the use of lemma
\ref{lm-nb-de-motifs} (providing an estimate of the number
$B(\ell)$) in the proof of this proposition.

\end{rmq}

We can now finish the proof of theorem A assuming propositions B,
C and D, which will be proved in next sections.

~

The aim is to apply Proposition D to the $m$ parallel
surfaces obtained in Proposition C, with
\begin{align*}
K &= a\chi_-(C)+b,\text{ and}\\
r&=\epsilon K = \epsilon (a \chi_-(C) + b).
\end{align*}
Set\begin{eqnarray*}
D &:=&\frac{8\epsilon \vol(M)}{\pi(\sinh(2\epsilon)-2\epsilon)},\\
\alpha &:=&\frac{\pi(\sinh(4D)-4D)} {\vol(M)} -1,\text{ and}\\
\sigma &:=&\frac{\pi(\sinh(2D+4\epsilon)-2D-4\epsilon)}{\vol(M)}.
\end{eqnarray*} From lemma \ref{lm-estimation_Dalpha&L},
$D$ is an upper bound for the diameter of the fundamental
polyhedron $\mathcal{D}$, and the number of 2-faces of
$\mathcal{D}$ is at most $\alpha$.

In addition, from lemma \ref{lm-estimation_Dalpha&L} again, $L=
\lfloor \frac{\pi(\sinh(2D+4\epsilon)-2D-4\epsilon)}
{\vol(M)}K\rfloor$. In particular,
\begin{eqnarray*}
L \leq \frac{\pi(\sinh(2D+4\epsilon)-2D-4\epsilon)}{\vol(M)}
(a\chi_-(C)+b)= \sigma (a \chi_-(C) + b).
\end{eqnarray*}

~

\noindent \textbf{Claim 1.} \emph{There exist $c_1\geq 2$ and
$k_1>0$, depending only on $\epsilon$ and $\vol(M)$, such that if
$\chi_-(C) \leq c_1$ and $k_1\,\chi_-(C)\ln \chi_-(C)\leq \ln\ln
d/q$, then assumptions of Proposition D are satisfied. In
particular, $M$ virtually fibers over the circle and the $m$
embedded surfaces in $M'$ are fibers.}

\emph{Furthermore, one can take $c_1:=\frac{1}{a} \left(\frac{16
\vol(M)}{\pi(\sinh(2\epsilon)-2\epsilon)}+\frac{1}{\epsilon}-b\right)$
and
\begin{eqnarray*}
k_1&:=& \frac{1}{2\ln 2} \ln \boldsymbol{(}72(2D+1) (c_1+4) \left(3+2(\alpha
\sigma)^2 (ac_1+b)^2 (\sqrt{2}\alpha \sigma
(ac_1+b))^{\alpha \sigma (ac_1+b)}\right)\\
 &+& 16\epsilon
(1+\frac{3a'}{a}) (ac_1+b) -\ln\left(\frac{\vol(M)}{2\pi}\right)
\boldsymbol{)}.
\end{eqnarray*}}

~

\noindent \underline{Proof of claim 1.}
\nopagebreak

Recall that $r=\epsilon (a \chi_-(C) + b)$ and $2D+1 =
\frac{16\epsilon \vol(M)}{\pi(\sinh(2\epsilon)-2\epsilon)}+1$.
Thus, if $\chi_-(C)\leq c_1=\frac{1}{a} \left(\frac{16
\vol(M)}{\pi(\sinh(2\epsilon)-2\epsilon)}+\frac{1}{\epsilon}-b\right)$,
then $r\leq 2D+1$. Assumptions of the second case of Proposition D
are then satisfied if $\left(\frac{r}{2D+1}m-1\right)
\frac{1}{\alpha^2L^2B(L)}\geq4$.

Taking lemma \ref{lm-minoration-m} and the expression of $r$ into
account, one obtains the sufficient condition:

\begin{eqnarray*}
\left(\frac{2\epsilon(a\chi_-(C)+b)}{(2D+1)(\chi_-(C)+4)}\left(
\frac{\ln\left(\frac{d}{q}\right)+
\ln\left(\frac{\vol(M)}{2\pi}\right)}{72\epsilon
(a\chi_-(C)+b)}-\frac{2}{9}(1+ \frac{3a'}{a})\right)
-6\right)\frac{1}{\alpha^2L^2B(L)}\geq 4.
\end{eqnarray*}

Replace $L$ by its upper bound $\sigma(a\chi_-(C)+b)$. From lemma
\ref{lm-nb-de-motifs}, one can chose for $B(L)$ the function $B(L)
=(\sqrt{2}\alpha L)^{\alpha L} \leq
(\sqrt{2}\alpha\sigma(a\chi_-(C)+b))^{\alpha\sigma(a\chi_-(C)+b)}$.

Thus one obtains the sufficient condition
\begin{eqnarray*}
\left(\frac{2\epsilon(a\chi_-(C)+b)}{(2D+1)(\chi_-(C)+4)}\left(
\frac{\ln\left(\frac{d}{q}\right)+
\ln\left(\frac{\vol(M)}{2\pi}\right)}{72\epsilon
(a\chi_-(C)+b)}-\frac{2}{9}(1+ \frac{3a'}{a})\right) -6\right) \geq\\
4 (\alpha\sigma)^2(a\chi_-(C)+b)^2
(\sqrt{2}\alpha\sigma(a\chi_-(C)+b))^{\alpha\sigma(a\chi_-(C)+b)}.
\end{eqnarray*}

Under assumptions of claim 1, $2\leq \chi_-(C) \leq c_1$. One
can then easily check that if $k_1\,\chi_-(C)\ln \chi_-(C)\leq
\ln\ln d/q$, the sufficient condition above is satisfied.\qed

~

\noindent \textbf{Claim 2.} \emph{Suppose that $\vol(M)\geq 2\pi$.
There exist $c_2\geq c_1$ and $k_2>0$, depending only on
$\epsilon$ and $\vol(M)$, such that if $\chi_-(C) \geq c_2$ and
$k_2\,\chi_-(C)\ln \chi_-(C)\leq \ln\ln d/q$, then assumptions of
Proposition D are satisfied. In particular, $M$ virtually fibers
over the circle and the $m$ embedded surfaces in $M'$ are fibers.}

\emph{Furthermore, one can take $k_2:=4\alpha \sigma a$, and
\begin{align*}
 c_2&=\max\{c_1, \frac{1}{a}\left( \frac{\ln 5
-\ln(4\alpha^2\sigma^2 (2a+b)^2)}
{\alpha \sigma  \ln(2\sqrt{2}\alpha \sigma a)}-b\right),\\
&\frac{1}{a}\left(\frac{\ln(1+\frac{3a'}{a}) -\ln
(108\alpha^2\sigma^2(2a+b)^2)}{\alpha \sigma \ln(\sqrt{2}
\alpha\sigma (2a+b))} -b\right), b/a, 4,
2\sqrt{2}\alpha \sigma a,\\
&\frac{b}{a}+ \frac{4}{\alpha \sigma a},
 \frac{\ln(18432 \epsilon \alpha^2\sigma^2a^3
(2\sqrt{2}\alpha\sigma a)^{\alpha \sigma b})}{\alpha \sigma a \ln2
},\\
&\frac{1}{a}
\left(\frac{1}{\alpha \sigma \ln(\sqrt{2}\alpha \sigma (2a+b))} \ln
\left(\frac{1}{4\alpha^2\sigma^2(2a+b)^2} (\frac{\left| -\ln
\frac{\vol(M)}{2\pi}
- \frac{2}{9}(1+\frac{3a'}{a}) \right|}{216 \epsilon (2a+b)}-5 )\right)
-b\right) \}.
\end{align*}}

~

\noindent \underline{Proof of claim 2.}
\nopagebreak

As $\chi_-(C)\geq c_2 \geq c_1$, from the proof of the first
claim, $r\geq 2D+1$. Assumptions of the first case of Proposition
D are then satisfied if $\frac{m}{\alpha^2L^2B(L)}\geq 4$. Now,
taking lemma \ref{lm-minoration-m} into account, together with the
inequalities $L\leq\sigma (a\chi_-(C)+b)$ and $B(L)\leq
(\sqrt{2}\alpha\sigma (a\chi_-(C)+b))^{\alpha  \sigma
(a\chi_-(C)+b)}$, one obtains the following sufficient condition:
\begin{eqnarray*}
\frac{2}{\chi_-(C)+4}\left( \frac{\ln\left(\frac{d}{q}\right)+
\ln\left(\frac{\vol(M)}{2\pi}\right)}{72\epsilon
(a\chi_-(C)+b)}-\frac{2}{9}(1+ \frac{3a'}{a})\right) -5 \geq\\
4\alpha^2\sigma^2 (a\chi_-(C)+b)^2 (\sqrt{2}\alpha\sigma
(a\chi_-(C)+b))^{\alpha  \sigma (a\chi_-(C)+b)},
\end{eqnarray*} which can also be written
\begin{eqnarray*}
\ln\left(\frac{d}{q}\right) &\geq &
72\epsilon (a\chi_-(C)+b) ( \frac{\chi_-(C)+4}{2} \boldsymbol{(}
4\alpha^2\sigma^2 (a\chi_-(C)+b)^2 \\
&&(\sqrt{2}\alpha\sigma (a\chi_-(C)+b))^{\alpha  \sigma (a\chi_-(C)+b)}+ 5
\boldsymbol{)} + \frac{2}{9}(1+ \frac{3a'}{a})  ) - \ln \frac{\vol(M)}{2\pi},
\end{eqnarray*} or also
\begin{eqnarray*}
\ln \ln \left(\frac{d}{q}\right) &\geq &
\ln \boldsymbol{(}72\epsilon (a\chi_-(C)+b) ( \frac{\chi_-(C)+4}{2}
\boldsymbol{(} 4\alpha^2\sigma^2 (a\chi_-(C)+b)^2 \\
&&(\sqrt{2}\alpha\sigma (a\chi_-(C)+b))^{\alpha  \sigma (a\chi_-(C)+b)}+ 5
\boldsymbol{)} + \frac{2}{9}(1+ \frac{3a'}{a}) )  - \ln
\frac{\vol(M)}{2\pi}\boldsymbol{)}.
\end{eqnarray*}

When $\chi_-(C)$ becomes very large, the dominant expression in
the right hand side of last inequality behaves like $\alpha \sigma
a \chi_-(C) \ln \chi_-(C)$. In fact, an explicit calculation shows
that if $\chi_-(C)\geq c_2$ and $\frac{\vol(M)}{2\pi} \geq 1$,
then
\begin{eqnarray*}
&&\ln  (72\epsilon (a\chi_-(C)+b) ( \frac{\chi_-(C)+4}{2} \boldsymbol{(}
4\alpha^2\sigma^2 (a\chi_-(C)+b)^2\\
&&(\sqrt{2}\alpha\sigma (a\chi_-(C)+b))^{\alpha  \sigma
(a\chi_-(C)+b)} + 5 \boldsymbol{)}  + \frac{2}{9}(1+ \frac{3a'}{a})
) - \ln\frac{\vol(M)}{2\pi})\\ &\leq & 4 \alpha \sigma a \chi_-(C)
\ln \chi_-(C).
\end{eqnarray*}
(see \cite[Chapter 1]{RenThese} for explicit details and calculations).

Thus, if $k_2:=4 \alpha \sigma a$, if $\chi_-(C)\geq c_2$, assuming that
$k_2\,\chi_-(C) \ln \chi_-(C)\leq \ln\ln d/q$ implies that the sufficient
condition above is satisfied, hence conclusions of the Pattern Proposition
D.\qed

~

\noindent \textbf{Claim 3.} \emph{If $c_1\leq \chi_-(C)\leq c_2$, conclusions
of Proposition D still hold if $k_3\,\chi_-(C)\\ \ln \chi_-(C)\leq
\ln\ln \frac{d}{q}$, with
\begin{align*}
k_3&=\frac{1}{c_1\ln c_1}\ln\boldsymbol{(} 36\epsilon(ac_2+b) (c_2+4)
\left( 4(\alpha\sigma)^2(ac_2+b)^2(\sqrt{2}\alpha
\sigma(ac_2+b))^{\alpha\sigma(ac_2+b)} +5 \right)\\
&+ 16 \epsilon (ac_2+b)(1+\frac{3a'}{a})-\ln\frac{\vol(M)}{2\pi} \boldsymbol{)}.
\end{align*}}

~

\noindent \underline{Proof of claim 3.}
\nopagebreak

As $\chi_-(C) \geq c_1$, it is the case where $r\geq 2D+1$, and we
proceed as above, using like during the proof of claim 1
that one has the bounds $c_1\leq \chi_-(C) \leq c_2$.\qed

~

Set $k:=\max\{k_0,k_1,k_2,k_3\}$. It follows from the last
three claims that if $k\, \chi_-(C)\\ \ln \chi_-(C) \leq \ln \ln
\frac{d}{q}$, then conclusions of Proposition D hold.
In particular, $M$ virtually fibers over the circle and the $m$
embedded surfaces in $M'$ are fibers. Furthermore, the
constant $k=k(\epsilon,\vol(M))$ depends only on $\epsilon =
\inj(M)/2$ and the volume $\vol(M)$, and its expression is
explicit. This ends the proof of theorem A.\qed

~

\noindent \underline{Proof of corollary
\ref{cor-thmA-volhandlebodies}.}

If $C_j$ is a handlebody and $\vol(C_j) \geq \vol(M) d/q$, the
proof of theorem A shows that one can construct in $C_j$ surfaces
that are fibers. In particular, the handlebody $C_j$ contains
incompressible surfaces, which is a contradiction.\qed

~

\subsection{Proof of Proposition B: finding nested surfaces.}\label{dem-lm4.5}

~

Let $C$ be the compression body of $M'$ obtained in lemma
\ref{lm-trouver-corps-a-anses}. The boundary of $C$ is a union of
pseudo-minimal surfaces, the genus of each boundary component is
at most $g$, and $\vol(C)\ge\vol(M)\frac{d}{q}$.

\subsubsection{Some modifications of the compression
body.}\label{par-modifications-cea}

~

Instead of the manifold with boundary $C$, we need to work in a
complete Riemannian manifold of sectional curvature at most $-1$.
This is the aim of the following lemma.

\begin{lm}\label{lm-completionC}

Up to modifying the compression body $C$ without significant
changes of volume, one can add collars to boundary components of
$C$ to obtain a (non compact) Riemannian 3-manifold homeomorphic
to the interior of $C$. This manifold is equipped with a complete
metric of sectional curvature at most $-1$, which coincides on $C$
with the induced metric given by the embedding of $C$ in $M'$.

\end{lm}

\noindent \underline{Proof of lemma \ref{lm-completionC}.}
\nopagebreak

We start with the compression body $C$ embedded in $M'$ and its
non complete induced hyperbolic metric. If necessary, we need to
modify slightly the compression body $C$ in order that each
boundary component of $C$ has its intrinsic sectional curvature at
most $-1$.

That is not a problem for boundary components which are minimal
surfaces, as their sectional curvature is always at most this of
the ambient hyperbolic manifold, i.e. $-1$.

If a boundary component of $C$ is the boundary of a small
neighborhood of a non-orientable minimal surface, we can choose
this neighborhood small enough in order that the sectional
curvature of this pseudo minimal surface is bounded from above by
$-1/2$. This is a consequence of the continuity of the intrinsic
sectional curvature in a neighborhood of the minimal surface
(because of the continuity of the Gauss curvature). By rescaling
the metric of the covering $M'$ and the metric of $M$ by a factor
$1/2$ (which multiply sectional curvatures by a factor $2$), we
can suppose that the intrinsic curvature of the boundary
components of $C$ which are the boundary of a small regular
neighborhood of a non-orientable minimal surface is at most $-1$.

If $\partial_+C$ is the boundary of a regular neighborhood $N(S)$
of a non orientable minimal surface $S$ with a small tube attached
vertically in the $I$-bundle structure, we have to consider two
cases. If this tube $\disque\times I$ belongs to the compression
body $C$, we can remove it. More precisely, we compress $C$ along
the disc $\disque\times\{1/2\}$ to get a new compression body of
lower genus. We lose the tube $\disque\times I$ during this
process, but as we can make this tube as small as we like, this
compression does not change significantly the volume of the
compression body. As the positive boundary of this new compression
body $C'$ is then the boundary of a small regular neighborhood of
the minimal non orientable surface $S$, the previous argument
shows that we can suppose that the intrinsic curvature of
$\partial_+C'$ is at most $-1$.

Otherwise, in the second case the tube $\disque\times I$ lies
outside $C$, meaning that $C$ is contained in $N(S)$. We can then
collapse the small tube to an arbitrarily small geodesic arc
$\gamma$ in the regular neighborhood of the minimal non orientable
surface $S$. The positive boundary $\partial_+C$ becomes the union
of the boundary of $N(S)$ and the arc $\gamma$. As before, we can
suppose that the sectional curvature of the surface $\partial
N(S)$ is at most $-1$.

For each boundary component $F$ of $C$, we glue a copy of $F\times
[0,+\infty)$ equipped with a warped product metric. A computation
of the sectional curvature of a warped product (see for example
Bishop and O'Neil \cite[p. 26]{BO'N}) shows that as we start from
a surface $F$ with sectional curvature at most $-1$, there exists
a warped product metric on $S\times [0,+\infty)$ such that this
Riemannian manifold is complete of sectional curvature at most
$-1$. If we are in the last case where $F$ is the boundary of a
regular neighborhood $N(S)$ of a minimal non orientable surface
$S$ with a small tube attached, and this tube is lying outside
$C$, then we forget the arc $\gamma$ for this construction and we
just glue a copy of $\partial N(S)\times [0,+\infty)$ with a
Riemannian metric of curvature at most $-1$. We perturb slightly
this metric to make it smooth, and we obtain thus a complete
Riemannian metric for the interior of $C$ (union $\gamma$ if we
are in this last case) with sectional curvature at most $-1$.\qed

~

The boundary surfaces of $C$ are pseudo
minimal surfaces. This fact is crucial as one can homotop a
minimal surface of genus $g$ to a simplicial surface not too far
away in $C$. This can be done by the following lemmas.

\begin{df}

Let $\epsilon >0$. The (intrinsic) $\epsilon$-diameter of a Riemannian surface
$S$ is the
minimal number of balls of radius $\epsilon$ for the metric of $S$ needed to
cover the surface $S$.

\end{df}

\begin{lm}\label{lm-controle-diam-surf-min}

Suppose $S$ is a pseudo minimal surface in a closed Riemannian 3-manifold
$N$ of sectional curvature at most $-1$. Let $\epsilon\leq \inj(N)$ and
$$a'=6\left(\frac{21}{4}+\frac{3}{4\pi}+ \frac{3}{4\epsilon}+
\frac{2}{\sinh^2(\frac{\epsilon}{4})} \right) .$$

Then the surface $S$ has $\epsilon$-diameter at most
$a'\left|\chi(S)\right|$, and it admits a one-vertex triangulation in which
each edge has length at most $2\epsilon a'\left|\chi(S)\right|$.

\end{lm}

\noindent \underline{Proof of lemma \ref{lm-controle-diam-surf-min}.}
\nopagebreak

This lemma is a direct consequence of \cite[Lemma 4.2 p.
2249]{Mah} and \cite[Proposition 6.1]{Lac} in the case the surface $S$ is
minimal and orientable, and we can take $a'/6$ instead of $a'$.
If $S$ is minimal, but not orientable, its homology class $[S]$ is non zero in
$H_2(N,\zz/2\zz)$. By Poincar\'e's duality, it corresponds to a
non-trivial element $\alpha\in H^1(N,\zz/2\zz)$. As the homology
class of the double cover of $S$ can be represented by the
boundary of a small regular neighborhood of the non-orientable
surface $S$, we have $2[S]=0$ in $H_2(N,\zz)$. If we take the
double cover $N'$ of $N$ corresponding to the kernel of $\alpha$,
the surface $S$ lifts to a minimal orientable surface $S'$. We can
apply the stronger version of lemma
\ref{lm-controle-diam-surf-min}, and bound the $\epsilon$-diameter
of $S'$ by $a'/6\left|\chi(S')\right| =a'/6
\times2\left|\chi(S)\right|= a'/3 \left|\chi(S)\right|$, and the
length of a one-vertex triangulation for $S'$ by $2\epsilon
a'/3\left|\chi(S)\right|$. As those numbers bound also from above
the $\epsilon$-diameter and the length of a one-vertex
triangulation of $S$, this proves the lemma for a minimal non orientable
surface, with $a'/3$ instead of $a'$.

If the surface $S$ is just pseudo minimal, it is the boundary of an arbitrarily
small regular neighborhood of a minimal surface $S'$. As the diameter and the
length of the edges of a one-vertex triangulation are at most $a'/3 \left|
\chi(S')\right|$ and $2\epsilon a'/3 \left| \chi(S') \right|$, with
$\left|\chi(S) \right| \leq 2 \left| \chi(S') \right|$, this ends the proof
of lemma \ref{lm-controle-diam-surf-min}. \qed

~

As from lemma \ref{lm-completionC}, the boundary components of $C$ are pseudo
minimal surfaces, lemma \ref{lm-controle-diam-surf-min} applies to bound from
above the $\epsilon$-diameter and the length of the edges of a one-vertex
triangulation those surfaces. Furthermore, if some geodesic arcs need to be
added, they can be made as small as necessary.

Recall some definitions and results of \cite[Sections
2 et 3]{Mah}.

\begin{df}

A \textbf{coned $n$-simplex} in a compact Riemannian manifold $N$
of sectional curvature at most $-1$ is defined inductively as
follows. A \textbf{coned 1-simplex} $\Delta^1=(v_0,v_1)$ is a
constant speed geodesic from $v_0$ to $v_1$. The speed is allowed
to be zero, and in this case the 1-simplex degenerates to the
point $v_0$. A \textbf{coned $n$-simplex} is a map
$\phi\,:\,\Delta^n\rightarrow N$ such that $\phi_{|\Delta^{n-1}}$
is a coned $(n-1)$-simplex and for all $x\in \Delta^{n-1}$,
$\phi_{|\{tx+(1-t)v_n\, |\,t\in[0,1]\}}$ is a constant speed
geodesic. The map $\phi$ depends on the order of the vertices
$(v_0,\ldots,v_n)$ and its image may not be embedded in $N$, just
immersed.

A \textbf{simplicial surface} is a continuous map
$\phi\,:\,S\rightarrow N$ where $S$ is a triangulated surface,
such that the restriction of the map $\phi$ to each triangle
$\Delta$ of $S$ is a coned 2-simplex.

\end{df}

\begin{lm}\label{lm-surf-min->surf-simpl}

Let $N$ be a complete Riemannian manifold with sectional curvature
at most $-1$. Suppose that $T$ is a connected and orientable
pseudo-minimal surface in $N$ with diameter bounded from above by
$\mathcal{N}$ and admitting a one-vertex triangulation in which
the length of the edges is at most $\mathcal{N}'$. Then $T$ can be
homotoped to a simplicial surface $T'$ with diameter at most
$2\mathcal{N}'$ and such that any point $x\in T$ and $x'\in T'$
are at distance at most $\mathcal{N}+ \mathcal{N}'$ from each
other. Furthermore, every point of $T'$ is at distance at most
$\mathcal{N}'$ from the vertex of the one-vertex triangulation of
$T$.

\end{lm}

\noindent \underline{Proof of lemma
\ref{lm-surf-min->surf-simpl}.}
\nopagebreak

Let $v$ be the vertex of the one-vertex triangulation of $T$.
First, we homotop each edge $e$ of the triangulation of $T$ to its
closed length-minimizing geodesic representative $e'$ in
$\pi_1(N,v)$. If the homotopy class of $e$ is zero (meaning that
the surface $T$ is compressible in $N$), we homotop $e$ to the
degenerate constant speed geodesic $\{v\}$.

Let $\mathcal{T}$ be a triangle in $T$. If all edges of
$\mathcal{T}$ are null-homotopic, $\mathcal{T}'$ is the degenerate
2-simplex corresponding to $\{v\}$. If at least one edge of
$\mathcal{T}$ corresponds to a null-homotopic curve, then we build a
simplicial triangle $\mathcal{T}'$ containing the closed geodesic
at $v$ corresponding to this edge, coned from $v$. More precisely,
the 1-skeleton of $\mathcal{T}'$ is the union of closed geodesics
corresponding to its non homotopically trivial edges. To build the
2-skeleton, we choose one of those non trivial edges and we cone
$v$ to this edge with constant speed geodesics. In this case, each
point in $\mathcal{T}'$ is at distance at most $\mathcal{N}'/2$
from the vertex $v$ (as it is on a closed geodesic of length at
most $\mathcal{N}'$).

If all the edges of $\mathcal{T}$ are non zero in $\pi_1(N,v)$,
they correspond to three non trivial closed geodesics $c_1$, $c_2$
and $c_3$, starting and ending at the point $v$. In the universal
cover $\widetilde{N}$ of $N$, we can choose lifts $a_1$, $a_2$ and
$a_3$ of $c_1$, $c_2$ and $c_3$ that bound a totally geodesic
triangle $\mathbb{T}$. By definition, the covering projection maps
$a_i$ to $c_i$ for $i=1,2,3$. The simplicial triangle
$\mathcal{T}'$ corresponding to $\mathcal{T}$ is the image under
the covering projection of the totally geodesic triangle
$\mathbb{T}$ in $\widetilde{N}$. As the covering projection is an
isometry from the interior of $\mathbb{T}$ to the interior of
$\mathcal{T}'$, and as each point in the interior of $\mathbb{T}$
lies at distance at most $\mathcal{N}'$ (which is an upper bound
for the maximum of the lengths of the sides $a_1$, $a_2$ and
$a_3$), each point $x'$ in the interior of $\mathcal{T}'$ lies at
distance at most $\mathcal{N}'$ from the vertex $v$.

Therefore, starting from the triangulated surface $T$, we can
build a simplicial surface $T'$ such that $v$ is the only vertex
of the simplicial structure of $T'$ and each point $x'$ in $T'$ is
at distance at most $\mathcal{N}'$ from $v$. In particular, the
diameter of $T'$ is at most $2\mathcal{N}'$. As the diameter of
$T$ is at most $\mathcal{N}$ and that $v$ is also a point of $T$,
for any points $x'\in T'$ and $x\in T$, we have:
\begin{eqnarray*}
d(x,x')&\leq & d(x,v)+d(x',v)\\
&\leq &\diam(T)+\mathcal{N}'\\
&\leq & \mathcal{N}+ \mathcal{N}',
\end{eqnarray*}
which proves lemma
\ref{lm-surf-min->surf-simpl}.\qed

~

Given a spine $\Gamma$ for the compression body $C$ which is a
union of simplicial surfaces corresponding to $\partial_-C$ joined
by geodesic arcs, there exists a simplicial surface homotopic to
this spine, by a homotopy that does not sweep out too much volume.
More precisely, this follows from the next lemma, proven in
\cite[Lemma 4.3 p. 2250]{Mah}.

\begin{lm}\cite[Lemma 4.3]{Mah}\label{lm4.3}

Let $\sigma_1,\ldots,\sigma_n$ be a collection of simplicial
surfaces, with basepoints $v_i$ in $N$, a complete Riemannian
3-manifold of sectional curvature at most $-1$. Join the basepoint
$v_1$ to each of the other basepoints by at least one geodesic arc
to obtain a geodesic 2-complex $\Gamma$ homotopic to a surface of
genus $g$. Then, there exists a homotopy of $\Gamma$ to a
simplicial surface $\Sigma_0$ of genus $g$, and this homotopy sweeps
out a volume of at most $3(2g+2)V_3$, where $V_3$ is the
maximal volume of an ideal hyperbolic tetrahedron.\qed

\end{lm}

Recall that $\epsilon\leq \inj(M)/2$, and
$a'=6\left(\frac{21}{4}+\frac{3}{4\pi}+ \frac{3}{4\epsilon}+
\frac{2}{\sinh^2(\epsilon/4)} \right) $. The constant $2\epsilon$
is a uniform lower bound for the injectivity radius of any finite
cover of $M$. In particular, $\inj(M')\geq 2\epsilon$.

~

\noindent\textbf{Lemme \ref{lm-chgts-negligeables}.}
\emph{Let $k_0=\max\left\{ \frac{\ln(4(2\epsilon a'+1))} {2\ln 2},
1+\frac{\ln(1+\ln(12V_3 /\vol(M))}{2\ln 2}\right\}$.}

\emph{If $k_0\,\chi_-(C)\ln \chi_-(C)\leq \ln \ln d/q $ et $\vol(M)\geq
\pi/2$, then applying lemmas \ref{lm-surf-min->surf-simpl} and
\ref{lm4.3} to replace the boundary surfaces of $C$ to simplicial
surfaces, one obtains a new compression body $C''$ with:
$$\vol(C'')\geq \frac{1}{4}\vol(C)\geq \frac{\vol(M) d}{4q}.$$}

~

\noindent \underline{Proof of lemma \ref{lm-chgts-negligeables}.}
\nopagebreak

Let $\partial_-C=T_1\cup\ldots\cup T_n$ be the components of
$\partial_-C$, with $g(T_1)+\ldots +g(T_n)\leq g(\partial_+C)$. As
in lemma \ref{lm-surf-min->surf-simpl}, replace $\partial_+C=:S_0$
and $\partial_-C=S_1\cup\ldots\cup S_n$ by simplicial surfaces
$S_0'$ and $T'_1\cup\ldots\cup T'_n$, close to the previous
surfaces. If $v_j\in T_j$ is the vertex of the one-vertex
triangulation of $T_j$, then lemmas \ref{lm-controle-diam-surf-min} and
\ref{lm-surf-min->surf-simpl}
show that every point of $T_j'$ lies at distance at most
$\mathcal{N}'=2\epsilon a'\left|\chi(T_j)\right| \leq 2\epsilon a'
\chi_-(C)$ from $v_j$. Thus, each new surface $T_j'$ is contained
in the ball of center $v_j$ and radius $2\epsilon a' \chi_-(C)$.
If $C'$ is the new compression body obtained by replacing the
surfaces $T_j$ by the surfaces $T_j'$, the modification of volume
is at most
\begin{eqnarray*}
\vol(C')&\geq & \vol(C)-\sum_{j=0}^n\vol(B(v_j,2\epsilon a' \chi_-(C)))\\
&\geq &\vol(C)-(g(C)+1) \vol(B_{\H3}(2\epsilon a' \chi_-(C)))\\
&\geq &\vol(C)\left(1- \frac{\pi(\chi_-(C)+4)(\sinh(4\epsilon a' \chi_-(C))
-4\epsilon a' \chi_-(C))}{2\vol(M)d/q} \right)
\end{eqnarray*}

Let us show that $\vol(C')\geq \vol(C)/2$, which is the same as
proving that
$$\frac{\pi(\chi_-(C)+4)(\sinh(4\epsilon a' \chi_-(C))-4\epsilon a'
\chi_-(C))}{2\vol(M)d/q} \leq \frac{1}{2}.$$ It suffices to prove
that $\ln \frac{\pi(\chi_-(C)+4)(\sinh(4\epsilon a'
\chi_-(C))-4\epsilon a' \chi_-(C))}{\vol(M)d/q} \leq 0$. But

\begin{eqnarray*}
\ln \left( \frac{\pi}{\vol(M)d/q} (\chi_-(C)+4)(\sinh(4\epsilon
a'\chi_-(C))
-4\epsilon a'\chi_-(C)) \right) &\leq &\\
\ln \left(  \frac{\pi}{2\vol(M)} \frac{(\chi_-(C)+4) \exp(4\epsilon
a'\chi_-(C))} {d/q} \right) &\leq &\\
\ln \left(  \frac{\pi}{2\vol(M)}\right)+ \ln \left( (\chi_-(C)+4)\exp(4\epsilon
a'\chi_-(C)) \right) -\ln(d/q) & \leq &\\
\ln (\chi_-(C)+4) + 4\epsilon a' \chi_-(C)&-&\ln(d/q),
\end{eqnarray*}
as by assumption, $\vol(M)\geq \pi/2$.

As for every $x\geq2$, $\ln(x+4)\leq 2x$, it suffices to prove that
$$ (2+4\epsilon a') \chi_-(C) \leq \ln(d/q),$$ which is the same as
$$\ln (2+4\epsilon a') + \ln \chi_-(C) \leq \ln \ln(d/q)
.$$ Now by assumption, $\frac{\ln\ln d/q}{\chi_-(C)\ln
\chi_-(C)}\geq k_0 \geq \frac{\ln(4(2\epsilon a'+1))}{2\ln 2}$. Thus,
\begin{eqnarray*}
\ln \ln d/q &\geq & \frac{\ln(4(2\epsilon a'+1))}{2\ln 2} \chi_-(C) \ln
\chi_-(C)\\
& \geq & \frac{\ln 2 + \ln (2+4\epsilon a')}{2\ln 2} \chi_-(C) \ln \chi_-(C)\\
& \geq & \frac{\chi_-(C)\ln \chi_-(C)}{2} + \frac{\ln (2+4\epsilon a')\chi_-(C)
\ln \chi_-(C)}{2\ln 2}\\
& \geq & \ln \chi_-(C) + \ln (2+4\epsilon a')
\end{eqnarray*}
as $\chi_-(C) \geq 2$, showing that $\vol(C')\geq
\vol(C)/2$.

~

From lemma \ref{lm4.3}, the volume swept out by the homotopy
between $\Gamma$ and $\Sigma_0$ is at most $3(\chi_-(C)+4)V_3$. As
the volume of $C$ is at least $\vol(M)d/q$ by lemma
\ref{lm-trouver-corps-a-anses}, the volume of what remains after
cutting the metric completion of $C'$ along $\Sigma_0$ and
throwing off components containing the infinite products to obtain
a new compression body $C''$ is at least
$$\vol(C')-3(\chi_-(C)+4)V_3
\geq\vol(C)/2\left(1-3V_3 (\chi_-(C)+4)
\frac{2q}{\vol(M)d}\right).$$

Therefore, it suffices to prove that $3V_3 (\chi_-(C)+4)
\frac{2q}{\vol(M)d} \leq \frac{1}{2}$, or $\ln (\frac{12V_3
(\chi_-(C)+4) }{\vol(M)d/q })\leq 0$, or also
$$\ln\left(\ln \frac{12 V_3}{\vol(M)} + \ln(\chi_-(C)+4)\right) \leq \ln
\ln(d/q).$$

As $\chi_-(C)\geq 2$, $\ln(\chi_-(C)+4)\geq \ln 6 >1$. Thus,
\begin{eqnarray*}
\ln\left( \ln \frac{12 V_3}{\vol(M)} + \ln(\chi_-(C)+4) \right)&=&\\
\ln\left( \ln(\chi_-(C)+4) (1+ \frac{\ln \frac{12 V_3}{\vol(M)}}
{\ln(\chi_-(C)+4)})  \right) &=&\\
\ln \ln (\chi_-(C)+4) + \ln \left(1+ \frac{\ln \frac{12 V_3}{\vol(M)}}
{\ln(\chi_-(C)+4)}\right) &\leq &
\ln \ln (\chi_-(C)+4) + \ln \left(1+ \ln \frac{12 V_3}{\vol(M)} \right).
\end{eqnarray*}

As soon as $c\geq 2$, $\frac{\ln\ln (c+4)}{c\ln c}\leq 1$. Then,
\begin{eqnarray*}
\frac{\ln\left(\ln \frac{12 V_3}{\vol(M)} + \ln(\chi_-(C)+4)\right)} {\chi_-(C)
\ln \chi_-(C)} &\leq &
\frac{\ln \ln (\chi_-(C)+4)}{\chi_-(C) \ln \chi_-(C)} + \frac{\ln \left(1+ \ln
\frac{12 V_3}{\vol(M)} \right)}{\chi_-(C) \ln \chi_-(C)}\\
& \leq & 1+ \frac{\ln \left(1+ \ln \frac{12 V_3}{\vol(M)} \right)}{2\ln 2}.
\end{eqnarray*}

Now, as $\frac{\ln \ln d/q}{\chi_-(C) \ln \chi_-(C)}\geq k_0 \geq
1+\frac{\ln(1+\ln(12V_3 /\vol(M))}{2\ln 2}$,
$$\frac{\ln \ln d/q}{\chi_-(C) \ln \chi_-(C)}\geq \frac{\ln\left(\ln \frac{12
V_3}{\vol(M)} + \ln(\chi_-(C)+4)\right)} {\chi_-(C) \ln \chi_-(C)},$$
and so $\vol(C'')\geq \vol(C)/4 \geq\vol(M)\frac{d}{4q}$, which
ends the proof of lemma \ref{lm-chgts-negligeables}.\qed

~

In the sequel, to simplify notations, we will still denote by $C$
the new compression body $C''$ and we work in the closure of the
region of $C$ bounded by the two connected simplicial surfaces
$\Sigma_0$ (corresponding to $\partial_-C$ union some arcs, and
forming a spine for $C$), and $\Sigma_1$ corresponding to
$\partial_+C$.

~

\subsubsection{Sweepouts.}

\begin{df}

Let $C$ be a compression body. Set $S=\partial_+C$. A
\textbf{sweepout} of the compression body $C$ is a 1-parameter
family of surfaces $\{S_t\}_{t\in[0,1]}$ such that $S_0$ is a
spine of $C$, $S_1=S=\partial_+C$, for all $t\in(0,1]$ the surface
$S_t$ is homeomorphic to $S$, and the application
$\Phi\,:\,S\times I\rightarrow C$ is of homological degree one.

\end{df}

There exists a sweepout $\{S_t\}_{t\in[0,1]}$ of the compression
body $C$ such that $S_0=\Sigma_0$ and $S_1=\Sigma_1$. The sweepout
surfaces $S_t$ for $t>0$ are of interest in order to construct a
long product in the compression body $C$. But, if we can control
the diameter of a minimal surface in terms of its genus and the
injectivity radius of the ambient manifold, we cannot control
uniformly the diameter of all the sweepout surfaces $S_t$: there
may appear some long and thin Margulis tubes, containing a closed
geodesic of the surface with length less than the injectivity
radius of $M'$.

To face this problem, we work with the notion of
$\epsilon$-diameter, for which non-connected surfaces with small
diameter components are considered as "small". Recall the definition.

\begin{df}

Let $\epsilon>0$. The \textbf{(intrinsic) $\epsilon$-diameter} of a
non-necessarily connected surface $F$ is the minimal number of
balls of radius $\epsilon$ for the metric of $F$ required to cover
the surface $F$.

If $F$ is immersed in a Riemannian 3-manifold $N$, the
\textbf{$\epsilon$-diameter} of $F$ in $N$ is the minimal number of 3-balls in
$N$ of radius $\epsilon$ for the metric of $N$ required to cover $F$.

\end{df}

\begin{rmq}

If $F$ is immersed in a Riemannian 3-manifold $N$, the $\epsilon$-diameter of
$F$ in $N$ is always at most the intrinsic $\epsilon$-diameter of $F$ with
respect to the induced metric.

\end{rmq}

At this point, we recall the technique of Maher to construct from
the original sweepout $\{S_t\}_{t\in I}$ of $C$ what he calls a
"generalized sweepout" $\{\widehat{S}_t\}_{t\in I}$ in which the
$\epsilon$-diameter of good sweepout surfaces is uniformly bounded
from above (see \cite[Sections 2 and 3]{Mah}).

The first step is to straighten the sweepout $\{S_t\}_{t\in I}$ to
a simplicial sweepout, using results of Bachman, Cooper and White
\cite{BCW}. We recall terminology and results stated in
\cite[Sections 2 and 3]{Mah}.

\begin{df}

A \textbf{simplicial sweepout} is a sweepout $\Phi\,:\,S\times
I\rightarrow N$ such that each surface $S_t$ is mapped to a
simplicial surface with at most $4g(S)$ triangles, and at most one
vertex of angle sum less than $2\pi$.

\end{df}

The following lemma ensures that we can homotop the sweepout
$\{S_t\}_{t\in[0,1]}$ between $\Sigma_0$ and $\Sigma_1$ to a
simplicial sweepout. It is an improvement of \cite[Theorem
2.3]{BCW}, and is proven by Maher \cite[Lemma 2.5 p. 2236]{Mah}.

\begin{lm}\cite[Lemma 2.5]{Mah}\label{lm2.5}

Let $N$ be a closed orientable Riemannian manifold of sectional
curvature at most $-1$. If $\Sigma_0$ and $\Sigma_1$ are
simplicial surfaces with one-vertex triangulations, which are
homotopic by a homotopy $\Phi\,:\,S\times I\rightarrow N$, then
there exists a simplicial sweepout $\Phi'\,:\,S\times I\rightarrow
N$ homotopic to $\Phi$ relative to $S\times \partial I$.\qed

\end{lm}

Therefore, we can suppose that the sweepout in the compression
body $C$ is simplicial between the simplicial surfaces
$\Sigma_0=S_0$ and $\Sigma_1=S_1$.

After getting this simplicial sweepout in the compression body
$C$, the next step will be to get rid of the long and thin tubes
in the sweepout surfaces to get a "generalized sweepout" in which
the $\epsilon$-diameter of all sweepout surfaces is uniformly
bounded from above.

\begin{df}\cite[Definition 3.2 p. 2237]{Mah}

Let $N$ be a compact, connected and oriented 3-manifold. A
\textbf{generalized sweepout} of $N$ is given by a triple
$(\Sigma,f,h)$, where $\Sigma$ is an orientable and compact
3-manifold, the map $h\,:\,\Sigma\rightarrow \rr$ is a Morse
function, constant on each boundary component of $\Sigma$ and such
that for all but finitely many $t\in\rr$, the set $f^{-1}(\{t\})$
is an immersed surface. Moreover, it is required that
$f\,:\,(\Sigma,\partial\Sigma)\rightarrow (N,\partial N)$ is of
homological degree one.

\end{df}

Of course, an ordinary sweepout $\Phi\,:\,S\times I\rightarrow N$
is an example of generalized sweepout: the Morse function $h\,:\,S\times
I\rightarrow\rr$ is given by the projection to the factor $I$, and
for all $t\in(0,1)$, $h^{-1}(\{t\})=S_t$ is an immersed surface in $N$. By
definition of a sweepout, $\Phi\,:\,(S\times I, S\times\partial
I)\rightarrow (N,\partial N)$ is of homological degree one.

For all $x\in\Sigma$, we think of $h(x)=t$ as the time coordinate.
A generalized sweepout can be seen as a one-parameter family of
immersed surfaces $S_t$ with singular times $t$ where the genus or
the number of components of those surfaces change.

Starting from the simplicial sweepout $\{S_t\}_{t\in I}$ of $C$,
we wish to obtain a generalized sweepout in which each sweepout
surface has bounded $\epsilon$-diameter. To this aim, we follow
Maher and introduce the notion of \textbf{surgery} of a
generalized sweepout.

\begin{df}

One can obtain from a generalized sweepout given by $(\Sigma,f,h)$
a new generalized sweepout $(\Sigma',f',h')$ by an operation
called a \textbf{surgery of generalized sweepouts}, as described
below. (In fact, it is a special case of a more general
construction called a modification of generalized sweepout,
described by Maher and Rubinstein in \cite{MR}.)

\end{df}

Let $(\Sigma,f,h)$ be a generalized sweepout of a 3-manifold $N$.
Take a submanifold in $\Sigma$ of the form $A\times[a,b]$ where
$0<a<b<1$ and $A$ is an annulus in the surfaces $S_t$ for
$t\in[a,b]$. We do $(0,1)$ surgery to this solid torus $A\times
[a,b]$ in the following way: choose two times $c$ and $d$ such
that $a<c<d<b$. Take a chore geodesic $\gamma$ for the annulus $A$
in the surface $S_a$. Shrink this geodesic : it gets shorter and
shorter, until it collapses to a point in a modification $S'_c$ of
the surface $S_c$. For all $t\in(c,d)$, replace the surface $S_t$
by the surface $S_t'$ obtained from $S_t$ by surgering along
$\gamma$, i.e. we cut $S_t$ along $\gamma$ and cap off the
resulting surface with two discs. Do this in a smooth way, such
that the two discs of $S'_t$ get closer and shrink to a single
point at time $d$. The new surface $S'_d$ is then singular, with a
singular point corresponding to the two former discs. This point
becomes again the geodesic $\gamma$ that gets larger when
$t\in(d,b]$ increases. Do this in such a way that you do not
modify $S_a$ nor $S_b$ nor $\partial A\times[a,b]$. In this way,
we get a new generalized sweepout $(\Sigma',f',h')$, where
$\Sigma'$ is obtained by replacing $A\times [a,b]\subset \Sigma$
by the new manifold where $S_t$ is replaced by $S_t'$ for all
$t\in[a,b]$. Let us denote by $T$ the small tube in $N$ bounded by
$A$, where the surgeries take place. The new maps $(f',h')$
coincide with $(f,h)$ outside $T\times[a,b]$ and in $\partial
(T\times [a,b])$. As the modification of the sweepout takes place
in a proper compact subset of $N$, there exists a point $x$ in the
interior of $N\setminus (T\times[a,b])$. As the map $f$ is not
modified in a neighborhood of $f^{-1}(\{x\})$, the homological
degree of $f'$ is the same as the homological degree of $f$, so it
is still equal to one. Thus the triple $(\Sigma',f',h')$ is still
a generalized sweepout.

\begin{center}

\includegraphics[width=1\textwidth]{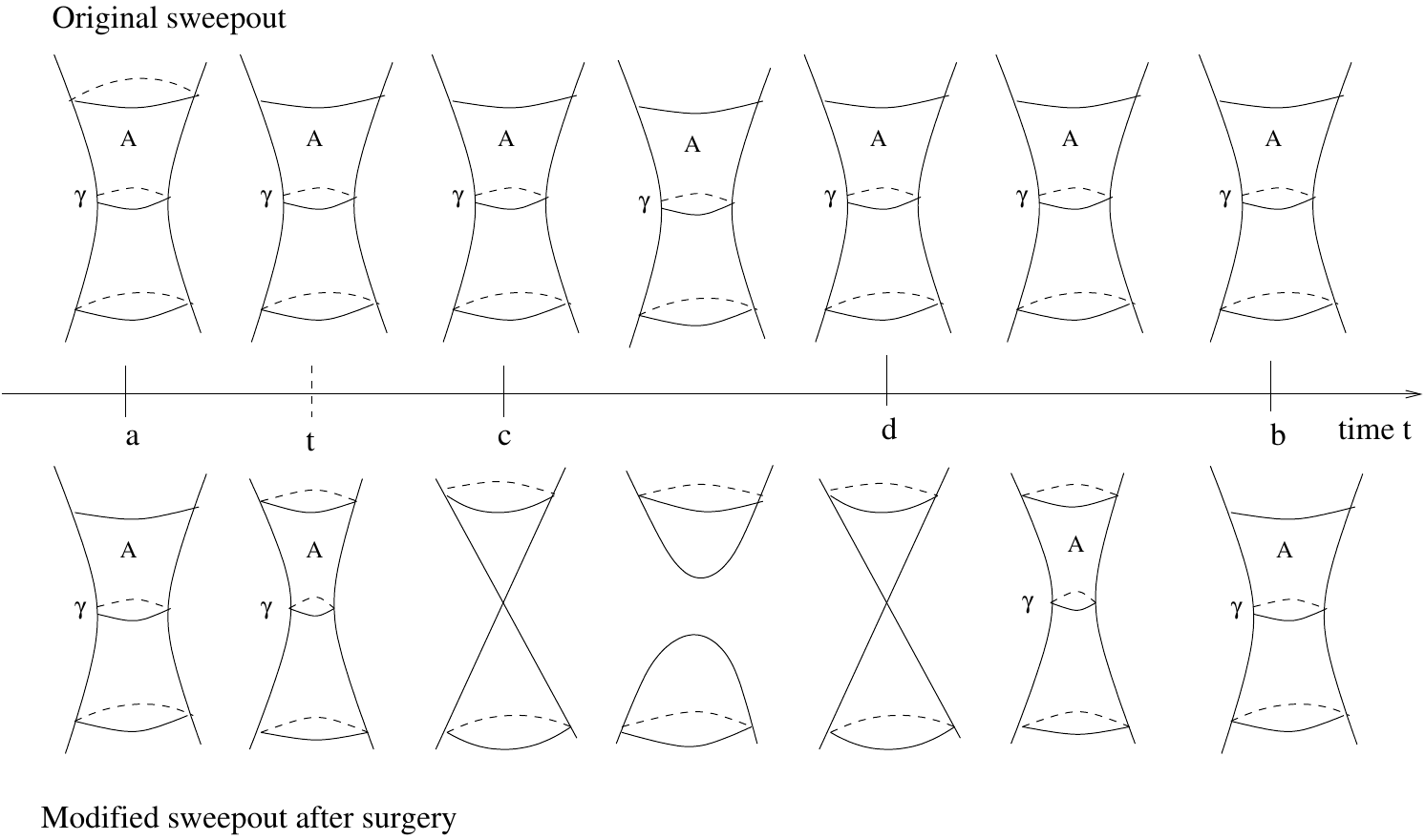}

\end{center}

Set $\boldsymbol{\mathcal{N}_+=\mathcal{N}(\partial_+C)}=\{x\in
C\,,\,d(x,\Sigma_1) \leq
\epsilon/2\}$.

Let $K= 4 \left(3+ 1/\sinh^2{(\epsilon/8)}\right) g(C)-10$ and
$K':=2a'\left|\chi(\partial_+C)\right|$.

\begin{prop}\label{prop-chirurgies-balayage}

Let $\mu>0$. There exists a constant $\eta>0$ as small as wanted,
depending only on the simplicial sweepout $\{S_t\}_{t\in I}$ and
$\mu$, and a finite sequence of surgeries of the simplicial
sweepout giving a generalized sweepout $\{\widehat{S}_t \}_{t\in
I}$ of $C$ and satisfying the following properties.

For every regular time $t\in [\eta,1-\eta]$, the intrinsic $\epsilon$-diameter
of every component of $\widehat{S}_t$ disjoint from $\mathcal{N}_+$ is
less than or equal to $K$. In every case, the diameter of any connected
component of $\widehat{S}_t$ in the compression body $C$ is at
most $\epsilon(1+2K'+2K)$. For $t\geq 1-\eta$, each point on the surface
$\widehat{S_t}$ lies at distance at most $\epsilon K'$ from
$\Sigma_1$. For $t\leq\eta$, any point on one of the original surfaces $S_t$ is
at distance at most $\mu/2$ from $\Sigma_0$. Furthermore, for each
regular time $t$, the surface $\widehat{S}_t$ is homotopic to an embedded
surface obtained from $\partial_+C$ by surgeries.

\end{prop}

\noindent \underline{Proof of proposition
\ref{prop-chirurgies-balayage}.}
\nopagebreak

In order to prove this proposition, the general idea is to cut the
simplicial surfaces $S_t$ along curves that are too short, namely
of length less than or equal to $\epsilon$, and to replace them by
ruled discs, to get rid of long and thin tubes. This is described
by Maher in the third section of \cite[p. 2238 to p. 2245]{Mah}.
We recall here the proof, and we bring some precisions when they
appear to be necessary.

Let $t$ be a regular time. The simplicial surface $S_t$ is
composed of ruled triangles with at most one vertex of angle sum
less than $2\pi$, denoted by $v_t$. Let $\overline{S}_t$ be the
completion of the universal cover $\widetilde{S_t}$ of $S_t\setminus
\{v_t\}$. As it is a metric 2-complex composed of triangles of
curvature at most $-1$ and with vertices whose cone angles are all
greater than or equal to $2\pi$, $\overline{S}_t$ is a complete
CAT($-1$) geodesic metric space. Those spaces satisfy some useful
properties, see \cite{BH} and \cite[p. 2239]{Mah}.

Let $\alpha$ be a homotopy class in $S_t\setminus \{v_t\}$. To
$\alpha$, we can associate a covering transformation of the
universal cover of $S_t\setminus \{v_t\}$, which can be extended
to an isometry of $\overline{S}_t$. As the completion of a
fundamental domain for $S_t\setminus \{v_t\}$ is compact, this
isometry cannot be parabolic. Thus it is hyperbolic or elliptic.
Let $\overline{\gamma}_t$ be the set of points in $\overline{S}_t$
which are moved the least distance by the isometry. This is a
geodesic if the isometry is hyperbolic, or isolated points if the
isometry is elliptic. We denote by $\gamma_t$ the projection of
$\overline{\gamma}_t$ under the covering map, in the sense that if
$\overline{\gamma}_t$ is a geodesic and does not meet
$\overline{S}_t\setminus \widetilde{S_t}$, $\gamma_t$ is a closed
piecewise geodesic homotopic to $\alpha$ in $S_t\setminus
\{v_t\}$. If $\overline{\gamma}_t$ is a geodesic meeting
$\overline{S}_t\setminus \widetilde{S_t}$, then we perturb it
slightly and in an equivariant way such that it avoids
$\overline{S}_t\setminus \widetilde{S_t}$ and its projection
$\gamma_t$ in $S_t\setminus \{v_t\}$ is an embedded closed curve
in the homotopy class of $\alpha$. Finally, if
$\overline{\gamma}_t$ is a set of points, it corresponds to the
constant loop $\gamma_t$ of length zero and equal to the point
$v_t$. By extension, in any case we will call $\gamma_t$
\textbf{the geodesic representative of $\alpha$}. Notice that
$\gamma_t$ is an embedded curve or a point in $C$.

As the negatively curved triangles that compose the surfaces $S_t$
vary continuously with the time $t$, we can expect the geodesic
representatives $\gamma_t$ to vary also continuously. This is
proven by Maher \cite[Lemma 3.4 p. 2240]{Mah}.

\begin{lm}\cite[Lemma 3.4]{Mah}

Let $\gamma$ be a simple closed curve in $S\setminus \{v\}$ where $v$
is a point of $S$ mapping to the point $v_t$ for each time $t$.
Then the geodesic representatives $\gamma_t$ of $\gamma$ vary
continuously with $t$.\qed

\end{lm}

\begin{df}

A geodesic representative $\gamma_t$ is said to be \textbf{short}
if its length is less than or equal to $\epsilon$.

\end{df}

For all $t$, let $\Gamma_t$ be the set of short geodesic
representatives of $S_t$. This is a finite set and it is not
empty, as the geodesic representative of the homotopy class of the
loop around $v_t$ has length zero.

Let $\gamma_t$ be a short geodesic representative. Pick up a
connected component $\widetilde{\gamma_t}$ of
$\overline{\gamma}_t$, the preimage of $\gamma_t$ in
$\overline{S}_t$. Choose an orientation for $\widetilde{\gamma_t}$
so that the distance function from $\widetilde{\gamma_t}$ has a
well defined sign. In the special case where the length of
$\gamma_t$ is zero, the distance from $\widetilde{\gamma_t}$ will
always be non negative. If $[p,q]$ is an interval of $\mathbb{R}$,
let $\widetilde{N}_{[p,q]}(\widetilde{\gamma_t})$ be the set of
points $x\in\overline{S}_t$ such that $p\leq
d(x,\widetilde{\gamma_t})\leq q$. Let $N_{[p,q]}(\gamma_t)$ be the
image in $S_t$ of $\widetilde{N}_{[p,q]}(\widetilde{\gamma_t})$
under the covering projection. If the interval is the single point
$\{r\}$, we will denote this neighborhood by $N_{[r]}(\gamma_t)$.

\begin{df}

Let $\mathcal{A}(\gamma_t)$ be the maximal neighborhood
$N_{[p,q]}(\gamma_t)$ such that for every length $r\in[p,q]$,
$N_{[r]}(\gamma_t)$ is an embedded simple curve of length at most
$\epsilon$. The set $\mathcal{A}(\gamma_t)$ is called the
\textbf{annular neighborhood} of $\gamma_t$.

Define
$\mathcal{E}(\gamma_t)=N_{[p+\epsilon/2,q-\epsilon/2]}(\gamma_t)$
to be the \textbf{surgery neighborhood} corresponding to
$\gamma_t$, with the convention that $\mathcal{E}(\gamma_t)$ is
the empty set if $q-p<\epsilon$. This neighborhood is the subset
of $\mathcal{A}(\gamma_t)$ corresponding to the union of all
curves $N_{[r]}(\gamma_t)$ lying at distance at least
$\frac{\epsilon}{2}$ from the boundary of $\mathcal{A}(\gamma_t)$.

\end{df}

As the annular neighborhood $\mathcal{A}(\gamma_t)$ contains
$\gamma_t=N_{[0]}(\gamma_t)$, it is not empty. The annular
neighborhood and the surgery neighborhood vary continuously with
$t$, but the surgery neighborhood $\mathcal{E}(\gamma_t)$ can be
empty, and it does not necessarily contain the geodesic
representative $\gamma_t$.

The following lemma is proven by Maher in \cite[Lemma 3.7 p. 2242]{Mah}.

\begin{lm}\cite[Lemma 3.7]{Mah}\label{lm3.7}

If $\alpha_t$ and $\beta_t$ are short geodesic representative of
distinct homotopy classes in $S_t\setminus v_t$, then their
surgery neighborhoods $\mathcal{E}(\alpha_t)$ and
$\mathcal{E}(\beta_t)$ are disjoint.\qed

\end{lm}

We notice that this lemma implies that for each time $t$, there
are at most $2g-1$ surgery neighborhoods, where $g$ is the genus
of the sweepout surface $S_t$.

Lemma \ref{lm3.7} allows us to do surgeries on the sweepout
surfaces $S_t$ to get a generalized sweepout where the diameter of
the thin tubes can be controlled. More precisely, the idea is to
remove the surgery neighborhoods from the sweepout surfaces each time it is
possible, to get a new generalized sweepout $(\widehat{S_t})_{t\in I}$. We
describe this construction in details.

Let $\mathcal{E}(\gamma_t)$ be a surgery neighborhood, and $[a,b]$
a maximal time interval on which $\mathcal{E}(\gamma_t)$ is not
empty. We have $0\leq a \leq b\leq 1$. First, suppose that
$0<a\leq b<1$, i.e. that $[a,b]$ is contained in the interior of
$I$. Then $\mathcal{E}(\gamma_a)$ and  $\mathcal{E}(\gamma_b)$ are
two embedded simple curves $N_{[r_a]}(\gamma_a)$ and
$N_{[r_b]}(\gamma_b)$ for lengths $r_a$ and $r_b$, and the union of the surgery
neighborhoods
$\mathcal{E}(\gamma_{[a,b]})=
\{\mathcal{E}(\gamma_t)\,,\,t\in[a,b] \}$ is a solid torus in
$\Sigma$, on which we wish to do a surgery of generalized
sweepouts. When it is possible, we follow Maher's construction.

There is a difficulty here, as the surgery of generalized sweepouts described
above is possible only if the geodesic $\gamma_t$ bounds an immersed disc in
the compression body $C$. Therefore, we need to make a distinction between two
cases of surgery neighborhoods.

\begin{df}

A \textbf{persistent surgery neighborhood} of $S_t$ is a surgery
neighborhood $\mathcal{E}(\gamma_t)$ for which the corresponding geodesic
$\gamma_t$ is not homotopically trivial in $C$.

\end{df}

\begin{lm}\label{lm-chirurgies_pas_possibles}

Let $t\in(0,1)$ and $\mathcal{E}(\gamma_t)$ be a persistent surgery
neighborhood for $S_t$.

If the corresponding surgery curve $\gamma_t$ is homotopically trivial in $C$,
then $\mathcal{E}(\gamma_t)$ is entirely contained in $\mathcal{N}_+$.

In particular, for every points $x$ and $y$ in the union of the persistent
surgery neighborhoods of $S_t$, the distance in $C$ between $x$ and $y$ is at
most $\epsilon+\diam(\Sigma_1)\leq \epsilon(1+2K')$.

\end{lm}

\noindent \underline{Proof of lemma \ref{lm-chirurgies_pas_possibles}.}

\nopagebreak

For each $r\in[p,q]$, as the curve $N_{[r]}(\gamma_t)$ is of length at most
$\epsilon < \inj(M')$, it is null-homotopic in $M'$ and is contained in a
hyperbolic 3-ball $B$, isometrically embedded in $M'$ and of
diameter $\epsilon/2$.

The curve $\gamma_t$ is not homotopically trivial in $C$ and
$N_{[r]}(\gamma_t)$ is homotopic to $\gamma_t$, so $B\cap
\partial C \neq \emptyset$. By assumption $(1)$ of theorem A, the
negative boundary $\partial_- C$ is a union of incompressible surfaces.
Necessarily, $B\cap \partial_+ C
\neq \emptyset$, which, with the simplicial surfaces, writes $B\cap \Sigma_1
\neq \emptyset$. Thus, as the curve $N_{[r]}(\gamma_t)$ is contained in $B$
and that $B$ intersects $\Sigma_1$, each point of $N_{[r]}(\gamma_t)$ is at
distance at most $\epsilon/2$ from $\Sigma_1$. It follows that the surgery
neighborhood $\mathcal{E}(\gamma_t)$ is entirely contained in $\mathcal{N}_+$.

Noticing that the diameter of $\mathcal{N}_+$ in the manifold $C$ is at
most $\epsilon+\diam(\Sigma_1)\leq \epsilon(1+2K')$ finishes to prove lemma
\ref{lm-chirurgies_pas_possibles}.\qed

\begin{rmq}\label{rmq-vois-persistant}

Let $\mathcal{E}(\gamma_{t_0})$ be a surgery neighborhood for a given time
$t_0\in(0,1)$ and $[a,b]$ a maximal time interval on which
$\mathcal{E}(\gamma_t)$ is not empty. If one of the surgery neighborhoods
$\mathcal{E}(\gamma_{t_1})$ is persistent, the all surgery neighborhoods
$\mathcal{E}(\gamma_t)$ are persistent for every time $t$ in $[a,b]$.

\end{rmq}

Indeed, for $t\in[a,b]$, the curves $\gamma_t$ are homotopic, so if one of them
is homotopically non trivial in $C$, it is the case for each of them.\qed

~

We can now carry on with the proof of proposition
\ref{prop-chirurgies-balayage}. Suppose that $\mathcal{E}(\gamma_t)$ is a non
persistent surgery neighborhood. We describe the operation of surgery on
$\mathcal{E}(\gamma_t)$, following Maher \cite[p. 2242 and 2243]{Mah}.

Choose a continuous family of basepoints on the boundary of
$\mathcal{E}(\gamma_t)$ such that the two basepoints agree at times $a$
and $b$. We modify the sweepout by expanding times $a$ and $b$ to
short intervals $I_a$ and $I_b$ on which the map is constant for
the moment.

On the interval $I_a$, the curve $N_{[r_a]}(\gamma_a)$ is an
embedded simple curve homotopic to the geodesic $\gamma_t$. As by assumption,
$\gamma_t$ is null-homotopic in $C$, the curve $N_{[r_a]}(\gamma_a)$ bounds an
immersed disc in $C$.

In the interval $I_a$, we replace
in a continuous way the curve
$N_{[r_a]}(\gamma_a)=\mathcal{E}(\gamma_a)$ by a pair of ruled
discs in $C$ coned from the basepoint $x_a$. More precisely, as
the metric of $C$ is complete, the ruled disc is the union of all
minimizing geodesics between each point of $N_{[r_a]}(\gamma_a)$
and the basepoint $x_a$. Its curvature is then at most $-1$. In
the interval $(a,b)$, we remove the surgery neighborhood
$\mathcal{E}(\gamma_t)$ and we replace it by a pair of such ruled
discs in $C$ coned from the basepoints of the boundary of the
surgery neighborhood. Finally, in the interval $I_b$ we paste the
discs together to come back to the original surface. This is a
surgery of a generalized sweepout as defined above.

The following lemma directly follows from Lemmas 3.8 to 3.10 and
is proven in \cite[p. 2243 to 2246]{Mah}.

\begin{lm}\label{lm-diam-surf-chirurgis\'ees}

Suppose that all surgery neighborhoods $\mathcal{E}(\gamma_t)$ of $S_t$ can be
replaced by pairs of ruled discs as described above and let $\widehat{S_t}$ be
the resulting surface. Then the $\epsilon$-diameter of $\widehat{S}_t$
is at most $K=4\left(3+1/\sinh^2{(\epsilon/8)}\right)
g(C)-10$.\qed

\end{lm}

The construction is now the following. Let $t\in [\eta,1-\eta]$. From remark
\ref{rmq-vois-persistant}, the fact that a given surgery neighborhood
$\mathcal{E}(\gamma_t)$ of $S_t$ is persistent or not depends only on the
maximal time interval $[a,b]$ on which it exists. If it is not persistent, then
apply the surgery procedure described above. If it is persistent, leave it
unchanged. Let $\widehat{S_t}$ be the new generalized sweepout surface obtained.

If none of the surgery neighborhoods $\mathcal{E}(\gamma_t)$ of $S_t$ are
persistent, they have been removed by the surgery procedure. From lemma
\ref{lm-diam-surf-chirurgis\'ees}, the intrinsic $\epsilon$-diameter of
$\widetilde{S}_t$ is at most $K$.

Otherwise, lemma \ref{lm-chirurgies_pas_possibles} ensures that the diameter of
the union of all persistent surgery neighborhoods in $C$ is at most
$\epsilon(1+2K')$, as they are contained in $\mathcal{N}_+$. As the intrinsic
$\epsilon$-diameter of each component of $\widehat{S_t}$ cut along persistent
surgery neighborhoods is at most $K$ from lemma
\ref{lm-diam-surf-chirurgis\'ees}, the diameter in the compression body $C$ of
each connected component of $\widehat{S_t}$ is at most $\epsilon(1+2K'+2K)$.

Furthermore, if a component of $\widehat{S_t}$ does not intersect
$\mathcal{N}_+$, it does not contain any persistent surgery neighborhood, and
its
intrinsic $\epsilon$-diameter is at most $K$ from lemma
\ref{lm-diam-surf-chirurgis\'ees}.

~

An other difficulty is that Maher's construction does not take the boundaries of
the time
interval $I$ into account. However, it may happen that $a=0$ or
$b=1$, and in this case we might be obliged to modify the starting
and finishing simplicial sweepout surfaces $S_0=\Sigma_0$ and
$S_1=\Sigma_1$, which we want to avoid. Therefore, if this case
occurs, we need to refine the construction to modify the
simplicial sweepout in a small regular neighborhood of $S_0\cup
S_1$ in such a way that we do not modify the surfaces $S_0$ and
$S_1$. As we will lose control on the diameter of the
sweepout surfaces in this regular neighborhood, we have to choose
it small enough in order that the sweepout surfaces we will pick
up later to be some of the nested surfaces are not in this
neighborhood. Thus we can control their diameter well. The constant
$\mu$ has been introduced in assumptions of proposition
\ref{prop-chirurgies-balayage} in order to take care of that, and
its value will be defined later.

To finish to modify the original simplicial sweepout to get the
desired generalized sweepout, there remains to consider the case
when $a=0$ or $b=1$. If $\mathcal{E}(\gamma_0)$ is a non persistent surgery
neighborhood and just a single closed
curve, we can apply the previous construction, replacing the time
$0$ by an interval $I_0$ and doing surgery on this interval,
without modifying the starting boundary surface $S_0=\Sigma_0$. It
works similarly if $\mathcal{E}(\gamma_1)$ is a single closed
curve. The problem is when $\mathcal{E}(\gamma_0)$ or
$\mathcal{E}(\gamma_1)$ have non empty interior and are non persistent surgery
neighborhoods. As the two cases
are similar, let us suppose for instance that the interior of
$\mathcal{E}(\gamma_0)$ is not empty. As everything is continuous,
there exists a maximal time $b\in(0,1]$ such that
$\mathcal{E}(\gamma_{t})$ is a non empty and non persistent surgery neighborhood
for all $t\in[0,b]$.

As the sweepout surfaces $(S_t)_{t\in I}$ vary continuously with
$t$, there exists a constant $\eta>0$ as small as we like,
depending only on the original simplicial sweepout $(S_t)_{t\in
I}$ and the choice of the point $x_0$ and the geodesic arc $c$,
such that for every $t\in[0,\eta]$, each point of $S_t$ lies at
distance at most $\mu/2$ from $\Sigma_0=S_0$, and that for every
$t\in[1-\eta,1]$, each point in $S_t$ is at distance at most
$\epsilon K'/2$ from $\Sigma_1=S_1$. If $b\leq\eta$, we do not
modify the sweepout. Otherwise, if $\eta<b<1$, we apply the
surgery construction for all $t\in[\eta,b]$: we replace the
surgery neighborhoods $\mathcal{E}(\gamma_t)$ by a pair of ruled
discs coned from basepoints in the boundary of
$\mathcal{E}(\gamma_t)$ in a continuous way. On the interval
$[0,\eta]$, we replace the surgery neighborhoods by a pair of
discs for $t$ near $\eta$, that get pasted to the initial surgery
neighborhood $\mathcal{E}(\gamma_0)$ as the time is decreasing to
$0$, not too far from the original surface $S_t$ and in a
continuous way. We can do this in such a way that it is still a
modification of a generalized sweepout. If $b=1$, do the same for
all $t\in[1-\eta,1]$. As the diameter of the ruled discs is less
than $\epsilon$ and $K'/2 \geq 1$, one can suppose that every
point in $\widehat{S_t}$ is at distance at most $\epsilon K'$ from
$\Sigma_1$ for all $t\in[1-\eta,1]$.

This ends the proof of proposition
\ref{prop-chirurgies-balayage}.\qed

\subsubsection{Sweepout surfaces and nested surfaces.}

~

To go on with the proof of proposition B, we
now need a lemma to precisely determine the constant $\mu$, which
corresponds to the size of the collar neighborhood of $S_0$ one
has to take into consideration. Set $K':=2 a'\chi_-( C)$. From lemma
\ref{lm-controle-diam-surf-min}, the
number $2\epsilon K'$ is an upper bound for the diameter of the
simplicial surface $\Sigma_1$, identified with $\partial_+C$. Let
$\delta$ be the diameter of the compression body $C$.

\begin{lm}\label{lm-point-loin-du-bord-plus}

There exists a point $x_0$ in the interior of $C$ and lying at
distance at least $(\frac{\delta}{2}-2\epsilon K')$ from
$\partial_+C$.

\end{lm}

\noindent \underline{Proof of lemma
\ref{lm-point-loin-du-bord-plus}.}
\nopagebreak

Suppose that the lemma is false: for every point $z$ in the
interior of $C$, $\dist(z,\partial_+C)<\frac{\delta}{2}-2\epsilon
K'$. For every point $z$ of $C$, the following inequality remains
true: $\dist(z,\partial_+C)\leq\frac{\delta}{2}-2\epsilon K'$.
Take two points $x$ and $y$ in $C$ such that
$d(x,y)=\diam(C)=\delta$. Then,

\begin{eqnarray*}
d(x,y)=\delta&\leq&
\dist(x,\partial_+C)+\diam(\partial_+C)+\dist(y,\partial_+C)\\
&\leq&(\frac{\delta}{2}-2\epsilon K')+2\epsilon K'+
(\frac{\delta}{2}-2\epsilon K')\\
&\leq &\delta-2\epsilon K' < \delta,
\end{eqnarray*}
which is a contradiction, proving lemma
\ref{lm-point-loin-du-bord-plus}.\qed

~

Let $c$ be a length-minimizing geodesic arc between $x_0$ and
$\partial_+C$. Let $\mu$ be the distance between the geodesic $c$
and $\Sigma_0$. As $c$ is embedded in the interior of $C$
(excepted for one extremity which belongs to
$\partial_+C=\Sigma_1$), the constant $\mu$ is strictly positive.

Now, for completeness of the proof of proposition B, we state and
prove a few lemmas which are implicit in \cite[proof of Lemma 4.5
p. 2251]{Mah}.

We recall from definition \ref{df-separer} that if $x$ is a point
in $C$ and $S$ an immersed surface of $C$, we say that
\textbf{$S$ separates $x$ from $\partial_+C$} if every oriented
path from $x$ to $\partial_+C$ has its algebraic intersection
number equal to $+1$.

If two surfaces $S$ and $T$ immersed in $C$ are such that $S$
separates every point of $T$ from $\partial_+C$, we say that
\textbf{$S$ separates $T$ from $\partial_+C$}. In this case, the
surfaces $S$ and $T$ are said to be \textbf{nested}.

\begin{lm}\label{lm-cns-separer}

A point $x$ lying in the interior of $C$ is separated from
$\partial_+C$ by $\widehat{S}_t$ if and only if there exists a
path $\gamma$ from $x$ to $\partial_+C$ intersecting the surface
$\widehat{S}_t$ with algebraic intersection number $+1$.

\end{lm}

\noindent \underline{Proof of lemma \ref{lm-cns-separer}.}
\nopagebreak

It suffices to show that if there exists a path $\gamma$ from $x$
to $\partial_+C$ with algebraic intersection number with
$\widehat{S}_t$ equal to $+1$, then every path $\gamma'$ from $x$
to $\partial_+C$ intersects $\widehat{S}_t$ with algebraic
intersection number $+1$.

Let $\gamma'$ be another path from $x$ to $\partial_+C$ in $C$. As
the immersed surface $\widehat{S}_t$ is homologous to
$\partial_-C$, the homology class of $[\widehat{S}_t]$ is equal to
zero in $H_2(C, \partial_-C)$. The composition
$\alpha=\gamma^{-1}\cdot \gamma'$ is a 1-cycle in $H_1(C,
\partial_+C)$.

As $\partial C=\partial_-C\cup \partial_+C$, $[\alpha]\cdot
[\widehat{S}_t]=[\alpha]\cdot 0=0$, and thus $\gamma'\cdot
\widehat{S}_t =\gamma\cdot \widehat{S}_t=+1$, proving lemma
\ref{lm-cns-separer}. \qed

~

For all $t\in[0,1]$, let $D_t$ be the closure of the set of points
$x\in C$ separated from $\partial_+C$ by $\widehat{S}_t$. As the
immersed surfaces $\left(\widehat{S}_t\right)_{t\in[0,1]}$ are
generalized sweepout surfaces of the compression body $C$, $D_0$
is the starting sweepout surface $\widehat{S}_0$, and $D_1$ is
equal to the whole compression body $C$. Let $E_t$ be the
component of $D_t$ containing $x_0$. As before, $E_0$ is a complex
of dimension at most 2, and $E_1=C$.

\begin{lm}\label{lm-bord-de-D_t}

The boundary of the set $D_t$ is the surface $\widehat{S}_t$.

\end{lm}

\noindent \underline{Proof of lemma \ref{lm-bord-de-D_t}.}
\nopagebreak

Let $x$ be a point in $\widehat{S}_t$. As $\widehat{S}_t$ is a
generalized sweepout surface for the compression body $C$, there
exists a path $c$ from $x$ to $\partial_+C$ such that for every
point $y$ of $c$ distinct from $x$, the path $c_{|y}$ obtained
from $c$ by deleting the interval $[x,y)$ does not intersect the
surface $\widehat{S}_t$. So in particular, for every point $y$ on
$c$ distinct from $x$, the algebraic intersection number between
$c_{|y}$ and $\widehat{S}_t$ is zero, and $y$ is in the complement
of $D_t$ in $C$. As the point $x$ is a limit of such points $y$,
$x$ is in the closure of the complement of $D_t$ in $C$. But as
$x$ is on the surface $\widehat{S}_t$ and that this surface
separates the compression body $C$, every point close enough to
$x$ and on the other side of $\widehat{S}_t$ with respect to $y$
is separated by $\widehat{S}_t$ from $\partial_+C$. As the set
$D_t$ is closed and that $x$ is a limit of points of $D_t$, the
point $x$ also belongs to $D_t$. Therefore, the point $x$ lies in
the boundary of $D_t$, and the surface $\widehat{S}_t$ is a subset
of the boundary of $D_t$.

To get the reverse inclusion, let us assume that there exists a
point $x$ in the boundary of $D_t$ which does not belong to the
surface $\widehat{S}_t$, and seek for a contradiction. The
distance $d=\dist(x,\widehat{S}_t)$ is then strictly positive. As
the point $x$ belongs to the boundary of $D_t$, there exists a
point $y$ in the complement of $D_t$ in $C$ such that
$d(x,y)\leq \frac{d}{2}$. As $y$ is in the complement of $D_t$,
there is a path $c$ from $y$ to the boundary $\partial_+C$ with
algebraic intersection number with $\widehat{S}_t$ different from
$+1$. Let $c'$ be a minimizing geodesic from $x$ to $y$: as the
length of $c'$, which is equal to the distance between $x$ and
$y$, is strictly less than the distance of $x$ to $\widehat{S}_t$,
the geodesic $c'$ does not intersect the surface $\widehat{S_t}$.
If $c''=c'\cup c$, $c''$ is a path from $x$ to $\partial_+C$
with algebraic intersection number with $\widehat{S}_t$ not equal
to $+1$. Therefore, the point $x$ is not separated from
$\partial_+C$ by $\widehat{S}_t$.

\begin{center}

\includegraphics[width=.7\textwidth]{separation}

\end{center}

But as the point $x$ belongs also to $D_t$, there exists a point
$z$ in $C$ separated from $\partial_+C$ by $\widehat{S}_t$ and
such that the distance between $z$ and $x$ is less than
$\frac{d}{2}$. Take a minimizing geodesic $a$ from $z$ to $x$. Let
us denote by $b=a\cup c''$. The path $b$ is linking $z$ to
$\partial_+C$, which implies that the algebraic intersection
number of $b$ with $\widehat{S}_t$ is equal to $+1$. From the
other hand, the distance between $z$ and $x$ is at most
$\frac{d}{2}<\dist(x,\widehat{S}_t)$, which implies that the
minimizing geodesic $a$ does not intersect the surface
$\widehat{S}_t$. But then, the algebraic intersection number of
the path $b=a\cup c''$ with the surface $\widehat{S}_t$ is not
equal to $+1$, which contradicts the fact that $z$ is separated
from $\partial_+C$ by $\widehat{S}_t$. Thus, the point $x$
necessarily belongs to the surface $\widehat{S}_t$, which ends the
proof of lemma \ref{lm-bord-de-D_t}.\qed

~

\begin{lm}\label{lm-connexite-bord-de-E_t}

For every time $t$, the boundary of $E_t$ is connected.

\end{lm}

\noindent \underline{Proof of lemma \ref{lm-connexite-bord-de-E_t}.}
\nopagebreak

Indeed, if the boundary of $E_t$ would not be connected, it would
have at least two components $S$ and $T$ of $\widehat{S}_t$. But
then, $S$ and $T$ would be two disjoint and separating surfaces in
the compression body $C$. If they are not nested, the set of the
points separated from $\partial_+C$ by $S$ is disjoint to the set
of points separated from $\partial_+C$ by $T$, which contradicts
the fact that $E_t$ is connected. Therefore, the surfaces $S$ and
$T$ are nested. But the surface $\widehat{S}_t$ is homotopic to a
surface obtained from $\partial_+C$ by surgeries and as surgeries
preserve the algebraic intersection number in homology, two
components of the same surface $\widehat{S}_t$ cannot be nested,
which ends the proof of lemma \ref{lm-connexite-bord-de-E_t}.\qed

~

To prove proposition B, we will pick up the desired
nested surfaces among the family of connected surfaces $(\partial
E_t)_{t\in[0,1]}$.

~

\noindent \underline{End of proof of proposition B.}
\nopagebreak

Let $c$ be the length-minimizing geodesic arc from the point $x_0$
obtained in lemma \ref{lm-point-loin-du-bord-plus} to
$\partial_+C$. As before, denote by $\mu$ the distance between the
geodesic $c$ and $\Sigma_0$. Let $L$ be the length of $c$. One has
$L\geq\frac{\delta}{2}-2\epsilon K'$. Take $\ell\mapsto c(\ell)$
an arc-length parameterization of  $c$, such that $c(0)=x_0$ and
$c(L)=y_0 \in\partial_+C$.

First, let us show that $E_t=\emptyset$ for $t\in[0,\eta]$, where
$\eta$ is the constant given by proposition
\ref{prop-chirurgies-balayage}. As every original sweepout surface
$S_t$ is contained in a $\mu/2$-neighborhood of
$S_0=\widehat{S_0}$ for all $t\leq \eta$, and that the distance
between $c$ and $S_0$ is at least $\mu$, the geodesic $c$ does not
meet the sweepout surfaces $S_t$ for every $t\leq \eta$. As the
new sweepout surfaces $\widehat{S_t}$ are obtained from the
surfaces $S_t$ by surgery, the intersection number between $c$ and
$\widehat{S_t}$ is the same as the intersection number between $c$
and $S_t$, so it is zero for $t\leq\eta$. Therefore, the geodesic
$c$ is an arc joining $x_0$ to $\partial_+C$ with intersection
number with $\widehat{S_t}$ equal to zero for $t\leq\eta$. By
definition, the surfaces $\widehat{S_t}$ do not separate $x_0$
from $\partial_+C$ for $t\leq \eta$, showing that there is no
component of $D_t$ containing $x_0$. Thus, $E_t=\emptyset$ for
every $t\in[0,\eta]$.

Let us assume that $\frac{\delta}{2}-6\epsilon K'\geq
5\epsilon K$. As the sets $E_t$ vary continuously with the time
$t$, the function $\mathcal{L}$ which maps the time $t$ to the
length of $c\cap E_t$ is a continuous map. From the fact that
$\mathcal{L}(\eta)=0$ and $\mathcal{L}(1)=L$ the length of $c$, we
deduce that there is a time $t_1\in(\eta,1)$ such that
$\mathcal{L}(t_1)=L-2\epsilon(1+K+K')$. Let $S_1$ be the
boundary of $E_{t_1}$. From lemma \ref{lm-connexite-bord-de-E_t},
the immersed surface $S_1$ is a connected component
of $\widehat{S}_{t_1}$. As $c$ is a minimizing arc-length
parametrized geodesic, for every $a$ and $b\in[0,L]$, we have
$d(c(a),c(b))=\left|b-a\right|$. Thus, the intersection point
$c(\mathcal{L}(t_1))$ between $S_1$ and $c$ is lying at distance
$2\epsilon(1+K+K')$ from $\partial_+C$. Since by
construction every point in the surface $\widehat{S_t}$ for $t\geq
1-\eta$ is at distance at most $\epsilon K'$ from
$\partial_+C$, necessarily $t_1<1-\eta$. As the sets
$E_t$ are empty for $t\leq\eta$, in fact $\eta<t_1< 1-\eta$. By
definition of $E_{t_1}$, the surface $S_1$ separates $x_0$ from
$\partial_+C$. By proposition \ref{prop-chirurgies-balayage}, $S_1$ is
connected. Let us show that its $\epsilon$-diameter is at most $K$. By
proposition \ref{prop-chirurgies-balayage}, the diameter in $C$ of a component
of $\widehat{S_t}$ is at most $\epsilon(1+2K+2K')$. Furthermore, if $S_1$
contains a persistent surgery neighborhood, it means that $S_1$ intersects
$\mathcal{N}_+$. That implies that every point of $S_1$ is at distance at most
$\epsilon(1+2K+2K')+\epsilon/2$ of $\partial_+C$, contradicting the fact that
the intersection point between $S_1$ and $c$ is at distance $2\epsilon (1+K+K')
> \epsilon(1+2K+2K')+\epsilon/2$ of $\partial_+C$. Thus, $S_1$ does not contain
any persistent surgery neighborhood. Proposition \ref{prop-chirurgies-balayage}
ensures that its intrinsic $\epsilon$-diameter is at most $K$ and its diameter
in $C$ is at most $2\epsilon K$. Therefore, the surface $S_1$ cannot meet
$\{c(\ell)\,,
\,0\leq\ell< L-2\epsilon (1+K+K')-2\epsilon K)\}\cup \{c(\ell)\,,
\,L-2\epsilon(K+ K') <\ell\leq L\}$. Let $\ell_1$ be the smallest
value of $\ell$ such that $c(\ell)\in S_1$. We have $L-2\epsilon(1+2K+K')
\leq\ell_1\leq L-2\epsilon (1+K')$. As $K'>1$, this implies that
$L-2\epsilon(1+2 K + K')\geq L-4\epsilon(K+K') \geq
\frac{\delta}{2}-4\epsilon K-6 \epsilon K'\geq\epsilon K>0$.

Let $c_1=\{c(\ell)\,,\,0\leq\ell\leq\ell_1-14\epsilon K\}$.
Replacing $c$ by $c_1$, we can iterate the previous process. If
$K$ is small enough compared to $\delta$, there exists a time
$t_2$ such that the length of $c_1\cap E_{t_2}$ is equal to:
$\lgr(c_1)-2\epsilon K=\ell_1-16\epsilon K\geq L-20\epsilon
K -2\epsilon(1+ K')\geq L-20 \epsilon K -4\epsilon K'$. For the same reasons as
before, the boundary of $E_{t_2}$ is a surface $S_2$ which is a connected
component of $\widehat{S}_{t_2}$ separating $x_0$ from $\partial_+C$, and it
intersects $c_1$ only on the set $\{c_1(\ell)\,,\,(\ell_1-14\epsilon
K)-4\epsilon K\leq \ell\leq \ell_1-14\epsilon K\}$. Furthermore, the surface
$S_2$ is too far from the boundary $\partial_+C$ to contain a persistent
surgery neighborhood, and its intrinsic $\epsilon$-diameter is at most $K$ by
proposition \ref{prop-chirurgies-balayage}.

Let us prove that the distance between the surfaces $S_1$ and
$S_2$ is less than or equal to $10\epsilon K$. Let $\ell_2$ be
the smallest real number $\ell$ such that $c(\ell)\in S_2$.
From the former discussion, $\ell_2\leq \ell_1-14\epsilon K$. As
$c(\ell_1)\in S_1$ and $c(\ell_2)\in S_2$, we have:
\begin{eqnarray*}
\dist(S_1,S_2)&\geq&\dist(c(\ell_1),c(\ell_2))
-\diam(S_1) -\diam(S_2)\\
&\geq&(\ell_1-\ell_2)-4\epsilon K\\
&\geq&14\epsilon K-4\epsilon K=10\epsilon K.
\end{eqnarray*}

We can iterate the process with $c_2=\{c(\ell)\,
,\,0\leq \ell\leq\ell_2-14\epsilon K\}$, on condition that
$\ell_2 -14\epsilon K> 4\epsilon K$, so for example if
$L-2\times 18\epsilon K-4\epsilon K'> 4\epsilon K$.

The iteration process stops when $L-18\epsilon K(n-1) -4\epsilon K'>
4\epsilon K$ but $L-18\epsilon K n -4\epsilon K'\leq
4\epsilon K$, so for
$n=\lceil\frac{L-4\epsilon(K+K')}{18\epsilon K}\rceil$. As
$L\geq\frac{\delta}{2}-\epsilon K'$, $n\geq\lceil
\frac{\delta}{36\epsilon K}-\frac{2}{9}
-\frac{K'}{3K}\rceil$, which proves proposition B.\qed

\subsection{Proof of Proposition C: from nested to parallel
surfaces.}\label{dem-prop-surf//}

~

With proposition B, we know that we can find
$n=\lceil \frac{\delta}{36\epsilon
K}-\frac{2}{9}-\frac{K'}{3K}\rceil$ immersed surfaces in the
compression body $C$ of the cover $M'$. All those surfaces are
nested, their $\epsilon$-diameter is at most $K$ and they are at
distance at least $10\epsilon K$ from each other, where $K=4
\left(3+1/\sinh^2(\epsilon/8)\right)g(C) -10$. Furthermore, all
those surfaces are homotopic to embedded surfaces obtained from
$\partial_+C$ by surgery.

Thus the genus of those immersed surfaces is between $0$ and
$g(C)= g(\partial_+C)$. So there are at least $n'=\lfloor
n/(g(C)+1) \rfloor$ surfaces $S_1,\ldots,S_{n'}$ with the same
genus, and this genus is at most $g(C)$. We take the indices $j$
such that $S_{j+1}$ separates $S_j$ from $\partial_+C$.

We then follow the proof of Maher \cite[p. 2252--2257]{Mah}. Let
$S=S_j$ be one of the previous immersed and nested surfaces with
the same genus. A \textbf{collection $\Delta_S$ of compression
discs of $\partial_+C$ to get $S$} is a finite set of properly
embedded discs in $C$, such that the sweepout gives a homotopy
from $S$ to a subset of $\partial_+C \cup \Delta_S$. The first
step is to show that for two connected and nested sweepout
surfaces, one can choose collections of compression discs such
that one of them is a subset of the other one. This is done in
\cite[Lemma 4.6 p. 2252]{Mah}. In particular, if the two surfaces
have the same genus, they are homotopic.

\begin{lm}\cite[Lemma 4.6]{Mah}\label{lm4.6}

Let $S_1$ and be two of the immersed surfaces obtained in
proposition B. Suppose for example that $S_2$ separates $S_1$ from
$\partial_+C$. Then we can choose a collection of compression
discs of $\partial_+C$, say  $\Delta_{S_1}$ to get $S_1$ and
$\Delta_{S_2}$ to get $S_2$, such that $\Delta_{S_2}$ is a subset
of $\Delta_{S_1}$. In particular, if the two surfaces $S_1$ and
$S_2$ have the same genus, $\Delta_{S_1}=\Delta_{S_2}$.\qed

\end{lm}

This lemma shows that all the nested surfaces $S_1, \ldots,
S_{n'}$ are homotopic, as they have the same genus.

The following lemma is crucial: we wish to replace the nested
immersed surfaces by embedded surfaces of the same genus in an
arbitrarily small neighborhood of the original immersed surfaces.
This lemma is proven in \cite[Lemma 4.7 p. 2253]{Mah}.

\begin{lm}\cite[Lemma 4.7]{Mah}\label{lm4.7}

Let $S$ be one of the surfaces obtained in proposition B. Let $T$
be a least genus, connected and embedded surface, separating $S$
from $\partial_+C$. Then $T$ is incompressible in $C\setminus S$
and the genus of $T$ is greater than or equal to the genus of $S$.

\end{lm}

\noindent \underline{Proof of lemma \ref{lm4.7}.}
\nopagebreak

We recall here Maher's proof.

If the surface $T$ were compressible in $C\setminus S$, it could
be compressed along embedded discs in $C\setminus S$ to obtain a
new surface $T'$ embedded in $C\setminus S$. But one component of
$T'$ would be an embedded surface in $C$ separating $S$ from
$\partial_+C$, with genus strictly less than the genus of $T$,
which is a contradiction. So the surface $T$ is incompressible in
$C\setminus S$.

The surface $S$ is homotopic to $\partial_+C$ compressed along a
collection $\Delta_S$ of embedded discs. Thus, if $C'$ is the
component of $C\setminus \Delta_S$ containing the surface $S$,
$C'$ is a compression body and we can find for it a spine $\Gamma$
that is homotopic to the immersed surface $S$. The map on first
homology $H_1(\Gamma)\rightarrow H_1(C)$ induced by the inclusion
of $\Gamma$ in $C$ is injective.

The surface $T$ is an embedded surface in the compression body
$C$, so it is separating and there exists a set $D_T$ of embedded
compression discs for $T$ such that $T$ compressed along $D_T$ is
parallel to some components of $\partial_-C$ (c.f. \cite[Lemma
2.3]{Bon1}). As $T$ is incompressible in $C\setminus S$, the
compression discs of $D_T$ for the surface $T$ are only in one
side of $T$. So the surface $T$ bounds a compression body $C''$ in
$C$. As the composition of the maps induced by the inclusions
$H_1(\Gamma)\rightarrow H_1(C'') \rightarrow H_1(C)$ is injective,
the map  $H_1(\Gamma)\rightarrow H_1(C'')$ is injective. Thus the
rank of $H_1(C'')$ is greater than or equal to the rank of
$H_1(\Gamma)$, and necessarily the genus of $T$ is greater than or
equal to the genus of $S$.\qed

~

A consequence of lemma \ref{lm4.6} is that all the nested and
immersed surfaces\\ $S_1, \ldots,S_{n'}$ are homotopic. We want a
little more: we need to find for all $j$ between $1$ and $(n'-1)$
a homotopy between $S_j$ and $S_{n'}$ that is disjoint from $S_k$
for all $k<j$. We follow the arguments of the proof of \cite[Lemma
4.8 p. 2254]{Mah}, but we compute precise upper bounds.

\begin{lm}\label{lm4.8}

From the surfaces $S_1, \ldots, S_{n'}$, one can construct a
collection of immersed surfaces $S_1', \ldots, S_{n'-1}', S_{n'}'$
which are disjoint, nested and homotopic, and the homotopy from
$S_{n'}'$ to $S_j'$ is disjoint from $S_k'$ for $1\leq k<j$.
Furthermore, the diameter of the surfaces $S_j'$ is at most
$8\epsilon K$, they are at distance at least $2\epsilon K$ from
each other, and the $\epsilon$-diameter of $S_2',\ldots,
S_{n'-1}'$ is at most $K$.

\end{lm}

\noindent \underline{Proof of lemma \ref{lm4.8}}
\nopagebreak

Each surface $S_j$ admits a one-vertex triangulation with
edge-length bounded by $4\epsilon K$, and its diameter is at most
$2\epsilon K$. Therefore, by lemma \ref{lm-surf-min->surf-simpl}
the surfaces $S_1$ and $S_{n'}$ are homotopic to simplicial
surfaces $S_1'$ and $S_{n'}'$ with diameter at most $4\epsilon K$
and such that for every points $x\in S_j$ and $x'\in S_j'$ (where
$j=1$ and $n'$), the distance between $x$ and $x'$ is at most
$6\epsilon K$. In fact, by construction of $S_j'$, each point of
$S_j'$ is at distance at most $4\epsilon K_i$ from the original
surface $S_j$.

The homotopy between the two simplicial surfaces $S_1'$ and
$S_{n'}'$ can be modified into a simplicial sweepout as in section
\ref{dem-lm4.5}. By proposition \ref{prop-chirurgies-balayage},
there exists a finite sequence of surgeries of generalized
sweepouts, starting from this simplicial sweepout and ending to a
generalized sweepout in which all the sweepout surfaces $S_t'$ for
$t\in[\eta,1-\eta]$ have $\epsilon$-diameter bounded above by $K$.
We can use the same constant $K$ as before since the genus of the
surfaces $S_j$ is at most $g(C)$. Moreover, the surfaces $S_t'$
are homotopic to the surface $S_{n'}$ after some compressions if
necessary. For $j$ between $2$ and $(n'-1)$, let $S_j'$ be the
first sweepout surface $S_t'$ intersecting $S_j$. As $S_j'$ is a
generalized sweepout surface, its $\epsilon$-diameter is at most
$K$.

We know from the construction of a generalized sweepout that the
genus of the surface $S_j'$ is at most the genus of the surface
$S_j$. In fact, we show that those two genera are equal.

\begin{claim}

For all $1\leq j \leq n'-1$, the genus of the surface $S_j'$ is
the same as the genus of the original sweepout surface $S_j$.

\end{claim}

Assuming the claim, since the modified sweepout surfaces $S_j'$
have the same genus as the original sweepout surfaces $S_j$, in
fact there is no compression to obtain the surfaces $S_j'$ and
they were already sweepout surfaces of the original simplicial
sweepout between $S_1'$ and $S_{n'}'$. So the surfaces $S_j'$
are homotopic to the surface $S_{n'}'$, and by definition of a
sweepout, this homotopy is disjoint from the surfaces $S_k'$ for
every $k<j$.

~

\noindent \underline{Proof of claim.}
\nopagebreak

Suppose that there exists some $j$ such that the genus of $S_j'$
is strictly less than the genus of $S_j$. By a result of Gabai, we
can then replace our simplicial surface $S_j'$ by an embedded
surface $T_j'$ in an arbitrarily small neighborhood of the
immersed surface $S'_j$. More precisely, take a small regular
neighborhood $N(S_j')$ of the immersed surface $S_j'$. This
neighborhood contains embedded surfaces in the same homology class
as $S_j'$ in $H_2(N(S_j'),
\partial N(S_j'))$. Gabai showed that the singular norm on
homology is the same as the embedded Thurston norm \cite{Ga},
hence there exists an embedded surface $T_j'$ in $N(S_j')$ with
the same homology class as $S_j'$ and of genus less than or equal
to the genus of $S_j'$. If we choose sufficiently small
neighborhoods $N(S_j')$, we can ensure that the diameter of the
embedded surface $T_j'$ is less than $3\epsilon K$. In particular,
as the surfaces $S_1'$ and $S_{n'}'$ are too far away, the
embedded surface $T_j'$ is disjoint from $S_1'$ and $S_{n'}'$, and
it is separating $S_1'$ from $S_{n'}'$. Applying lemma
\ref{lm4.7}, we see that the genus of $T_j'$ must be at least the
genus of $S_1'$: $g(T_j)\geq g(S_1')$. But as the genus of $S_1'$
is the same as the genus of $S_j$, and that the genus of $T_j'$ is
at most the genus of $S_j'$, which we have supposed strictly less
than the genus of $S_j$, we have $g(T_j')<g(S_1')$, which is a
contradiction.\qed

~

As the surfaces $S_j$ were at distance at least $10\epsilon K$
from each other and that every point of $S_j'$ is at distance at
most $4\epsilon K$ from the original surface $S_j$ for all
$j=1,\ldots, n'$, the new surfaces $S_j'$ are at distance at most
$2\epsilon K$ from each other (which also shows that the surfaces
$S_j'$ are all disjoint). Furthermore, their diameter is bounded
from above by $8\epsilon K$ and the $\epsilon$-diameter of $S_2',
\ldots, S_{n'-1}'$ is at most $K$.

There remains to show that the surfaces $S_1',\ldots,S_{n'}'$ are
nested. In the spirit of the proof of proposition B, let us denote
by $D_{n'}$ the closure of the subset of the points of $C$
separated from $\partial_+C$ by $S_{n'}'$. For all $j<n'$, the
surface $S_j'$ intersects the surface $S_j$, which lies in
$D_{n'}$. As $S_j'$ is at distance at least $2\epsilon K$ from
$S_{n'}'=\partial D_{n'}$, $S_j'$ is contained in the interior of
$D_{n'}$. So it is separated from $\partial_+C$ by $S_{n'}'$.
Therefore, if we denote by $D_j$ the closure of the points of $C$
separated from $\partial_+C$ by $S_j'$, $D_j\subset D_{n'}$. Let
$1\leq k <j<n'$. If we take a point $x$ in $D_k$, as $D_k\subset
D_{n'}$, every path $\gamma$ from $x$ to $\partial_+C$ has its
algebraic intersection number with $\partial_+C$ equal to $+1$. As
the surface $S_j'$ is homotopic to $S_{n'}'$ by a homotopy that is
disjoint from $S_k'$, this homotopy does not change the
intersection number, so the intersection number of $\gamma$ with
$S_j'$ is still equal to $+1$, and $x$ is in $D_j$. Thus
$D_k\subset D_j$ for $1\leq k<j\leq n'$, showing that the surfaces
$S_1', \ldots, S_{n'}'$ are nested. This ends the proof of lemma
\ref{lm4.8}.\qed

~

In the sequel, we replace the family $S_1,\ldots, S_{n'}$ by the
new family $S_1', \ldots, S_{n'-1}', S_{n'}'$ of surfaces obtained
by lemma \ref{lm4.8}, and for simplicity, we will still denote
this family by $S_1,\ldots, S_{n'}$.

We then wish to replace our immersed surfaces by embedded surfaces
in an arbitrarily small neighborhood of the immersed surfaces. It
is the aim of the following lemma.

\begin{lm}\label{lm-surf-plongees}

For every $j$ from $1$ to $n'$, there exists an embedded surface
$T_j$ in a small regular neighborhood of $S_j$, with the same
genus as $S_j$, and which can be covered by at most
$\diam_\epsilon(S_j)\leq K$ embedded balls in $M'$ of radius
$2\epsilon$. Furthermore, two surfaces $T_j$ and $T_k$ for $j\neq
k$ are at distance at least $\epsilon K$ from each other.

\end{lm}

\noindent \underline{Proof of lemma \ref{lm-surf-plongees}.}
\nopagebreak

Take a small regular neighborhood $N(S_j)$ of one of the immersed
and nested surfaces $S_j$. As in the proof of the claim, by Gabai
\cite{Ga}, this neighborhood contains an embedded surface $T_j$ in
the same homology class as $S_j$ in $H_2(N(S_j), \partial N(S_j))$
and of genus less than or equal to the genus of $S_j$. If we
choose sufficiently small neighborhoods $N(S_j)$, we can ensure
that the diameter of the embedded surfaces $T_j$ in the ambient
manifold $M'$ is less than $9\epsilon K$, and two embedded
surfaces $T_j$ and $T_k$ are at distance at least $\epsilon K$.
Furthermore, if we take a set $\mathcal{B}$ of
$\diam_\epsilon(S_j)$ embedded balls of radius $\epsilon$ and
centers on the surface $S_j$, one can choose $N(S_j)$ small enough
such that it is contained in the union of corresponding balls with
the same center and radius $2\epsilon$. Thus, the surface $T_j$
can be covered by at most $\diam_\epsilon(S_j)$ embedded balls of
$M'$ with radius $2\epsilon$.

The genus of $T_j$ is at most the genus of $S_j$, but we wish to
show that in fact, the genus of $T_j$ is the same as the genus of
$S_j$.

With lemma \ref{lm4.7}, we know that the genus of the embedded
surface $T_j$ for $j=2,\ldots,n'$ is greater than or equal to
the genus of the immersed surface $S_1$ that it separates from
$\partial_+C$. But as the genus of $T_j$ is at most the genus
of $S_j$, which is equal to the genus of $S_1$, in fact the genus
of $T_j$ is equal to the genus of $S_j$: the surfaces $T_2,\ldots,
T_{n'}$ have the same genus as the immersed surfaces
$S_2,\ldots, S_{n'}$. This proves lemma \ref{lm-surf-plongees}.\qed

~

The final step in the proof of proposition C is to show that some
of the embedded surfaces are actually parallel.

\begin{lm}\label{lm-surf//}

The embedded surfaces $T_4,\ldots, T_{n'-1}$ are parallel.

\end{lm}

\noindent \underline{Proof of lemma \ref{lm-surf//}.}
\nopagebreak

This lemma relies on homological arguments, see \cite[Lemmas 4.9
to 4.11]{Mah}. For completeness, we give here a shorter proof,
based on classical 3-manifold topological results.

Let $V$ be the 3-complex in $C$ bounded by the immersed surfaces
$S_1$ and $S_{n'}$. There is a sweepout $\phi$ between $S_1$ and
$S_{n'}$ such that for each $1\leq j \leq n'$, the surface $S_j$
is a sweepout surface. In other words, the application induced by
the map $\phi \, :\,S\times I\rightarrow V$ in homology
$\phi_*\,:\, H_3(S\times I,\partial (S\times I)) \rightarrow
H_3(V, \partial V)$ is an isomorphism and for each $j$, there
exists a time $t_j\in I$ such that $S_j=\phi(S\times \{t_j\})$.
Moreover, we have $0=t_1< t_2< \ldots <t_{n'}=1$.

By a classical construction (see \cite[point 3. p. 96]{Sta} for
example), we can homotop the sweepout $\phi$ to a map $\phi'$
which is still degree one, and such that for every $2\leq j\leq
n'$, $\phi'^{-1}(T_j)$ is an embedded incompressible surface (not
necessarily connected) in $S\times I$.

Take $3<j<k\leq n'-1$. As the homology class of the surfaces $T_j$
and $T_k$ is the same as the homology class of $S_3$, the homology
class of the preimages $\phi'^{-1}(T_j)$ and $\phi'^{-1}(T_k)$ in
$H_2(S\times [t_3,1], \partial (S\times [t_3,1]))$ is the same as
the homology class of the fiber $S\times \{t\}$. As those
preimages are incompressible embedded surfaces, $\phi'^{-1}(T_j)$
and $\phi'^{-1}(T_k)$ are each composed of an odd number of
connected surfaces isotopic to the fiber $S\times \{t\}$ with
total algebraic intersection number with any path from
$S\times\{t_3\}$ to $S\times\{1\}$ equal to $+1$. Up to isotopy,
we can suppose that there exist times $t_3<t_1^j< \ldots
<t_{2n_j+1}^j$ and $t_3<t_1^k<\ldots < t_{2n_k+1}^k$ such that
$\phi'^{-1}(T_j)=\cup_{\ell=1}^{2n_j+1} \epsilon_\ell^j (S\times
\{t_{\ell}^j\})$ and $\phi'^{-1}(T_k)=\cup_{\ell=1}^{2n_k+1}
\epsilon_\ell^k (S\times \{t_{\ell}^k\})$, with $\epsilon_\ell^j$
and $\epsilon_\ell^k$ equal to $+1$ or $-1$, depending on the
orientation of the component of $\phi'^{-1}(T_j)$ or
$\phi'^{-1}(T_k)$ corresponding to the fiber $S\times
\{t_\ell^j\}$ or $S\times\{t_\ell^k\}$. As
$\sum_{\ell=1}^{n_j+1}\epsilon_\ell^j=+1$ and
$\sum_{\ell=1}^{n_k+1}\epsilon_\ell^k=+1$, there exists $\ell$ and
$\ell'$ such that $\epsilon_\ell^j=+1=\epsilon_{\ell'}^k$. Suppose
for example that $t_\ell^j<t_{\ell'}^k$. Then
$\phi'\,:\,S\times[t_\ell^j,t_{\ell'}^k]\rightarrow V$ is a
homotopy between the embedded surfaces $T_j$ and $T_k$ contained
in the region in $V$ bounded by $S_3$ and $S_{n'}$. As the
embedded surface $T_2$ is not in this region, if we denote by $Y$
the submanifold of $C$ bounded by $T_2$ and $\partial_+C$, the two
embedded surfaces $T_j$ and $T_k$ are homotopic in the interior of
$Y$.

By lemma \ref{lm4.7}, the surfaces $T_j$ and $T_k$ are
incompressible in $C\setminus S_1$. As they are contained in the
interior of $Y$ and $Y$ is included in the component of
$C\setminus S_1$ containing $T_j$ and $T_k$, the surfaces $T_j$
and $T_k$ are incompressible in $Y$. Thus, by a result of
Waldhausen \cite[Corollary 5.5 p. 76]{Wal}, they are in fact
isotopic in $Y$. Therefore, $T_j$ and $T_k$ are parallel in $C$,
for $3<j<k\leq n'-1$. Thus we have $m=n'-4$ embedded surfaces
$T_4,\ldots , T_{n'-1}$ parallel in the compression body $C$,
which ends the proof lemma \ref{lm-surf//}. As the
$\epsilon$-diameter of $S_2,\ldots, S_{n'-1}$ is at most $K$ and
the surfaces $T_4,\ldots, T_{n'-1}$ can be covered by at most $K$
embedded balls in $M'$ of radius $2\epsilon$, this ends also the
proof of proposition C.\qed

~\qed

\subsection{Proof of Proposition D: from patterns of fundamental domains to
virtual fibration.}\label{dem-prop-motifs-dom-fond}

~

This part is dedicated to the proof of Proposition D,
which is based on \cite[Lemma 4.12 p.2258]{Mah}. This proof is
much involved than the one of Lemma 4.12 in \cite{Mah}, which is
too quick for our purpose since we need explicit bounds and
precise constants.

~

Assume that there are $m$ connected, orientable, embedded and
disjoint parallel surfaces in $M'$. Furthermore, suppose that each
of those surfaces can be covered by at most $K$ embedded
balls in $M'$ of radius $2\epsilon$ and that any two surfaces are
at distance at least $r>0$ from each other. In particular, there
exists an embedded product $T\times[0,m-1]$ in the manifold $M'$
in which the surface $T_j$ coincides with the fiber $T\times\{j\}$
for all $j$ from $0$ to $m-1$.

Let $\mathcal{D}$ be a Dirichlet fundamental domain for the
manifold $M$ in its universal cover $\widehat{M}\simeq
\mathbb{H}^3$. The translates of $\mathcal{D}$ by the covering
maps form a tiling of the universal cover $\widehat{M}$. This
tiling descends to a tiling of the cover $M'$ by $d$ copies of
$\mathcal{D}$. Each of the $m$ embedded and parallel surfaces
$T_1,\ldots, T_{m}$ in $M'$ intersects some copies of
$\mathcal{D}$.

\begin{df}

The union in $M'$ of copies intersected by one of the
surfaces $S_j$ is called a \textbf{pattern (of fundamental
domains) for $S_j$} and denoted by $P_j$.

\end{df}

As the surface is connected, a pattern is a connected 3-complex.
We can suppose that each of the embedded surfaces intersects the
2-skeleton of the tiling transversally. More precisely, we can
suppose that each surface does not meet the vertices of the
fundamental polyhedra, that it intersects the edges in isolated
points and it is transverse to the 2-dimensional faces of the
polyhedra. Therefore, a pattern is a connected union of some
copies of $\mathcal{D}$ glued along their 2-dimensional faces. Let
$D$ be an upper bound for the diameter of $\mathcal{D}$, and
$\alpha$ an upper bound for the number of its 2-dimensional faces.

For all $\ell\in\nn$, we recall that $B(\ell)$ is an upper bound
for the number of possibilities of patterns obtained by gluing
together at most $\ell$ fundamental domains. Let $L= \lfloor
\frac{\pi(\sinh(2D+4\epsilon)-2D-4\epsilon)}
{\vol(M)}K\rfloor$ as in lemma \ref{lm-estimation_Dalpha&L}.
The integer $L$ is an upper bound for the number of fundamental
domains a given surface can intersect. Thus, a pattern is the
union of at most $L$ fundamental domains.

Suppose that $r/ (2D+1)\geq 1$ and $\frac{m}{\alpha^2L^2B(L)}\geq
4$ (which will be called condition (a)), or that $r/(2D+1)<1$ and
$\left(\frac{r}{2D+1}m-1\right) \frac{1}{\alpha^2L^2B(L)}\geq 4$
(called condition (b)).

\begin{lm}\label{lm-motifs-disjoints}

If conditions (a) or (b) are satisfied, there are at least $ 4
\alpha^2 L^2 B(L)$ surfaces for which the corresponding patterns
of fundamental domains are disjoints.

\end{lm}

\noindent \underline{Proof of lemma \ref{lm-motifs-disjoints}.}
\nopagebreak

If two surfaces $T_j$ and $T_k$ are at distance strictly more than
$2D$, the patterns of fundamental domains associated to $T_j$ and
$T_k$ are necessarily disjoint, as the diameter of a fundamental
domain is at most $D$.

If $r/(2D+1) \geq 1$ as in condition (a), any pair of surfaces
$T_j$ and $T_k$ with $j\neq k$ are at distance strictly more than
$2D$, and  all the $m$ patterns associated to the parallel
surfaces are disjoint.

Otherwise, $r/2D <1$. In this case, there are at least
$\lfloor\frac{r}{2D+1}m\rfloor \geq \frac{r}{2D+1}m-1$ surfaces
$T_j$ which are separated from each other by a distance at least
$2D+1>2D$. Thus, every corresponding patterns of fundamental
domains are disjoint.

As in condition (a), $m\geq 4 \alpha^2 L^2 B(L)$, or in condition
(b),  $\frac{r}{2D+1}m-1\geq 4 \alpha^2 L^2 B(L)$, there are at
least $4 \alpha^2 L^2 B(L)$ surfaces whose corresponding patterns
are disjoint.\qed

\begin{lm}\label{lm-motifs-identiques}

There exist an "abstract" pattern of fundamental domains $P$ and
at least $4 \alpha^2 L^2$ patterns of fundamental domains $P_j$,
which are disjoint and homeomorphic to $P$. More precisely, for at
least $4 \alpha^2 L^2$ of the previous indices $j$ for which the
corresponding patterns of fundamental domains are disjoint, there
exists a homeomorphism $\varphi_j\,:\,P_j\rightarrow P$ which
preserves polyhedral decomposition and gluing isometries between
the faces of the fundamental domains belonging to the patterns.

\end{lm}

\noindent \underline{Proof of lemma \ref{lm-motifs-identiques}.}
\nopagebreak

The proof is straightforward. Indeed, as a pattern is the union of
at most $L$ fundamental domains, there are at most $B(L)$ possible
patterns. Among the $4 \alpha^2 L^2 B(L)$ disjoint previous
patterns, there are at least $4 \alpha^2 L^2$ of them
corresponding to the same "abstract" pattern $P$.\qed

~

From now on, we only consider $4 \alpha^2 L^2$ indices $j$
satisfying the conclusions of last lemma.

\begin{lm}\label{lm-motif-au-moins-deux-bords}

The number of boundary components of the pattern $P$ is between
$2$ and $\alpha L$.

\end{lm}

\noindent \underline{Proof of lemma \ref{lm-motif-au-moins-deux-bords}.}
\nopagebreak

Each fundamental polyhedron in the pattern $P$ has $\alpha$
2-faces. As $P$ is the union of at most $L$ polyhedra, it has at
most $\alpha L$ 2-faces. It is an upper bound for the number of
boundary components of $P$.

To see that there is at least two boundary components in $P$, it
suffices to show that for example $P_1$ has at least two boundary
components. But as the surface $T_1$ is contained in the interior
of the pattern $P_1$, $P_1\cap (T\times[0,1])\neq\emptyset$ and
$P_1\cap (T\times[1,2])\neq\emptyset$. The pattern $P_1$ is
disjoint to $T_0$ and $T_2$, so the product regions $T\times[0,1]$
and $T\times[1,2]$ are not contained in $P_1$. By connexity of
$T\times[0,1]$ and $T\times[1,2]$, the boundary of the pattern
$P_1$ has at least two components, one as a subset of
$T\times(0,1)$ and the other one in $T\times(1,2)$. This proves
lemma \ref{lm-motif-au-moins-deux-bords}.\qed

~

Set $\partial P=E_1\cup E_2\cup\ldots\cup E_s$, where the immersed
surfaces $E_j$ are the boundary components of the pattern $P$,
with $2\leq s\leq\alpha L$.

\begin{df}

For every index $j$ between 1 and $4\alpha^2L^2-2$, the pattern
$P_j$ intersects $T\times (j-1,j)$ and $T\times(j,j+1)$. At least
one component of the boundary of $P_j$ is in the boundary of the
component of $(T\times[j-1,j])\setminus (T\times[j-1,j])\cap P_j$
containing the fiber $T\times\{j-1\}$, which we will call a
\textbf{"left" component of the boundary of the pattern $P_j$}.
Similarly, at least one component of the boundary of $P_j$ is in
the boundary of the connected component of
$(T\times[j,j+1])\setminus (T\times[j,j+1])\cap P_j$ containing
the fiber $T\times\{j+1\}$. We will call this component a
\textbf{"right" component for the boundary of $P_j$}.

\end{df}

\begin{lm}\label{lm-motifs-m-comp}

For every index $j$ between $1$ and $4\alpha^2L^2-2$, choose a
left and a right component for the pattern $P_j$  (arbitrarily if
there exist at least two such components). Those two component
correspond to components $E_j^-$ and $E_j^+$ in the boundary of
the abstract pattern $P$. There are at least two indices $j$ and $k$ for
which the pairs of left and right components corresponding to
the patterns $P_{j}$ and $P_k$ coincide in $\partial P$.

\end{lm}

\noindent \underline{Proof of lemma \ref{lm-motifs-m-comp}.}
\nopagebreak

As there are at most $s(s-1)\leq \alpha L(\alpha L-1)<\alpha^2L^2$
pairs of left and right boundary components of $P$, there are at
least $(4 \alpha^2 L^2 -2)/(\alpha^2 L^2)\geq 2$ surfaces $T_j$
and $T_k$ with $1\leq j<k \leq 4\alpha^2L^2-2$, for which the
pairs of left and right components corresponding to the patterns
$P_j$ and $P_k$ coincide.\qed

~

In the sequel, in order to simplify notations, let us denote by
$T_1$ the surface $T_j$, $T_2$ the surface $T_k$ and $T_3$ the
last surface $T_{4\alpha^2 L^2 -1}$. The surfaces $T_0$ and $T_3$
bound a product $T\times[0,3]$ in $M'$, such that
$T_1=T\times\{1\}$ and $T_2=T\times\{2\}$. The two patterns $P_1$
and $P_2$ are contained in the interior of the product
$T\times[0,3]$. Denote by $\psi:=\varphi_{2}^{-1}\circ \varphi_1$
the composed homeomorphism between patterns $P_{1}$ and $P_{2}$.
Let $T_{1}'$ be the image of the surface $T_{1}$ in the interior
of the pattern $P_{2}$ under the action of $\psi$~: $T_1'=
\varphi_{2}^{-1}\circ \varphi_1(T_1)=\psi(T_1)$. It is an embedded
surface in the product $T\times[0,3]$. Clearly, the surfaces $T_1$
and $T_2$ are parallel, but they may not be embedded in their
patterns in the same way. However, the surfaces $T_1$ and $T_1'$
are embedded in the patterns $P_1$ and $P_2$ in exactly the same
way, but there is no evidence to say a priori that those two
surfaces are parallel. It is in fact true, thanks to the following
lemma.

\begin{lm}\label{lm-T_1&T_'//}

The surfaces $T_1$ and $T_1'$ are parallel in $M'$.

\end{lm}

\noindent \underline{Proof of lemma \ref{lm-T_1&T_'//}.}
\nopagebreak

\begin{center}

\includegraphics[width=1\textwidth]{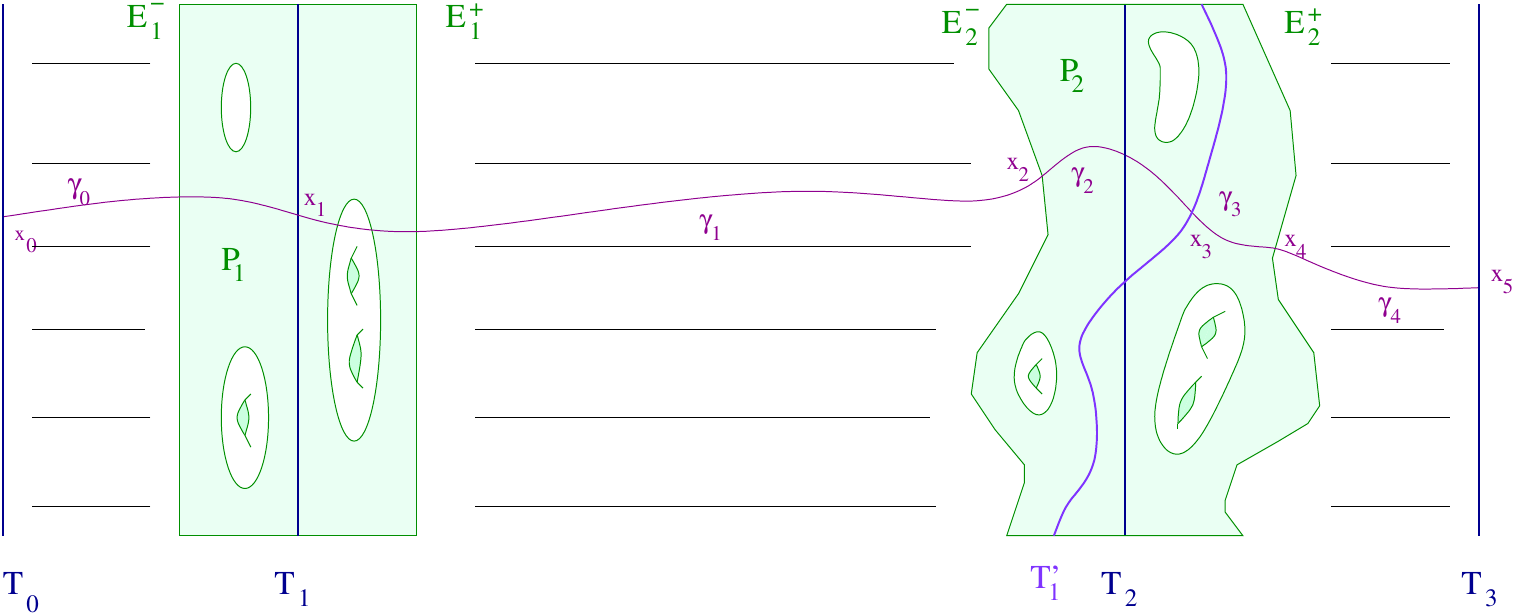}

\end{center}

\begin{claim}

The homology class of $T_{1}'$ in the product $T\times[0,3]$
is equal to the homology class of the fiber $[T]=[T_{1}]=[T_{2}]$.

\end{claim}

\noindent \underline{Proof of claim.}
\nopagebreak

By choice of the surfaces $T_1$ and $T_2$, the left component
$E_1^-$ of the boundary of the pattern $P_1$ and the left
component $E_2^-$ of the boundary of the pattern $P_2$ have the
same image in the pattern $P$:
$\varphi_1(E_1^-)=\varphi_2(E_2^-)$, so $E_2^-=
\varphi_2^{-1}\circ \varphi_1(E_1^-)$. By definition, $E_2^-$ is a boundary
component of the connected component of $(T\times[1,2])\setminus
(T\times[1,2])\cap P_2$ containing the fiber $T_1$, and the
component $E_1^-$ is a boundary component
of the pattern $P_1$ in the boundary of the component of
$(T\times[0,1])\setminus (T\times[0,1])\cap P_1$
containing the fiber $T_0$. As $P_1\cap(T \times[0,1])$ is
connected, there exists a path $\gamma'_2$, properly embedded in
$P_1\cap(T \times[0,1])$ and joining the component $E_1^-$ to the
surface $T_1$. The image by the homeomorphism $\varphi_2^{-1}\circ
\varphi_1$ between the patterns $P_{1}$ and $P_{2}$ of the path
$\gamma'_2$ is a path $\gamma_2= \varphi_{2}^{-1}\circ \varphi_1
(\gamma_2')$ in $P_{2}$ from the boundary component $E_2^-$ to the
surface $T_1'$. The interior of the path $\gamma_2$ is contained
in the interior of the component of $P_2\setminus T_1'$ containing
$E_2^-$. Let $x_2$ be the extremity of $\gamma_2$ belonging to the
boundary component $E_2^-$, and $x_3$ the other one, on the
surface $T_1'$.

Similarly, there exists a path $\gamma_3$ from $x_3$ to a point
$x_4$ lying on the right component $E_2^+$ of the boundary of
$P_2$, and such that its interior is contained in the interior
of the component of $P_2\setminus T_1'$ containing
$E_2^+$.

As $E_2^-$ is in the boundary of the connected component of
$(T\times[1,2])\setminus P_2\cap (T\times[1,2])$
containing the fiber $T_1$, there exists a path $\gamma_1$ with
its interior contained in the interior of this component, and
joining a point $x_1$ of the fiber $T_1$ to the point $x_2$ of
$E_2^-$. Similarly, by choice of $E_2^+$, there exists a
path $\gamma_4$ with interior contained in the interior of the
component of $(T\times[2,3])\setminus P_2\cap
(T\times[2,3])$ containing the fiber $T_3$ and linking the
point $x_4$ of $E_2^+$ to a point $x_5$ of $T_3$.
Eventually, as the product $T\times[0,1]$ is connected, there
exists a path $\gamma_0$ with interior contained in
$T\times(0,1)$ joining the point $x_1$ of $T_1$ to a point
$x_0$ of $T_0$.

Let $\gamma$ be the path obtained by concatenating the paths
$\gamma_0$, $\gamma_1$, $\gamma_2$, $\gamma_3$ and $\gamma_4$. The
path $\gamma$ joins the point $x_0$ of $T_0$ to the point
$x_5$ of $T_3$ and intersects the surface $T_1'$ only once,
at the point $x_3$. As the orientations of the patterns $P_1$
and $P_2$ coincide, the intersection number of $\gamma$ with
the surface $T_1'$ is $+1$. So it is the same as the
intersection number of $\gamma$ with the fiber $T$. By Poincar\'e
duality, as the homology group $H_2(T\times[0,3],\zz)$ is of
rank one, generated by the class of the fiber $[T]$, the class of
the surface $T_1'$ is equal to the class of the fiber in the
homology of the product, showing the claim.\qed

~

As the surface $T_{1}'$ is embedded in the product $T\times[0,3]$,
by a result of Waldhausen \cite{Wal}, it follows that the surface
$T_1'$ is parallel to the fiber $T_{1}$, possibly after performing
a finite number of compressions on $T_1'$. But as the surface
$T'_{1}$ is homeomorphic to $T_1$, it is of the same genus as the
fiber $T_{1}$. So in fact there is no compression. Therefore,
those two surfaces bound a product in $M'$.\qed

\begin{lm}\label{lm-fibre-virtuelle}

The manifold $M$ admits a cover $N$ of finite degree at most $d$
which fibers over the circle, and the embedded surface $T_1$ in
$M'$ is an (incompressible) virtual fiber.

\end{lm}

\noindent \underline{Proof of lemma \ref{lm-fibre-virtuelle}.}
\nopagebreak

One can cut the manifold $M'$ open along those two disjoint
surfaces $T_1$ and $T'_1$. We keep only the component
corresponding to the product region between the two parallel
surfaces, and we identify the two surfaces via the homeomorphism
$\psi={(\varphi_2^{-1}\circ \varphi_1)}_{|T_1}$ to obtain a
manifold $N$ fibering over the circle, with fiber
$\widehat{T_1} =(T_1 \sim T_1')$. The homeomorphism
$\varphi_{2}^{-1}\circ \varphi_1$ identifies the "left" part of
the pattern $P_2$ with the "left" part of the pattern $P_1$,
so we get a pattern $\widehat{P_1}$ corresponding to
$\widehat{T_1}$ in $N$ homeomorphic to the pattern $P$: the
"left" part of this pattern corresponds to the left part of the
pattern $P_2$ via the homeomorphism $\varphi_2$, and the
"right" part of the pattern corresponds to the right part of $P_1$
via the homeomorphism $\varphi_1$. As those homeomorphisms
preserve the gluings between the 2-dimensional faces of the
fundamental domains, the gluings between the fundamental domains
in the pattern $\widehat{P_1}$ are the same as the gluings in
the model pattern $P$. Therefore, we obtain a tiling of $N$ by
finitely many copies of fundamental domains homeomorphic to
$\mathcal{D}$ and with matching gluings. Thus, $N$ is a finite
cover of the original manifold $M$, and $N$ is fibered over the
circle, with fiber $\widehat{T_1}$.

\begin{center}
\begin{tikzpicture}[scale=.3]

\draw [very thick] (0,3) .. controls +(2,0) and +(0,1) .. (3,0) ..
controls +(0,-1) and +(2,0) .. (0,-3) .. controls +(-2,0) and
+(0,-1) .. (-3,0) .. controls +(0,1) and +(-2,0) .. (0,3);

\draw (0,2) -- (1.5,2) -- (2,3.5) -- (0,4) -- (0,2) -- (-1.5,2) --
(-2,3.5) -- (0,4);

\draw [dashed] (1.5,2) -- (2.2,.7);

\draw [dashed] (-1.5,2) -- (-2.2,.7);

\draw [dashed] (2,3.5) -- (3,2.5);

\draw [dashed] (-2,3.5) -- (-3,2.5);

\draw [very thick] (0,6) .. controls +(4,0) and +(0,2) .. (6,0) ..
controls +(0,-2) and +(4,0) .. (0,-6) .. controls +(-4,0) and
+(0,-2) .. (-6,0) .. controls +(0,2) and +(-4,0) .. (0,6);

\draw (0,5) -- (1.5,5) -- (2,6.5) -- (0,7) -- (0,5) -- (-1.5,5) --
(-2,6.5) -- (0,7);

\draw [dashed] (1.5,5) -- (2.7,4.2);

\draw [dashed] (-1.5,5) -- (-2.7,4.2);

\draw [dashed] (2,6.5) -- (3.5,5.7);

\draw [dashed] (-2,6.5) -- (-3.5,5.7);

\node at (1,3.2) {$\star$};

\node at (1,6.2) {$\star$};

\node at (-1,3.2) {$\times$};

\node at (-1,6.2) {$\times$};

\node at (-3.8,0) {$T_{1}'$};

\node at (-6.8,0) {$T_{1}$};

\node at (-6,6) {$M'$};

\draw [->] (0,-7) -- (0,-11);

\draw (-1,-14.5) -- (1,-14.5) -- (1.4,-12.5) -- (-1.2,-12) --
(-1,-14.5);

\node at (-2.5,-13.5) {$M$};

\draw [very thick] (-1.05,-13.5) -- (1.2,-13.3);

\draw [very thick] (1.2,-13.3) -- (-1.07,-13);

\draw [very thick] (-1.07,-13) -- (1.05,-13.9);

\draw [very thick] (1.05,-13.9) -- (-1.05,-13.5);

\draw [very thick] (11,-3) .. controls +(2,0) and +(0,1) ..
(14,-6) .. controls +(0,-1) and +(2,0) .. (11,-9) .. controls
+(-2,0) and +(0,-1) .. (8,-6) .. controls +(0,1) and +(-2,0) ..
(11,-3);

\draw (11,-4) -- (12.5,-4) -- (13,-2.5) -- (11,-2) -- (11,-4) --
(9.5,-4) -- (9,-2.5) -- (11,-2);

\draw [dashed] (12.5,-4) -- (13.2,-5.3);

\draw [dashed] (9.5,-4) -- (8.8,-5.3);

\draw [dashed] (13,-2.5) -- (14,-3.5);

\draw [dashed] (9,-2.5) -- (8,-3.5);

\node at (10,-2.8) {$\times$};

\node at (12,-2.8) {$\star$};

\node at (16.5,-7) {$N$};

\draw [->,dotted] (13,-2.5) .. controls +(.1,.5) and +(-.5,0) ..
(14.5,-1) .. controls +(1,0) and +(0,1) .. (16,-4) .. controls
+(0,-1) and +(1,0) .. (14,-7) .. controls +(-1,0) and +(.1,-1) ..
(12.5,-4);

\node at (16,-1) {$\cercle$};

\node at (7,-6) {$\widehat{T}_{1}$};

\draw [->] (8,-8) -- (3,-12);

\draw [very thick] (8,15) .. controls +(2,0) and +(0,1) .. (13,12)
.. controls +(0,-1) and +(2,0) .. (8,9) .. controls +(-2,0) and
+(0,-1) .. (3,12) .. controls +(0,1) and +(-2,0) .. (8,15);

\draw [->,dotted] (8,15) .. controls +(0,1) and +(-1,0) .. (9,18)
.. controls +(1,0) and +(0,1) .. (10,15) .. controls +(0,-1) and
+(1,0) .. (9,12) .. controls +(-1,0) and +(0,-1) .. (8,15);

\node at (10.5,18) {$\cercle$};

\node at (2.3,12) {$F$};

\node at (13,15) {$W$};

\draw [->] (3.5,10.5) -- (0.5,7.5);

\draw [->] (11,9) -- (11,-1);

\end{tikzpicture}
\end{center}

~

The two covers $M'$ and $N$ admit a common regular finite cover
$W$, which fibers over the circle as it is a finite cover of $N$.
A component of the preimage of $\widehat{T_1}$ by the covering
projection $W \rightarrow N$ is a fiber $F$ for the fibration of
$W$ over the circle. As the diagram is commutative, it is also a
component the preimage of the embedded surface $T_1$ in $M'$, as
$T_1$ and $\widehat{T_1}$ have the same image in $M$, which is an
immersed surface. As $F$ is incompressible in $W$, the surface
$T_1$ embedded in $M'$ we started from is also incompressible.
 Thus the embedded surface $T_1$ is a virtual fiber
in $M'$, and is incompressible.\qed

~

Therefore, the $m$ initial parallel surfaces are virtual fibers
for the manifold $M'$. In fact, they are fibers of a bundle over the circle or
of a twisted $I$-bundle. Indeed, if $T$ is one of those surfaces, the complement
$M'_T$ of an open neighborhood of $T$ in $M'$ admits a finite cover that is the
product of a $T'$ by an interval $I$. In particular, the fundamental group of
the compact manifold $M'_T$ contains a finite index surface subgroup. By
\cite[Theorem 10.6]{Hem}, it is a $I$-bundle, possibly twisted. This ends the
proof of the Pattern Proposition D.\qed

\section{Heegaard genera and fibration.}\label{section-grad-s-log}

The proof of theorem \ref{thm-gHetfibration} is the starting point
for the proof of the main theorem A. The aim was to establish a
virtual fibration criterion standing between Lackenby's conjecture
\ref{cj-GH} and Maher's theorem. Maher himself suggested in
\cite{Mah} the possibility to get explicit constants and upper
bounds at each stage of the proof of Theorem 1.1 of \cite{Mah},
but without precise statements.

This section is dedicated to the proof of theorem
\ref{thm-gHetfibration} and corollary \ref{cor-gHslog}.

\subsection{Proof of theorem \ref{thm-gHetfibration} (1) and corollary
\ref{cor-gHslog} (1): Heegaard genus.}\label{dem-thm-glfaible}

~

\noindent \underline{Proof of theorem \ref{thm-gHetfibration}
(1).} \nopagebreak

Suppose that $M' \rightarrow M$ is a cover of $M$ with finite
degree $d$. Let $S\subset M'$ be a minimal genus Heegaard surface for
$M'$: $g(S)=g(M')$. The aim is to construct from $S$ a
pseudo-minimal surface which satisfies assumptions of theorem A.
We start with the following lemma.

\begin{lm}\label{lm-SHetSurfPM}

Let $N$ be a connected, oriented and closed hyperbolic 3-manifold.
Let $S$ be a minimal genus Heegaard surface for $N$ and
$\mathcal{H}$ the Heegaard splitting of $N$ with Heegaard surface
$S$. Let $F$ be the union of the even and odd surfaces of a
$\mathcal{H}$-thin generalized Heegaard splitting for $N$. Then
$F$ is a pseudo-minimal surface, which divides the manifold $N$ in
$q \leq\chi_-^h(N)+2$ compression bodies $C_1,\ldots , C_q$ with
$\chi_-(C_j)\leq \chi_-^h(N)$ for all $j$ between $1$ and $q$.

Furthermore, if $F^-$ is the union of the negative boundary components
$\partial_-C_j$, then it is a union of incompressible surfaces.

\end{lm}

\noindent \underline{Proof of lemma \ref{lm-SHetSurfPM}.}
\nopagebreak

The topological part (1) of the following theorem is a consequence
of works of Casson and Gordon, Scharlemann and Thompson
(\cite{CaGo} and \cite{ST2}). The metric part (2) comes from
results of Frohman, Freedman, Hass and Scott about incompressible
surfaces (\cite{FHS} and \cite{FrHa}). The last part (3) is a
result of Pitts and Rubinstein (\cite{PiRu}, see also \cite{Sou},
\cite{CoDL} and \cite{Pi}).

\begin{thm}\label{thm-surf-min}

Let $N$ be a connected, oriented and closed hyperbolic 3-manifold,
and $\mathcal{H}$ a $\mathcal{H'}$-thin generalized Heegaard splitting for some
Heegaard decomposition $\mathcal{H}'$. Then
$\mathcal{H}$ satisfies the following properties.

\begin{enumerate}
  \item Each of the even surfaces is incompressible in $N$ and the odd
surface are strongly irreducible Heegaard surfaces for the
components of the manifold $N$ cut along the even surfaces.
  \item Each even surface can be isotoped to a minimal surface or the
boundary of a small neighborhood of a non-orientable minimal
surface.
  \item each odd surface can be isotoped to a minimal surface,
or to the boundary of a small regular neighborhood of a
non-orientable minimal surface, with a small tube attached
vertically in the $I$-bundle structure.
\end{enumerate}

\end{thm}\qed

Thanks to theorem \ref{thm-surf-min}, up to isotopy, one can
assume that the surface $F$ is pseudo-minimal, and it is immediate that $F^-$
is a union of incompressible surfaces.

As described in section \ref{rappels}, surgeries of generalized
Heegaard splittings are a modification in the order of attachment
of the 1- and 2-handles of a corresponding handle decomposition of
the manifold. Therefore, surgeries do not change the number of 1-
and 2-handles. As it is equal in the starting Heegaard splitting
to $2g(S) = \chi_-^h(N)+2$, there are also $(\chi_-^h(N)+2)$ 1-
and 2-handles in a handle decomposition associated to the surface
$F$. As this number is an upper bound for the number of
compression bodies in the complement of $F$, the inequality $q\leq
\chi_-^h(N) +2$ holds.

Furthermore, as each component of $F$ is obtained from $S$ by
surgery, the characteristic $\chi_-(C)$ of each compression body
is at most $\left|\chi(S)\right| = \chi_-^h(N)$.\qed

~

\noindent \underline{End of the proof of theorem
\ref{thm-gHetfibration} (1).} \nopagebreak

Recall that $S$ is a minimal genus Heegaard surface for the cover
$M'\rightarrow M$ of finite degree $d$. Let $F$ be the
pseudo-minimal surface obtained in lemma \ref{lm-SHetSurfPM}. The
aim is to apply the main theorem A to $F$. With notations of
theorem A and this choice of surface $F$, one has $c=\chi_-^h(M')$
and $q=\chi_-^h(M')+2$.

Set $\epsilon = \inj(M)/2$ and let $k=k(\epsilon, \vol(M))$ be the
constant obtained in theorem A. To satisfy assumptions of theorem
A, one needs to have
$$k\, \chi_-^h(M') \ln \chi_-^h(M') \leq \ln
\ln \frac{d}{\chi_-^h(M')+2}.$$

If the ratio $\chi_-^h(M')\ln \chi_-^h(M') /\ln \ln d$ tends to
zero, then the ratio $\chi_-^h(M') /\sqrt{d}$ tends also to zero.
Therefore, there exists an explicit constant $\overline{k_1}>0$
such that if $\overline{k_1}\, \chi_-^h(M') \ln \chi_-^h(M') \leq
\ln \ln d$, then $\chi_-^h(M')+2 \leq \sqrt{d}$. Under this
assumption, one has
\begin{eqnarray*}
\ln \ln \frac{d}{\chi_-^h(M')+2} & \geq & \ln \ln \sqrt{d}\\
& =& \ln \left(\frac{1}{2}\ln d\right)\\
&=& \ln \ln d - \ln 2 \geq \frac{1}{2} \ln \ln d
\end{eqnarray*}
if $\ln \ln d\geq 2\ln 2$, which is the fact for example if $\ln
\ln d \geq \chi_-^h(M') \ln \chi_-^h(M')$ as $\chi_-^h(M') \geq
2$.

Therefore, if $\chi_-^h(M') \ln \chi_-^h(M') \leq
\ln \ln d$, $\overline{k_1}\, \chi_-^h(M') \ln \chi_-^h(M') \leq
\ln \ln d$ and\\ $2k\, \chi_-^h(M') \ln \chi_-^h(M') \leq \ln \ln
d$, then
$$k\, \chi_-^h(M') \ln \chi_-^h(M') \leq \ln
\ln \frac{d}{\chi_-^h(M')+2}$$ and assumptions of theorem A are
satisfied. This proves theorem \ref{thm-gHetfibration} with
$\overline{k}=\max\{1,2k, \overline{k_1}\}$.\qed

~

\noindent \underline{Proof of corollary \ref{cor-gHslog} (1).}
\nopagebreak

It is obvious that if $M$ virtually fibers over the circle, then
the $\eta$-sub-logarithmic Heegaard gradient of $M$ is zero for
every $\eta \in (0,1)$, as $M$ admits an infinite family of finite
degree covers with bounded Heegaard genus.

If the $\eta$-sub-logarithmic Heegaard gradient of $M$ is zero for
some $\eta \in (0,1)$, this means that $M$ admits an infinite
family of covers $(M_i\rightarrow M)_{i\in \nn}$ with finite
degrees $d_i$, and such that
$$\lim_{i\div} \frac{\chi_-^h(M_i)}{(\ln \ln d_i)^\eta} =0,$$
which can be also written
$$\lim_{i\div} \frac{\chi_-^h(M_i)^{1/\eta}}{\ln \ln d_i} =0.$$ As
$1/\eta >1$, this implies that
$$\lim_{i\div} \frac{\chi_-^h(M_i) \ln \chi_-^h(M_i)}{\ln \ln d_i} =0.$$

Thus, for $i$ large enough, on has $\overline{k}\, \chi_-^h(M_i)
\ln \chi_-^h(M_i) \leq \ln \ln d_i$ and the assumptions of theorem
\ref{thm-gHetfibration} are satisfied. In particular, the manifold
$M$ virtually fibers over the circle, which proves corollary
\ref{cor-gHslog} (1).\qed

\subsection{Proof of theorem \ref{thm-gHetfibration} (2) and
corollary \ref{cor-gHslog} (2): strong Heegaard
genus.}\label{dem-thm-glfort}

~

\noindent \underline{Proof of theorem \ref{thm-gHetfibration}
(2).}

Suppose by contradiction that in a finite cover $M'\rightarrow M$
of degree $d$, one has $\bar{k}\, \chi_-^{sh}(M') \ln
\chi_-^{sh}(M') \leq \ln \ln d$. Let $F$ be a strongly irreducible
Heegaard surface for $M'$ such that
$\chi_-^{sh}(M')=\chi_-^{}(F)$.

Thanks to theorem \ref{thm-surf-min}, up to isotopy, one can
assume that $F$ is pseudo-minimal. This surface separates $M'$
into two handlebodies, so the volume of one of those handlebodies
$C$ must be at least $\vol(M) d/2$. But as $\bar{k}\, \chi_-^{}(F)
\ln \chi_-^{}(F) \leq \ln \ln d$, the proof of theorem
\ref{thm-gHetfibration} (1) shows that the surface $F$ satisfies
the assumptions of theorem A. This is in contradiction with
corollary \ref{cor-thmA-volhandlebodies}. This proves theorem
\ref{thm-gHetfibration} (2).\qed

~

\noindent \underline{Proof of corollary \ref{cor-gHslog} (2).}

To prove corollary \ref{cor-gHslog} (2), just notice that as
$\bar{k}\, \chi_-^{sh}(M') \ln \chi_-^{sh}(M') > \ln \ln d$, for
$\theta \in (0,1)$, there is a constant $\bar{k}_\theta>0$ such
that $\chi_-^{sh}(M')/(\ln \ln d)^\theta \geq \bar{k}_\theta$,
proving that the strong $\eta$-sub-logarithmic Heegaard gradient
of $M$ is strictly positive.\qed

\section{Circular decomposition and fibered homology classes.}

The aim of this section is to consider the case of circular
decompositions, and to prove corollaries \ref{cor-dec-circ} and
\ref{cor-surf-incomp&fibre-virtuelle}.

\subsection{Circular decomposition and thin position.}

~

A circular decomposition is the equivalent of a Heegaard
decomposition, but this decomposition is associated to a Morse
function that no longer takes values in $I=[0,1]$ but in the
circle $\cercle$.

\begin{df}

A \textbf{circular Morse function} is a Morse function $f\,:\, M
\rightarrow \cercle$.

If $f\,:\, M \rightarrow \cercle$ is a circular Morse function,
the handle decomposition of $M$ given by the function $f$ is
called the \textbf{circular decomposition associated to $f$}.

\end{df}

See F. Manjarrez-Guti\'errez \cite{MG}, Matsumoto \cite{Mat} and
Milnor \cite{Mil} for further details about circular Morse
functions. Let $f\,:\,M\rightarrow \cercle$ be a circular Morse
function. If we remove a small open neighborhood of a regular
value $x \in \cercle$, by restriction of $f$, we obtain a Morse
function $g$ of $M_R=M\setminus \mathcal{N}(R)$, which is the
manifold $M$ minus a small regular open neighborhood of the
surface $R:=f^{-1}(\{x\})$, on the interval $I$. Thus, the theory
of Heegaard splittings and generalized Heegaard splittings applies
to the function $g$, as recalled in section \ref{rappels}.

~

An other viewpoint is to see a circular decomposition as a handle
decomposition of the cobordism $(M\setminus\mathcal{N}(R),
R\times\{1\}, R\times \{-1\})$. Starting with a Heegaard splitting
of Heegaard surface $S$ for $M_R=M\setminus\mathcal{N}(R)$, one
can change the order in which $1$- and $2$-handles are attached to
get a new generalized Heegaard splitting
$(F_1=R\times\{1\},S_1,F_2,\ldots,S_n,F_{n+1}=R\times\{-1\})$ for
$(M_R, R\times\{1\}, R\times \{-1\})$. Gluing back $R\times\{1\}$
to $R\times\{-1\}$, one obtains a circular decomposition for the
manifold $M$. Denote it by
$\mathcal{H}=(F_1,S_1,F_2,\ldots,S_n,F_{n+1})$, with
$F_1=F_{n+1}=R$. The surfaces $F_j$ divide $M$ into $n$
3-manifolds with boundary $W_1,\ldots,W_n$, and surfaces $S_j$ are
Heegaard surfaces for those manifolds. For $1\leq j\leq n$, $S_j$
divides the manifold $W_j$ into two compression bodies $A_j$ and
$B_j$, such that $\partial_+A_j=\partial_+B_j=S_j$,
$\partial_-A_j=F_j$ and $\partial_-B_j=F_{j+1}$.

Let $S$ be a closed surface. If $S$ is connected, recall that the
\textbf{complexity} of $S$ is $c(S)=\max(0,2g(S)-1)$. If $S$ is
the union of several connected components, the complexity of $S$
is the sum of the complexities of the components of $S$. There is
a definition of the complexity of a circular decomposition
analogous to the complexity of a generalized Heegaard splitting.

\begin{df}

The \textbf{circular width} of a circular decomposition
$\mathcal{H}=(F_1,S_1,F_2,\\ \ldots,S_n,F_{n+1})$ is the set of
the $n$ integers $(c(S_1),\ldots,c(S_n))$, with repetitions and
arranged in monotonically non-increasing order. Widths are
compared using the lexicographic order.

The integer $n\geq 1$ is called the \textbf{length} of the
circular decomposition
$\mathcal{H}=(F_1,S_1,\\F_2,\ldots,S_n,F_{n+1})$.

\end{df}

\begin{prop}\label{prop-dec-mince}

Let $M$ be a hyperbolic, connected, oriented and closed
3-manifold. Let $R$ be an orientable, closed, non-separating,
incompressible and embedded surface in $M$. Denote by $S$ a
Heegaard surface for $M\setminus \mathcal{N}(R)$. Starting from
the circular decomposition $\mathcal{H}=(R,S,R)$ of $M$, there
exists a finite number of surgeries to get a circular
decomposition $\mathcal{H'}=(F_1,S_1,F_2,\ldots,S_n,F_{n+1})$ with
$F_1=F_{n+1}=R$, such that:
\begin{enumerate}
  \item the circular width of $\mathcal{H'}$ is minimal among the
  widths of such circular decompositions obtained by a finite number
  of surgeries of $\mathcal{H}$,
  \item each surface $F_j$ is incompressible, no component of $F_j$
  is a sphere, and $g(F_j)\leq g(S)$,
  \item each surface $S_j$ is a strongly irreducible Heegaard
  surface for the Heegaard decomposition $(A_j,B_j)$ of $W_j$
  and $g(S_j)\leq g(S)$,
  \item $n\leq \frac{1}{2}(\chi(R) - \chi(S))$,
\end{enumerate}

\end{prop}

\begin{df}

Let $\mathcal{H}$ be a circular decomposition. A circular
decomposition $\mathcal{H'}=(F_1,S_1,F_2,\ldots,S_n,F_{n+1})$ that
is circular-length-minimizing among all circular decompositions
obtained from $\mathcal{H}$ by a finite number of surgeries is
said to be a \textbf{thin position}. We will call such a
decomposition a \textbf{thin circular decomposition associated to
$\mathcal{H}$.}

\end{df}

\noindent \underline{Proof of proposition \ref{prop-dec-mince}.}
\nopagebreak

The proof of this proposition is essentially the same as the proof
of \cite[Theorem 3.2]{MG}, which is itself an adaptation of
techniques of \cite{ST2} to the case of circular decompositions
(see also \cite{Lac}). See \cite{Ren2}, Proposition 1.1 and its
proof. The proof is based on an operation called a
\textbf{surgery} of circular decompositions, which is analogous to
the surgery of generalized Heegaard splittings described in
section \ref{rappels}. Again, the crucial fact is that a surgery
procedure strictly decreases the complexity of the circular
decomposition.\qed

\begin{cor}\label{cor-dec-mince&surfPM}

Let $M$ be a hyperbolic, connected, oriented and closed
3-manifold. Take $\mathcal{H}=(F_1,S_1,F_2,\ldots,S_n,F_{n+1})$ a
thin circular decomposition of $M$. Then, up to isotopy, one can
assume that all surfaces $F_j$ and $S_j$ are pseudo-minimal.

\end{cor}

\noindent \underline{Proof of corollary \ref{cor-dec-mince&surfPM}.}
\nopagebreak

From proposition \ref{prop-dec-mince} points (2) and (3), the
surfaces $F_j$ are incompressible for each $j$ and the surfaces
$S_j$ correspond to strongly irreducible Heegaard surfaces. The
proof is then the same as for theorem \ref{thm-surf-min} (2) and
(3) of section \ref{section-grad-s-log}.\qed

\subsection{Circular characteristic and fibered homology
classes.}\label{grad-cir&homol-fibree}

~

Recall definition \ref{def-deccirc} from the introduction.
Corollary \ref{cor-dec-circ} is analogous to theorem
\ref{thm-gHetfibration} for circular decompositions associated to
a non trivial cohomology class.

~

\noindent \underline{Proof of corollary \ref{cor-dec-circ}.}
\nopagebreak

Let $M'\rightarrow M$ be a cover of $M$ with finite degree $d$,
and a non-trivial cohomology class $\alpha'\in H^1(M')$. The aim
is to show that if the ratio $\chi_-^c(\alpha') \ln
\chi_-^c(\alpha') / \ln\ln d$ is small enough, then the
assumptions of theorem A are satisfied.

Let $R'$ be an embedded surface in $M'$ and
$\|\alpha'\|$-minimizing. First, suppose that in addition
$h(M',\alpha')=h(M',\alpha',R')$. Take $S'$ a minimal genus
Heegaard surface for $M'_{R'}$. By construction,
$\chi_-^c(\alpha')=\left|\chi(S')\right|$.

From proposition \ref{prop-dec-mince}, starting from the circular
decomposition $(R',S',R')$ of $M'$, we can construct a thin
circular decomposition
$\mathcal{H}=(F_1,S_1,F_2,\ldots,S_n,F_{n+1})$. Set
$F:=\bigcup_jF_j\cup \bigcup_jS_j$. From corollary
\ref{cor-dec-mince&surfPM}, one can assume that $F$ is a
pseudo-minimal surface.

Still from proposition \ref{prop-dec-mince}, the surface $F$
separates the manifold $M'$ into $q\leq \frac{1}{2} (\chi(R') -
\chi(S')) \leq \chi_-^c(\alpha')/2$ compression bodies
$C_1,\ldots, C_q$, with $\chi_-(C_j)\leq
\left|\chi(S')\right|=\chi_-^c(\alpha')$ for every $j$.

As the surfaces $F_j$ are incompressible, assumption $(1)$ of theorem A is
satisfied.

Let $k=k(\epsilon,\vol(M))$ be the constant given by theorem A. To
satisfy assumptions of theorem A, there remains to show that $k\,
\chi_-^c(\alpha') \ln \chi_-^c(\alpha') \leq \ln \ln
\frac{2d}{\chi_-^c(\alpha')}$. But as in section
\ref{section-grad-s-log}, one can find a constant
$\ell'=\ell'(\epsilon,\vol(M))$ such that if $\ell'
\,\chi_-^c(\alpha') \ln \chi_-^c(\alpha') \leq \ln \ln  d$, then
$k\, \chi_-^c(\alpha') \ln \chi_-^c(\alpha') \leq \ln \ln
\frac{2d}{\chi_-^c(\alpha')}$ and all the assumptions of theorem A
are satisfied.

Therefore, if $\ell' \,\chi_-^c(\alpha') \ln \chi_-^c(\alpha')
\leq \ln \ln  d$, then from theorem A, the manifold $M'$ contains
an embedded surface that is a virtual fiber.

Furthermore, all the constructions take place in fact in
${M'}_{R'}=M'\setminus \mathcal{N}(R')$. Thus, the virtual fiber
built in theorem A is in the complement of $R'$ in $M'$. This
virtual fiber lifts to a connected fiber $\overline{T}$ in a
fibered finite cover $\overline{M'}\rightarrow M'$ of $M'$. In
this cover, the incompressible surface $R'$ lifts to a surface
$\overline{R'}$ in the complement of the fiber. Cutting along
$\overline{T}$, this shows that the connected components of
$\overline{R'}$, which are all incompressible, are parallel in the
product to the fiber $\overline{T}$ (see \cite{Wal}). Thus, the
homology class of $\overline{R'}$ is fibered. Still from Gabai
\cite[Lemma 2.4]{Ga1}, this implies that the homology class of
$R'$ is fibered. As the surface $R'$ minimizes Thurston's norm, it
is also a fiber.

To end the proof of corollary \ref{cor-dec-circ}, there remains to
show that if  $R'$ be an embedded surface in $M'$,
$\|\alpha'\|$-minimizing, but that does not necessarily satisfy
$h(M',\alpha')=h(M',\alpha',R')$, then $R'$ is still a fiber. But
if one takes an embedded surface $R''$ such that $R''$ is
$\|\alpha'\|$-minimizing, and satisfies
$h(M',\alpha')=h(M',\alpha',R'')$, the proof above shows that
$R''$ is a fiber. As $R'$ is norm-minimizing, it is an
incompressible surface in the homology class of $R''$, hence also
a fiber. This ends the proof of corollary \ref{cor-dec-circ}.\qed

~

The following corollary is immediate from corollary
\ref{cor-dec-circ}.

\begin{cor}\label{cor-complexitecirculaire-asymptotique}

Let $M$ be a hyperbolic, connected, oriented and closed
3-manifold. Suppose that there exists an infinite family of covers
$(M_i\rightarrow M)_{i\in\nn}$ with finite degrees $d_i$, and for
each $i\in\nn$, a non-trivial cohomology class $\alpha_i\in
H^1(M_i)$ such that:
$$\inf_{i\in\nn}\frac{\chi_-^c(\alpha_i)\ln \chi_-^c(\alpha_i)}{\ln\ln d_i}=0.$$
Then, for infinitely many indices $i\in\nn$, every embedded
surface $R_i$ in $M_i$, $\|\alpha_i\|$-minimizing and such that
$h(M_i,\alpha_i) = h(M_i,\alpha_i,R_i)$ is a fiber. In particular,
the manifold $M$ virtually fibers over the circle $\cercle$.

\end{cor}

\subsection{Incompressible surfaces and fibrations.}

~

In section \ref{grad-cir&homol-fibree}, we have established
criteria in order to show that a non-trivial cohomology class of a
hyperbolic 3-manifold $M$ lifts to fibered classes in finite
covers. Now, if $R$ is a non separating embedded surface in $M$,
there is a dual cohomology class associated to $R$. In some cases,
we have seen that $R$ could then be a fiber. But the
question can be asked for any embedded, incompressible and
connected surface $R$ in $M$, separating or not.

Recall definition \ref{def-caracteristique-surface} from the
introduction. Corollary \ref{cor-surf-incomp&fibre-virtuelle} is
different from last section as the surface $R$ is a priori not
supposed to be non-separating.

~

\noindent \underline{Proof of corollary \ref{cor-surf-incomp&fibre-virtuelle}.}
\nopagebreak

In the case where the surface $R'$ is not separating, it is a
generalization of corollary \ref{cor-dec-circ}. Indeed, if $S'$ is
a minimal genus Heegaard surface for $M'_{R'}$, $\chi_-^h(R') =
\left|\chi(S')\right|$ and $(R',S',R')$ is a circular
decomposition of $M'$. As the starting surface $R'$ is
incompressible, we apply then proposition \ref{prop-dec-mince} to
build a thin circular decomposition. From the proof of corollary
\ref{cor-dec-circ}, assumptions of theorem A are satisfied if
$\ell' \chi_-^h(R') \ln \chi_-^h(R') \leq \ln \ln d$, and in this
case, the surface $R'$ is a virtual fiber. But as $R'$ belongs to
the preimage of $R$, the surface $R$ is also a virtual fiber.
Furthermore, if the surface $R$ is not separating, its homology
class is non zero, and the same argument applies to prove that
this class is fibered. As the surface $R$ is incompressible, it is
itself a fiber.

In the case where the surface $R'$ separates the manifold $M'$
into two connected components $M_l$ and $M_r$, let $S_l$ and $S_r$
be minimal genus Heegaard surfaces for $M_l$ and $M_r$
respectively. By definition, $\chi_-^h(R') =
\max\{\left|\chi(S_l)\right|, \left|\chi(S_r)\right|\}$.

In each side of $R'$, we can then build a generalized Heegaard
decomposition for $M_l$ and $M_r$ in thin position starting from
surfaces $S_l$ and $S_r$. We get then a surface $F$ with one
component which is the incompressible surface $R'$, separating the
manifold $M'$ into $q\leq 2g(S_l) +2g(S_r) \leq 2\chi_-^h(R') +4$
compression bodies $C_1,\ldots ,C_q$ with $\chi_-(C_j) \leq
\chi_-^h(R')$ for all $j$.

As $F$ is the union of incompressible surfaces and strictly
irreducible Heegaard surfaces, we may assume that $F$ is
pseudo-minimal, and $F^-$ is the union of incompressible surfaces.

Thus, to satisfy assumptions of theorem A, it suffices to show
that $k\, \chi_-^h(R')\\ \ln \chi_-^h(R') \leq \ln \ln
\frac{d}{2\chi_-^h(R') +4}$. But as before, one can find a
constant $\ell''\geq \ell'$ such that if $\ell''\,\chi_-^h(R') \ln
\chi_-^h(R') \leq \ln \ln d$, then $k\, \chi_-^h(R') \ln
\chi_-^h(R') \leq \ln \ln \frac{d}{2\chi_-^h(R') +4}$. From
theorem A, in this case, the manifold $M'$ contains an
incompressible surface that is a virtual fiber. But as in the
proof of corollary \ref{cor-dec-circ}, this incompressible surface
is built in the complement of the incompressible surface $R'$.
Thus, the surface $R'$ is also a virtual fiber. As $R'$ is a lift
of $R$, the starting surface $R$ is a virtual fiber, hence the fiber of a
twisted $I$-bundle, which ends
the proof of corollary \ref{cor-surf-incomp&fibre-virtuelle}.\qed

\bibliographystyle{short}

\bibliography{biblio}

~

Claire \textsc{Renard},

\'Ecole Normale Sup\'erieure de Cachan,

Centre de Math\'ematiques et de Leurs Applications.

61 avenue du pr\'esident Wilson

F-94235 CACHAN CEDEX.

\emph{claire.renard@normalesup.org}

\end{document}